\documentclass[reqno,11pt]{article}
\usepackage{psfig, amsmath, amstexnb, amsthm}
\usepackage{amssymb}  

\newtheorem{theorem}{Theorem}[section]
\newtheorem{prop}[theorem]{Proposition}
\newtheorem{lemma}[theorem]{Lemma}

\newtheorem{notation}[theorem]{Notation}

\theoremstyle{definition}
\newtheorem{remark}[theorem]{Remark}

\numberwithin{equation}{section}

\makeatletter

\def\enum{\ifnum \@enumdepth >3 \@toodeep\else
        \advance\@enumdepth \@ne 
        \edef\@enumctr{enum\romannumeral\the\@enumdepth}\list
        {\csname label\@enumctr\endcsname}
        {\setlength{\topsep}{1mm}
        \setlength{\parsep}{0mm}
        \setlength{\itemsep}{0mm}
        \setlength{\labelsep}{2mm}
        \settowidth{\leftmargin}{M.}
        \addtolength{\leftmargin}{\labelsep}
        \usecounter{\@enumctr}
        \def\makelabel##1{\hss\llap{##1}}}\fi}

\def\itemiz{\ifnum \@itemdepth >3 \@toodeep\else \advance\@itemdepth \@ne
        \edef\@itemitem{labelitem\romannumeral\the\@itemdepth}%
        \list{\csname\@itemitem\endcsname}{
        \setlength{\topsep}{1mm}
        \setlength{\parsep}{0mm}
        \setlength{\parsep}{0mm}
        \setlength{\itemsep}{0mm}
        \setlength{\labelsep}{2mm}
        \settowidth{\leftmargin}{M.}
        \addtolength{\leftmargin}{\labelsep}
        \def\makelabel##1{\hss\llap{##1}}}\fi}

\def\captionheadfont@{\scshape}
\def\captionfont@{\small}
\long\def\@makecaption#1#2{%
  \setbox\@tempboxa\vbox{\color@setgroup
    \advance\hsize-3pc\noindent
    \captionfont@\captionheadfont@#1\@xp\@ifnotempty\@xp
        {\@cdr#2\@nil}{.\captionfont@\upshape\enspace#2}%
    \unskip\kern-3pc\par
    \global\setbox\@ne\lastbox\color@endgroup}%
  \ifhbox\@ne 
    \setbox\@ne\hbox{\unhbox\@ne\unskip\unskip\unpenalty\unkern}%
  \fi
  \ifdim\wd\@tempboxa=\z@ 
    \setbox\@ne\hbox to\columnwidth{\hss\kern-3pc\box\@ne\hss}%
  \else 
    \setbox\@ne\vbox{\unvbox\@tempboxa\parskip\z@skip
        \noindent\unhbox\@ne\advance\hsize-3pc\par}%
\fi
  \ifnum\@tempcnta<64 
    \addvspace\abovecaptionskip
    \moveright 1.5pc\box\@ne
  \else 
    \moveright 1.5pc\box\@ne
    \nobreak
    \vskip\belowcaptionskip
  \fi
\relax
}

\def\overbar#1{\skewbar{#1}{-1}{-1}{.25}}
\def\skewbar#1#2#3#4{\preciseskewbar{#1}{#2}{#3}{#2}{#3}{#2}{#3}{#4}1}
\def\preciseskewbar#1#2#3#4#5#6#7#8#9{{\mathchoice
 {\makeoverbar\textfont\displaystyle{#1}1{#2}{#3}{#8}{#9}}
 {\makeoverbar\textfont\textstyle{#1}1{#2}{#3}{#8}{#9}}
 {\makeoverbar\scriptfont\scriptstyle{#1}{.7}{#4}{#5}{#8}{#9}}
 {\makeoverbar\scriptscriptfont
  \scriptscriptstyle{#1}{.5}{#6}{#7}{#8}{#9}}}#1}
\def\makeoverbar#1#2#3#4#5#6#7#8{{%
 \setbox0=\hbox{$\m@th#2\mkern#5mu{#3}\mkern#6mu$}%
 \setbox1=\null \dimen@=#4\fontdimen8#13 \dimen@=3\dimen@ 
 \advance\dimen@ by \ht0 \dimen@=-#7\dimen@ \advance\dimen@ by \wd0
 \wd1=\dimen@ \dp1=\dp0 
 \dimen@=#4\fontdimen8#13
 \dimen@i=\fontdimen8#13
 \fontdimen8#13=#8\dimen@
 \advance\dimen@ by -\fontdimen8#13 \dimen@=3\dimen@
 \advance\dimen@ by \ht0 \ht1=\dimen@ 
 \rlap{\hbox to \wd0{$\m@th\hss#2{\overline{\box1}}\mkern#5mu$}}
 \fontdimen8#13=\dimen@i}}

\makeatother

\DeclareMathSymbol{\leqsymb}{\mathalpha}{AMSa}{"36}
\def\leqs{\mathrel\leqsymb}
\DeclareMathSymbol{\geqsymb}{\mathalpha}{AMSa}{"3E}
\def\geqs{\mathrel\geqsymb}
\DeclareMathSymbol{\gtreqqlesssymb}{\mathalpha}{AMSa}{"54}

\def\,{\ifmmode\mskip\thinmuskip\else\kern\fontdimen3}

\def\figref#1{Fig.$\mskip-1mu$~\ref{#1}}

\newcommand{\field}[1]{\mathbb{#1}}

\newcommand{\R}{\field{R}\,}    
\newcommand{\E}{\field{E}}      
\newcommand{\fP}{\field{P}}     

\newcommand{\cA}{{\mathcal A}}  
\newcommand{\cB}{{\mathcal B}}  
\newcommand{\cC}{{\mathcal C}}  
\newcommand{\cD}{{\mathcal D}}  
\newcommand{\cF}{{\mathcal F}}  
\newcommand{\cM}{{\mathcal M}}  
\newcommand{\cS}{{\mathcal S}}  

\DeclareMathOperator{\e}{e}             
\DeclareMathOperator{\dd}{d}            

\def\math#1{\ifmmode
\mathchoice{\mbox{$\displaystyle\rm#1$}}
{\mbox{$\textstyle\rm#1$}}
{\mbox{$\scriptstyle\rm#1$}}
{\mbox{$\scriptscriptstyle\rm#1$}}\else
{\mbox{$\rm#1$}}\fi}            
\DeclareMathOperator{\defby}{\raisebox{0.35pt}{\math{:}}\!\!=}

\def\nbh{neighbourhood}

\def\eps{\varepsilon}
\def\w{\omega}

\def\z{\zeta}
\def\t{t}
\def\xdet{x^{\det}}
\def\xhatdet{\widehat x^{\mskip2mu\det}}
\def\ba{\bar a}
\def\bb{\bar b}

\def\balpha{\skewbar{\alpha}{.5}{-1}{.25}}
\def\bx{\bar x}
\def\ta{\tilde a}
\def\ha{\widehat a}
\def\hb{\widehat b}
\def\tx{\tilde x}
\def\talpha{\tilde\alpha}
\def\halpha{\widehat\alpha}

\def\xc{x_{\math c}}
\def\tc{t_{\math c}}
\def\lc{\lambda_{\math c}}

\def\subas{\stackrel{{\rm a.s.}}{\subset}}      
\def\indexfct#1{1_{\set{#1}}}           

\def\6#1{\dd\!#1}                       
\def\dpar#1#2{\frac{\partial #1}{\partial #2}}  
\def\sdpar#1#2{\partial_{#2}#1}         
\def\dtot#1#2{\frac{\6{#1}}{\6{#2}}}  

\def\bigpar#1{\bigl(#1\bigr)}                   
\def\biggpar#1{\biggl(#1\biggr)}        
\def\Bigpar#1{\Bigl(#1\Bigr)} 
\def\Biggpar#1{\Biggl(#1\Biggr)} 

\def\bigbrak#1{\bigl[#1\bigr]}          
\def\Bigbrak#1{\Bigl[#1\Bigr]}          
\def\biggbrak#1{\biggl[#1\biggr]}       
\def\Biggbrak#1{\Biggl[#1\Biggr]}       

\def\set#1{\{#1\}}                              
         
\def\Bigset#1{\Bigl\{#1\Bigr\}}         
\def\biggset#1{\biggl\{#1\biggr\}} 
\def\Biggset#1{\Biggl\{#1\Biggr\}} 

\def\setsuch#1#2{\{#1\colon #2\}}                
\def\bigsetsuch#1#2{\bigl\{#1\colon #2\bigr\}}  
\def\Bigsetsuch#1#2{\Bigl\{#1\colon #2\Bigr\}}  
\def\biggsetsuch#1#2{\biggl\{#1\colon #2\biggr\}}       

\def\abs#1{\lvert#1\rvert}                      

\def\intpartplus#1{\lceil#1\rceil}              

\def\biggintpartplus#1{\biggl\lceil#1\biggr\rceil}

\def\Order#1{{\mathcal O}(#1)}                  
\def\bigOrder#1{{\mathcal O}\bigl(#1\bigr)}     
\def\BigOrder#1{{\mathcal O}\Bigl(#1\Bigr)}     
\def\biggOrder#1{{\mathcal O}\biggl(#1\biggr)}  

\def\orderone#1{{\scriptstyle\mathcal O}_{#1}(1)}      

\def\bigprobin#1#2{\fP^{\mskip1.5mu #1}\bigl\{#2\bigr\}}    
\def\Bigprobin#1#2{\fP^{\mskip1.5mu #1}\Bigl\{#2\Bigr\}}

\def\Bigexpecin#1#2{\E^{\mskip1.5mu #1}\Bigl\{#2\Bigr\}}

\def\bibtitle#1#2{#1, {\em #2}}                       
\def\bibref#1#2#3#4#5{#1 {\bf #2}:#3--#4 (#5)}        
\def\bibarticle#1#2#3#4#5#6#7{\bibtitle{#1}{#2},
\bibref{#3}{#4}{#5}{#6}{#7}.}

\def\JPA{J.\ Phys.\ A}
\def\PD{Physica D}
\def\PRA{Phys.\ Rev.\ A}

\def\PRE{Phys.\ Rev.\ E}
\def\PRL{Phys.\ Rev.\ Letters}
\def\SIAM{SIAM J.\ Appl.\ Math.}

\begin{document}


\title{A sample-paths approach to noise-induced synchronization: 
Stochastic resonance in a double-well potential}
\author{Nils Berglund and Barbara Gentz}
\date{}

\maketitle

\begin{abstract}
\noindent
Additive white noise may significantly increase the response of bistable
systems to a periodic driving signal. We consider two classes of
double-well potentials, symmetric and asymmetric, modulated periodically in
time with  period $1/\eps$, where $\eps$ is a moderately (not
exponentially) small parameter. We show that the response of the system
changes drastically when the noise intensity $\sigma$ crosses a threshold
value. Below the threshold, paths are concentrated near one potential well,
and have an exponentially small probability to jump to the other well.
Above the threshold, transitions between the wells occur with probability
exponentially close to $1/2$ in the symmetric case, and  exponentially
close to $1$ in the asymmetric case. The transition zones are localised in
time near the points of minimal barrier height. We give a mathematically
rigorous description of the behaviour of individual paths, which allows us,
in particular, to determine the power-law dependence of the critical noise
intensity on $\eps$ and on the minimal barrier height, as well as the
asymptotics of the transition and non-transition probabilities.
\end{abstract}

\leftline{\small{\it Date.\/} December 29, 2000.}
\leftline{\small 2000 {\it Mathematics Subject Classification.\/} 37H99,
60H10 (primary), 34F05, 34E15 (secondary).}


\section{Introduction}
\label{sec_in}

Since its introduction as a model for the periodic appearance of the
ice ages \cite{Benzi}, stochastic resonance has been observed in a
large number of physical  and biological systems, including lasers,
electronic circuits and the sensory system of crayfish (for reviews of
applications, see for instance \cite{MW}). 

The mechanism of stochastic resonance can be illustrated in a simple
model. Consider the overdamped motion of a particle in a double-well
potential. The two potential wells describe two macroscopically
different states of the unperturbed system, for instance cold and warm
climate. The particle is subject to two different kinds of
perturbation: a deterministic periodic driving force (such as the
periodic variation of insulation caused by the changing eccentricity
of the earth's orbit), and an additive noise (modeling the random
influence of the weather). 
Each of these two perturbations, taken by itself, does not produce any
interesting dynamics (from the point of view of resonance). Indeed,
the periodic driving is assumed to have too small an amplitude to
allow for any transitions between the potential wells in the absence
of noise. On the other hand, without periodic forcing, additive noise
will cause the particle to jump from one potential well to the other
at random times. The expected time between transitions is given
asymptotically, in the small noise limit, by Kramers' time, which is
proportional to the exponential of the barrier height $H$ over the
noise intensity squared, namely $\e^{H/\sigma^2}$. 
When both perturbations are combined, however, and their amplitudes
suitably tuned, the particle will flip back and forth between the
wells in a close to periodic way. Thus the internal noise can
significantly enhance the weak external periodic forcing, by producing
large amplitude oscillations of the system, hence the name of
resonance. 

The choice of the term \lq\lq resonance\rq\rq\ has been questioned, as
\lq\lq it would be more appropriate to refer to {\em noise-induced
signal-to-noise ratio enhancement\/}\rq\rq\ \cite{Fox}. In the regime
of a periodic driving whose amplitude is not a small parameter, one
also speaks of {\em noise-induced synchronization\/}
\cite{Neiman}. Appreciable, though still sub-threshold amplitudes of
the periodic driving have the advantage to enable transitions for
small noise intensities, without requiring astronomically long driving
periods.

While the heuristic mechanism of stochastic resonance is rather well
understood, a complete mathematical description is still lacking,
though important progress has been made in several limiting
cases. Depending on the regime one is interested in, several
approaches have been used to describe the phenomenon
quantitatively. The simplest ones use a discretization of either time
or space. When the potential is considered as piecewise constant in
time, the generator of the autonomous case can be used to give a
complete solution \cite{Benzi}, showing that resonance occurs when
driving period and Kramers' time are equal. Alternatively, space can
be discretized in order to obtain a two-state model, which is
described by a Markovian jump process \cite{ET}. The two-state model
has also been realised experimentally by an electronic circuit, called
the Schmitt trigger \cite{FH,McW}. 

In physical experiments, one has often access to indirect
characteristics of the dynamics, such as the power spectrum, which
displays a peak at the driving frequency. The strength of the
resonance is quantified through the signal-to-noise ratio (SNR), which
is proportional to the area under the peak (this definition obviously
leaves some liberty of choice). The SNR has been estimated, in the
limit of small driving amplitude, by using spectral theory of the
Fokker--Planck equation \cite{Fox,JH}, or a \lq\lq rate\rq\rq\
equation for the probability density \cite{McW}. The signal-to-noise
ratio is found to behave like $\e^{-H/\sigma^2}/\sigma^4$, which
reaches a maximum for $\sigma^2=H/2$. 

The probability density of the process, however, only gives part of the
picture, and a more detailed understanding of the behaviour of individual
paths is desirable. Some interesting progress in this direction is found in
\cite{Freidlin}. The approach applies to a very general class of dynamical
systems, in the limit of vanishing noise intensity. When the period of
the forcing scales like Kramers' time, solutions of the stochastic
differential equation are shown to converge to periodic functions in
the following sense: The $L^p$-distance between the paths and the
periodic limiting function converges to zero in probability as the
noise intensity goes to zero. Due to its generality, however, this
approach does not give any information on the rate of convergence of
typical paths to the periodic function, nor does it estimate the
probability of atypical paths. Also, since the period of the forcing
must scale like Kramers'  time, the assumed small noise intensity goes
hand in hand with exponentially long waiting times between interwell
transitions.

In the present work, we provide a more detailed description of the
individual paths' behaviour, for small but finite noise intensities
and driving frequencies.  We consider two classes of one-dimensional
double-well potentials, symmetric and asymmetric ones. The height of
the potential barrier is assumed to become small periodically, which
allows us to consider situations where the period need not be
exponentially large in $1/\sigma^2$ for transitions between the wells
to be likely.

In the case of an asymmetric potential, we are interested, in
particular, in determining the optimal noise  intensity as a function
of the driving frequency and the minimal barrier height, guaranteeing
a close-to-periodic oscillation between both wells. We will estimate
both the deviation (in space and time) of typical paths from the
limiting periodic function, and the asymptotics of the probability of
exceptional paths. The case of a symmetric potential shows an
additional feature. For the right choice of the noise intensity,
transitions become likely once per period, at which time the \lq\lq
new\rq\rq\ well is chosen at random. We will again estimate the
deviation from a suitable reference process and the asymptotics of the
probability of exceptional paths.

The systems are described by stochastic differential equations (SDEs)
of the form 
\begin{equation}
\label{i1}
\6x_s = -\dpar{}{x} V(x_s,s) \6s + \sigma \6W_s,
\end{equation}  
where $W_s$ is a Brownian motion. The potential $V(x,s)$ is
$1/\eps$-periodic in $s$, and admits two minima for every value of
$s$. The frequency $\eps$, the minimal barrier height between the
wells and the noise intensity $\sigma$ are considered as (moderately)
small parameters, the relation between which will determine the
transition probability. 

The first class of potentials we consider is symmetric in $x$. A
typical representative of this class is the potential
\begin{equation}
\label{i2}
V(x,s) = -\frac12 a(\eps s) x^2 + \frac14 x^4, 
\qquad\text{with $a(\eps s) = a_0 + 1 - \cos(2\pi\eps s)$.}
\end{equation} 
Here $a_0\geqs 0$ is a parameter controlling the minimal barrier height. 
We introduce the slow time $t=\eps s$ for convenience. The potential has
two wells, located at $\pm\sqrt{a(t)}$, separated by a barrier of height
$\frac14 a(t)^2$. The distance between the wells and the barrier
height become small simultaneously, at integer values of $t$. 
\begin{figure}
 \centerline{\psfig{figure=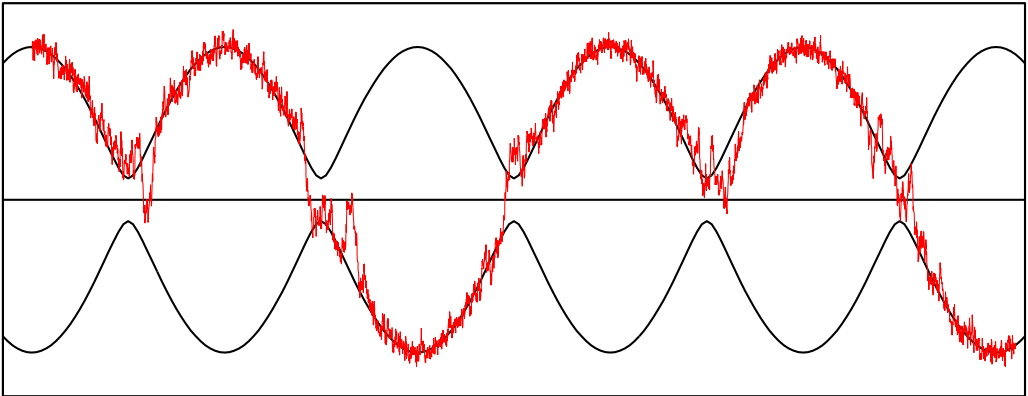,width=137mm,clip=t}}
 \caption[]
 {A typical solution of the SDE~\eqref{i1} in the case of the symmetric
 potential~\eqref{i2}. Heavy curves indicate the position of the potential
 wells, which approach each other at integer times. The straight line is
 the location of the saddle. Parameter values $\eps=0.01$, $\sigma=0.08$
 and $a_0=0.02$ belong to the regime where the transition probability
 between wells is close to $1/2$. We show that transitions are concentrated
 in regions of order $\sqrt\sigma$ around the instants of minimal barrier
 height.}
\label{fig1}
\end{figure}

\goodbreak

Our results for symmetric potentials can be summarized as follows:
\begin{itemiz}
\item	In the deterministic case $\sigma=0$, we describe the
dependence of solutions on $t$, $a_0$ and $\eps$
(Theorem~\ref{thm_sdet}). Solutions starting at $x>0$ are attracted by
the potential well at $\sqrt{a(t)}$, which they track with a small
lag. If $a_0\geqs \eps^{2/3}$, this lag is at most of the order
$\eps/a_0$; if $a_0\leqs \eps^{2/3}$, it is at most of the order
$\eps^{1/3}$, but solutions never approach the saddle closer than a
distance of order $\eps^{1/3}$ (even if $a_0=0$).  

\item	When noise is present, but $\sigma$ is small compared to the
maximum of $a_0$ and $\eps^{2/3}$, the paths are likely to track the
solution of the corresponding deterministic differential equation at a
distance of order $\sigma/\max\set{\abs{t},\sqrt{a_0},\eps^{1/3}}$
(Theorem~\ref{thm_snearxdet}). The probability to reach the saddle
during one time period is exponentially small in
$\sigma^2\!/\mskip-.5mu(\max\set{a_0,\eps^{2/3}})^2$.

\item	If $\sigma$ is larger than both $a_0$ and $\eps^{2/3}$,
transitions between potential wells become likely, but are
concentrated on the time interval
$[-\sqrt{\sigma},\sqrt{\sigma}\mskip1.5mu]$ (repeated
periodically). During this time interval, the paths may jump back and
forth frequently between both potential wells, and they have a typical
spreading of the order $\sigma/\max\set{\sqrt{a_0},\eps^{1/3}}$. After
time $\sqrt{\sigma}$, the paths are likely to choose one of the wells
and stay there till the next period (Theorem~\ref{thm_sescape}). The
probability to choose either potential well is exponentially close to
$1/2$, with an exponent of order $\sigma^{3/2}/\eps$, which is
independent of $a_0$ (Theorem~\ref{thm_strans}). 

\item	This picture remains true when $\sigma$ is larger than both
$\sqrt{a_0}$ and $\eps^{1/3}$, but note that the spreading of paths
during the transition may become very large. Thus increasing noise
levels will gradually blur the periodic signal.
\end{itemiz}
These results show a rather sharp transition to take place at 
$\sigma=\max\set{a_0,\eps^{2/3}}$, from a regime where the paths are
unlikely to switch from one potential well to the other one, to a
regime where they do switch with a probability exponentially close to
$1/2$ (\figref{fig1}). 

The second class of potentials we consider is asymmetric, a
typical representative being
\begin{equation}
\label{i3}
V(x,s) = -\frac12 x^2 + \frac14 x^4 - \lambda(\eps s)x.
\end{equation}
This is a double-well potential if and only if $\abs{\lambda}<\lc =
2/(3\sqrt3)$. We thus choose $\lambda(\eps s)=\lambda(t)$ of the form
\begin{equation}
\label{i4}
\lambda(t) = -(\lc-a_0)\cos(2\pi t).
\end{equation}
Near $t=0$, the right-hand potential well approaches the saddle at a
distance of order $\sqrt{a_0}$, and the barrier height is of order
\smash{$a_0^{3/2}$}. A similar encounter between the left-hand
potential well and the saddle occurs at $t=1/2$. 

\begin{figure}
 \centerline{\psfig{figure=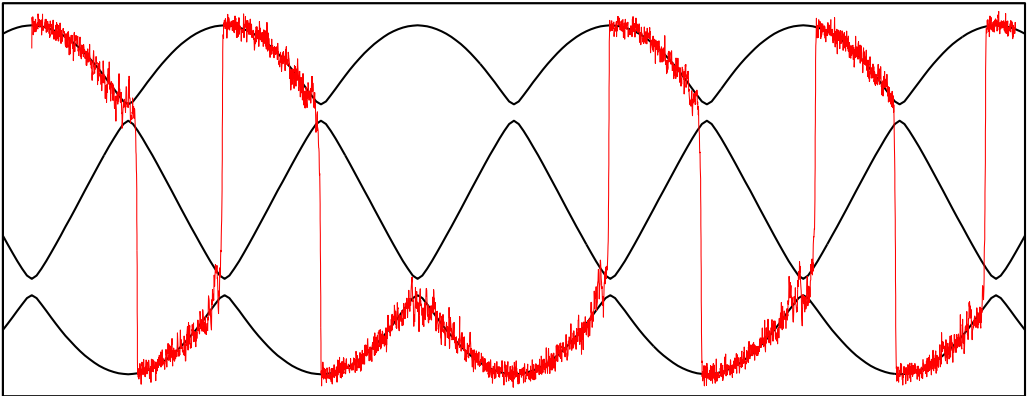,width=137mm,clip=t}}
 \caption[]
 {A typical solution of the SDE~\eqref{i1} in the case of the asymmetric
 potential~\eqref{i3}. The upper and lower heavy curves indicate the
 position of the potential wells, while the middle curve is the location of
 the saddle. Parameter values $\eps=0.005$, $\sigma=0.08$ and $a_0=0.005$
 belong to the regime where the transition probability between wells is
 close to $1$. We show that transitions are concentrated in
 regions of order $\sigma^{2/3}$ around the instants of minimal barrier
 height.}
\label{fig2}
\end{figure}
Our results for asymmetric potentials can be summarized as follows:
\begin{itemiz}
\item	In the deterministic case, solutions track the potential wells
at a distance at most of order
$\min\set{\eps/a_0,\sqrt\eps\mskip1.5mu}$. If $a_0\leqs\eps$, they
never approach the saddle closer than a distance of order $\sqrt\eps$
(Theorem~\ref{thm_adet}).  

\item	When $\sigma$ is small compared to the maximum of
\smash{$a_0^{3/4}$} and $\eps^{3/4}$, paths are likely to track the
deterministic solutions at a distance of order
$\sigma/\max\set{\sqrt{\abs{t}},\smash{a_0^{1/4},\eps^{1/4}}}$
(Theorem~\ref{thm_anearxdet}). The probability to overcome the barrier
is exponentially small in
$\sigma^2\!/(\max\!\smash{\set{a_0^{3/4},\eps^{3/4}}})^2$.

\item	For larger $\sigma$, transitions become probable during the
time interval $[-\sigma^{2/3},\sigma^{2/3}]$. Due to the asymmetry,
the probability to jump from the less deep potential well to the
deeper one is exponentially close to one, with an exponent of order
$\sigma^{4/3}/\eps$, while paths are unlikely to come back
(Theorem~\ref{thm_atrans}).  

\item	This picture remains true when $\sigma$ is larger than both
\smash{$a_0^{1/4}$} and \smash{$\eps^{1/4}$}, but the spreading of
paths during the transition may become very large. 
\end{itemiz}
Again, we find a rather sharp transition to take place, this time at
$\sigma=\max\smash{\set{a_0^{3/4},\eps^{3/4}}}$. In contrast to the
symmetric case, for large $\sigma$ the paths are likely to jump from 
one potential well to the other at every half-period (\figref{fig2}). 
\goodbreak

In both the symmetric and the asymmetric case, we thus obtain a high
switching probability between the potential wells even for small noise
intensities, provided minimal barrier height and driving frequency are
sufficiently small. They only need, however, to be smaller than a
power of $\sigma$: $a_0\ll\sigma$ and $\eps\ll\sigma^{3/2}$ in the
symmetric case, and $a_0\ll\sigma^{4/3}$, $\eps\ll\sigma^{4/3}$ in the
asymmetric case are sufficient conditions for switching dynamics.

Our results require a precise understanding of dynamical effects, and
the subtle interplay between the probability to reach the potential
barrier, the time needed for such excursions, and the total number of
excursions with a chance of success. In this respect, they  provide a
substantial progress compared to the \lq\lq quasistatic\rq\rq\
approach, which considers potentials that are piecewise constant in
time. Note that some of our results may come as a surprise. In
particular, neither the width (in time) of the transition zone nor the
asymptotics of the transition probability  depend on the minimal
barrier height $a_0$. In fact, the picture is independent of $a_0$ as
soon as $a_0$ is smaller than $\eps^{2/3}$ (in the symmetric case) or
$\eps$ (in the asymmetric case), even for $a_0=0$. This is due to the
fact that when $a_0$ is small, the time during which the potential
barrier is low is too short to contribute significantly to the
transition probability. 

The remainder of this paper is organized as follows. The results are
formulated in detail in Section~\ref{sec_res},
Subsection~\ref{ssec_rsym} being devoted to symmetric potentials, and
Subsection~\ref{ssec_rasym} to asymmetric
potentials. Section~\ref{sec_sym} contains the proofs for the
symmetric case, while Section~\ref{sec_asym} contains the proofs for
the asymmetric case.

\nobreak
\subsubsection*{Acknowledgements:}  
The fascinating problems related to stochastic resonance were brought 
to our attention by Anton Bovier, whom we would like to thank for
stimulating discussions. N.\,B. thanks the WIAS for kind hospitality.


\section{Results}
\label{sec_res}


\subsection{Preliminaries}
\label{ssec_rpre}

We consider non-autonomous SDEs of the form~\eqref{i1}. Introducing
the slow time $t=\eps s$ allows to study the system on a time interval
of order one. When substituting $t$ for $\eps s$, Brownian motion is
rescaled and we obtain an SDE
\begin{equation}
\label{SDE}
\6x_t = \frac{1}{\eps} f(x_t,t) \6t + \frac{\sigma}{\sqrt{\eps}} \6W_t,
\qquad x_{t_0}=x_0,
\end{equation}
where $f$ is the force, derived from the potential $V$, and
$\{W_t\}_{t\geqs t_0}$ is a standard Wiener process on some
probability space $(\Omega, \cF, \fP)$. Initial conditions $x_0$ are
always  assumed to be square-integrable with respect to $\fP$ and
independent of $\{W_t\}_{t\geqs t_0}$. Without further mentioning we
always assume that $f$ satisfies the usual (local) Lipschitz and
bounded-growth conditions which guarantee existence and pathwise
uniqueness of a strong solution $\{x_t\}_t$ of~\eqref{SDE}. Under
these conditions, there exists a continuous version of
$\{x_t\}_t$. Therefore we may assume that the paths $\omega\mapsto
x_t(\omega)$ are continuous for $\fP$-almost all $\omega\in\Omega$.  

We introduce the notation $\fP^{\mskip1.5mu t_0,x_0}$ for the law of
the process $\{x_t\}_{t\geqs t_0}$, starting in $x_0$ at time $t_0$,
and use $\E^{\mskip1.5mu t_0,x_0}$ to denote expectations with respect
to $\fP^{\mskip1.5mu t_0,x_0}$. Note that the stochastic process
$\{x_t\}_{t\geqs t_0}$ is an inhomogeneous Markov process. We are
interested in first exit times of $x_t$ from space--time sets. Let $\cA
\subset \R\times[t_0,t_1]$ be Borel-measurable. Assuming that $\cA$
contains $(x_0,t_0)$, we define the first exit time of $(x_t,t)$ from
$\cA$ by 
\begin{equation}
\label{pre1}
\tau_{\cA}=\inf\bigsetsuch{t\in[t_0,t_1]}{(x_t,t)\not\in\cA},
\end{equation}
and agree to set $\tau_{\cA}(\omega)=\infty$ for those
$\omega\in\Omega$ which satisfy $(x_t(\omega),t)\in\cA$ for all
$t\in[t_0,t_1]$. For convenience, we shall call $\tau_\cA$ the {\it
first exit time of $x_t$ from $\cA$}. Typically, we will consider sets
of the form $\cA=\setsuch{(x,t)\in\R\times[t_0,t_1]}{g_1(t)<x<g_2(t)}$
with continuous functions $g_1<g_2$.  Note that in this case,
$\tau_\cA$ is a stopping time\footnote{%
For a general Borel-measurable set $\cA$, the first exit time
$\tau_\cA$ is still a stopping time with respect to the canonical
filtration, completed by the null sets.}  
with respect to the canonical filtration
of $(\Omega, \cF, \fP)$ generated by $\{x_t\}_{t\geqs t_0}$.

Before turning to the precise statements of our results, let us
introduce some notations. We shall use
\begin{itemiz}
\item 
$\intpartplus{y}$ for $y\geqs 0$ to denote the smallest integer which
is greater than or equal to $y$, and 
\item 
$y\vee z$ and $y \wedge z$ to denote the maximum or minimum,
respectively, of two real numbers $y$ and $z$.
\item
If $\varphi(t,\eps)$ and $\psi(t,\eps)$ are defined for small $\eps$ and
for $t$ in a given interval $I$, we write
$\psi(t,\eps)\asymp\varphi(t,\eps)$ if there exist strictly positive
constants $c_\pm$ such that $c_-\varphi(t,\eps) \leqs \psi(t,\eps) \leqs
c_+\varphi(t,\eps)$ for all $t\in I$ and all sufficiently small
$\eps$. The constants $c_\pm$ are understood to be independent of $t$
and $\eps$ (and hence also independent of quantities like $\sigma$ and
$a_0$ which we consider as functions of $\eps$). 
\item 
By $g(u)=\Order{u}$ we indicate that there exist $\delta>0$ and $K>0$
such that $g(u)\leqs K u$ for all $u\in[0,\delta]$, where $\delta$ and
$K$ of course do not depend on $\eps$ or on the other small parameters
$a_0$ and $\sigma$. Similarly, $g(u)=\orderone{}$ is to be understood
as $\lim_{u\to0} g(u)=0$. 
\end{itemiz}
Finally, let us point out that most estimates hold for small enough
$\eps$ only, and often only for $\fP$-almost all $\omega\in\Omega$. We
will stress these facts only where confusion might arise.  


\subsection{Symmetric case}
\label{ssec_rsym}

We consider in this subsection the SDE \eqref{SDE} in the case of $f$
being periodic in $t$, odd in $x$, and admitting two stable
equilibrium branches, with a \lq\lq barrier\rq\rq\ between the
branches becoming small once during every time period. A typical
example of such a function is   
\begin{equation}
\label{sres1}
f(x,t) = a(t)x - x^3 
\qquad\text{with}\qquad 
a(t) = a_0 + 1 - \cos 2\pi t.
\end{equation}
We will consider a more general class of functions $f:\R^2\to\R$, which we
assume to satisfy the following hypotheses:
\begin{itemiz}
\item	{\it Smoothness:} 
$f \in \cC^4(\cM,\R)$, where $\cM = [-d,d\mskip1.5mu]\times\R$ and 
$d>0$ is a constant;

\item	{\it Periodicity:}
$f(x,t+1) = f(x,t)$ for all $(x,t)\in\cM$;

\item	{\it Symmetry:}
$f(x,t) = - f(-x,t)$ for all $(x,t)\in\cM$;

\item	{\it Equilibrium branches:}
There exists a continuous function $x^\star:\R\to(0,d\mskip1.5mu]$
with the property that $f(x,t)=0$ in $\cM$ if and only if $x=0$ or
$x=\pm x^\star(t)$; 

\item	{\it Stability:}
The origin is unstable and the equilibrium branches $\pm x^\star$ are
stable, that is, for all $t\in\R$,  
\begin{equation}
\label{sres2}
\begin{split}
a(t) &\defby \sdpar fx(0,t) > 0 \\
a^\star(t) &\defby \sdpar fx(x^\star(t),t) < 0.
\end{split}
\end{equation}

\item	{\it Behaviour near $t=0$:}
We want the three equilibrium branches to come close at integer times.
Given the symmetry of $f$, the natural assumption is that we have an \lq\lq
avoided pitchfork bifurcation\rq\rq, that is, 
\begin{equation}
\label{sres3}
\begin{split}
&\sdpar f{xxx}(0,0) < 0 \\
&a(t) = a_0 + a_1 t^2 + \Order{t^3},
\end{split}
\end{equation}
where $a_1>0$ and $\sdpar f{xxx}(0,0)$ are fixed (of order one), while $a_0
= a_0(\eps) = \orderone{\eps}$ is a positive small parameter. 
Is is easy to show that $x^\star(t)$ behaves like $\sqrt{a(t)}$ for small
$t$, and admits a quadratic minimum at a time $t^\star=\Order{a_0}$.
Moreover, $a^\star(t)\asymp -a(t)$ near $t=0$. 

We can choose a constant $T\in(0,1/2)$ such that the derivatives of
$a(t)$ and $x^\star(t)$ vanish only once in the interval $[-T,T]$. We
finally require that $x^\star(t)$, $a(t)$ and $a^\star(t)$ are bounded away
from zero outside this interval. We can summarize these properties as
\begin{align}
\label{sres4a}
x^\star(t) &\asymp 
\begin{cases}
\sqrt{a_0} & \text{for $\abs{t}\leqs\sqrt{a_0}$} \\
\abs{t} & \text{for $\sqrt{a_0}\leqs\abs{t}\leqs T$} \\
1 \phantom{MM}& \text{for $T\leqs t \leqs 1-T$,} 
\end{cases}
\\
\label{sres4b}
a(t) &\asymp 
\begin{cases}
a_0 & \text{for $\abs{t}\leqs\sqrt{a_0}$} \\
t^2 & \text{for $\sqrt{a_0}\leqs\abs{t}\leqs T$} \\
1 \phantom{MM}& \text{for $T\leqs t \leqs 1-T$,} 
\end{cases}
\\
\label{sres4c}
a^\star(t) &\asymp -a(t) \qquad\text{for all $t$.}
\end{align}
\end{itemiz}
We start by considering the deterministic equation
\begin{equation}
\label{sres5}
\eps\dtot{\xdet_t}t = f(\xdet_t,t).
\end{equation}
Without loss of generality, we may assume that $\xdet_t$ starts at time
$-1+T$ in some $\xdet_{-1+T}>0$. Tihonov's theorem~\cite{Grad,Tihonov},
applied on the interval $[-1+T,-T]$, implies that $\xdet_t$ converges
exponentially fast to a \nbh\ of order $\eps$ of $x^\star(t)$. We may thus
assume that $\xdet_{-T} = x^\star(-T)+\Order{\eps}$. In fact, since
$x^\star$ is decreasing at time $-T$, we may even assume that $\xdet_{-T} -
x^\star(-T)\asymp\eps$. 

The motion of $\xdet_t$ in the interval $[-T,T]$ is described in the
following theorem.

\begin{theorem}[Deterministic case]
\label{thm_sdet}
The solution $\xdet_t$ and the curve $x^\star(t)$ cross once and only
once during the time interval $[-T,T]$. This crossing occurs at a time
$\tilde t$ satisfying $\tilde t -t^\star \asymp
(\eps/a_0)\wedge\eps^{1/3}$. There exists a constant $c_0>0$ such that  
\begin{equation}
\label{sres6}
\xdet_t - x^\star(t) \asymp
\begin{cases}
\vrule height 14pt depth 14pt width 0pt
\phantom{-{}}\dfrac\eps{t^2} 
&\text{for $-T\leqs t\leqs -c_0(\sqrt{a_0}\vee\eps^{1/3})$} \\
\vrule height 14pt depth 12pt width 0pt
-\dfrac\eps{t^2} 
&\text{for $c_0(\sqrt{a_0}\vee\eps^{1/3})\leqs t\leqs T$,} 
\end{cases}
\end{equation}
and thus $\xdet_t\asymp\abs{t}$ in these time intervals. 
For $\abs{t}\leqs c_0(\sqrt{a_0}\vee\eps^{1/3})$, 
\begin{equation}
\label{sres7}
\xdet_t \asymp
\begin{cases}
\sqrt{a_0} 
&\text{if $a_0\geqs\eps^{2/3}$} \\
\eps^{1/3} 
&\text{if $a_0\leqs\eps^{2/3}$.} 
\end{cases}
\end{equation}
Finally, the linearization of $f$ at $\xdet_t$ satisfies
\begin{equation}
\label{sres8}
\ba(t) \defby \sdpar fx(\xdet_t,t) \asymp -(t^2\vee a_0\vee\eps^{2/3}).
\end{equation}
\end{theorem}

We give the proof in Subsection \ref{ssec_symdet}. 
The relation \eqref{sres7} may be surprising, since it means that no
matter how small we make $a_0$, $\xdet_t$ never approaches the saddle
at $x=0$ closer than a distance of order $\eps^{1/3}$. This fact can
be intuitively understood as follows. Even if $a_0=0$ and near $t=0$,
we have 
\begin{equation}
\label{sres9}
\eps\dtot{\xdet}{t} \geqs -\text{{\it const}} \;(\xdet)^3
\qquad\Rightarrow\qquad
\xdet_t \geqs \text{{\it const}} 
\frac{\xdet_{t_0}}{\sqrt{1+(\xdet_{t_0})^2(t-t_0)/\eps}}.
\end{equation}
Since \smash{$\xdet_{t_0}\asymp\eps^{1/3}$} for $t_0\asymp-\eps^{1/3}$,
$x_t$ cannot approach the origin significantly during any time interval of
order $\eps^{1/3}$. After such a time, however, the repulsion of the saddle
will make itself felt again, preventing the solution from further
approaching the origin. In other words, the time interval during which
$a(t)$ is smaller than $\eps^{2/3}$ is too short to allow the deterministic
solution to come close to the saddle. 

\begin{figure}
 \centerline{\psfig{figure=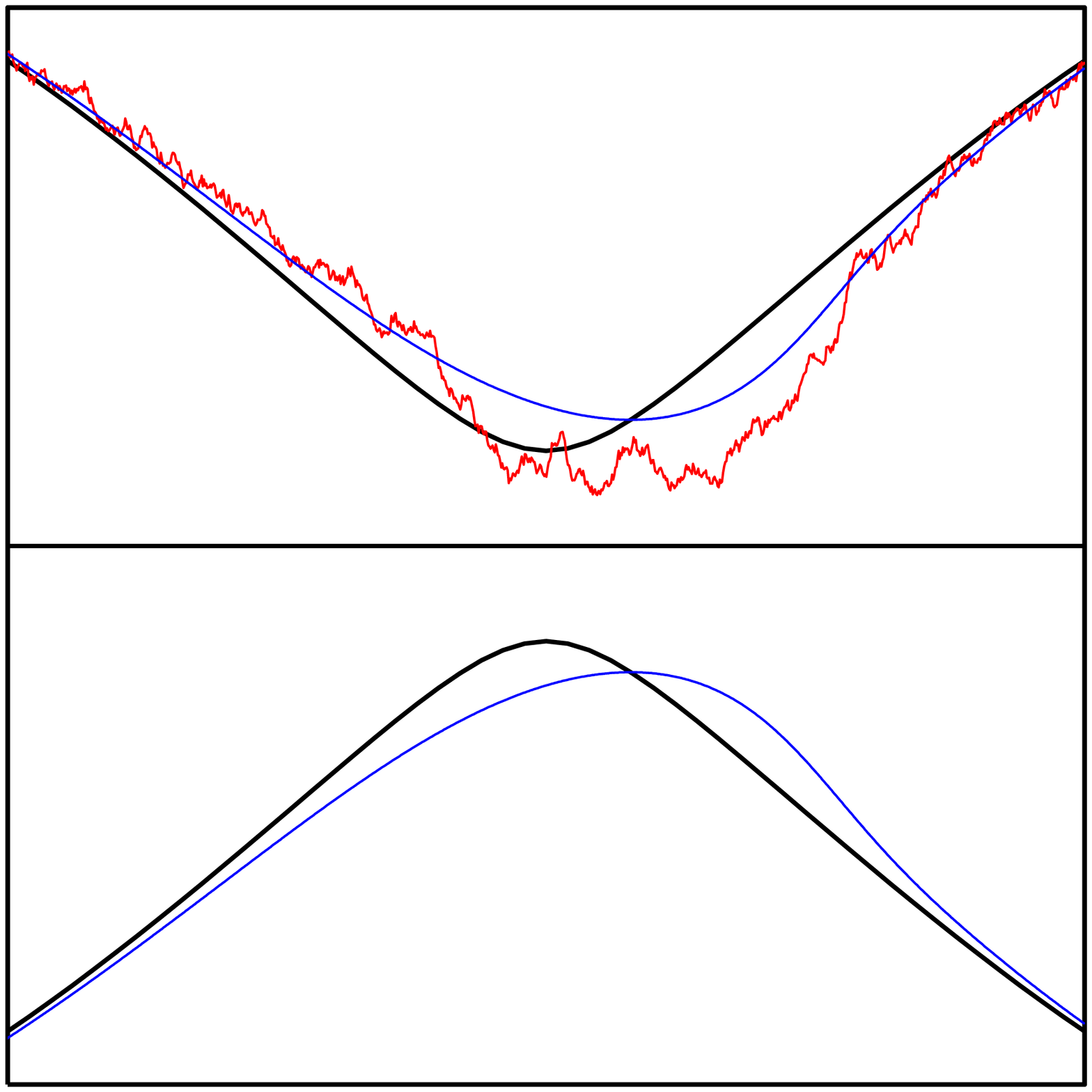,height=50mm,clip=t}
 \hspace{10mm}
 \psfig{figure=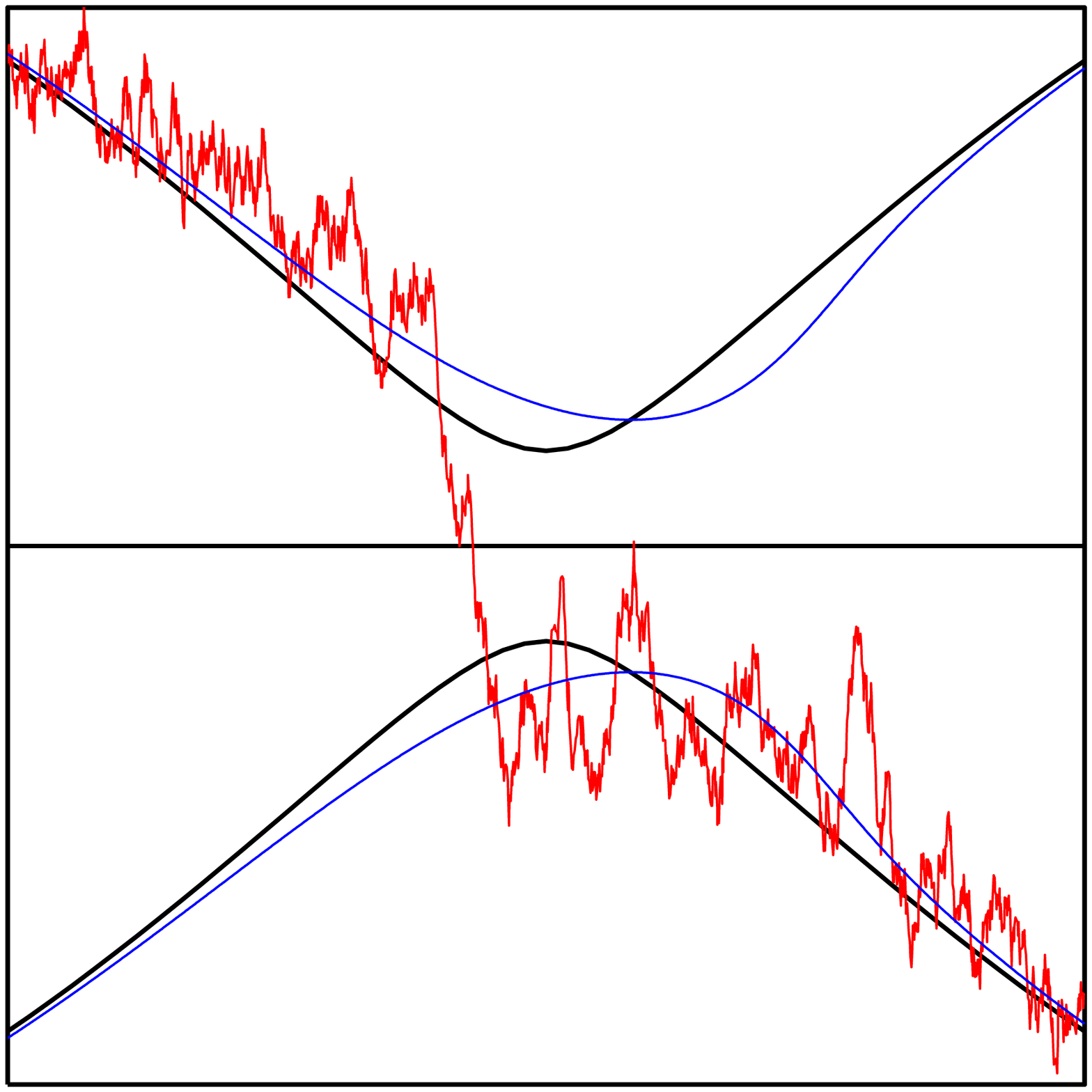,height=50mm,clip=t}
 }
 \caption[]
 {Solutions of the SDE \eqref{SDE} with symmetric drift term \eqref{sres1},
 shown for two different noise intensities, but for the same
 realization of  Brownian motion. Heavy curves represent the
 equilibrium branches $\pm  x^\star(t)$, and the straight line
 represents the saddle. Smooth light curves are solutions of the
 deterministic equation \eqref{sres5} tracking  the potential wells
 with a small lag, while rugged curves are paths of the
 SDE. Parameter values are $\eps=0.01$, $a_0=0.02$, $\sigma=0.02$
 (left)  and $\sigma=0.08$ (right).}
\label{fig3}
\end{figure}
We return now to the SDE \eqref{SDE} with $\sigma>0$. 
Assume that we start at some deterministic $x_{-1+T}>0$. Theorem 2.3 in
\cite{BG} shows that the paths are likely to track the deterministic
solution $\xdet_t$ with the same initial condition at a distance of order
$\sigma^{1-\delta}$ for any $\delta>0$ (with probability $\geqs 1-
(1/\eps^2) \exp\{-\text{\it const}/\sigma^{2\delta}\}$), as long as
the equilibrium branches are well separated, that is, at least for
$-1+T\leqs t\leqs -T$. A transition between the potential wells is
thus unlikely if $\sigma = \Order{\abs{\log \eps}^{-1/2\delta}}$,
and interesting phenomena can only be expected between the times $-T$
and $T$. Upon completion of one time period, i.\,e., at time $T$, the
Markov property allows to repeat the above argument. Hence there is no
limitation in considering the SDE \eqref{SDE} on the time interval
$[-T,T]$, with a fixed initial condition $x_{-T}$ satisfying
$x_{-T}-x^\star(-T)\asymp \eps$. We will denote by $\xdet_t$ and
$x_t$, respectively, the solutions of \eqref{sres5} and \eqref{SDE}
with the same initial condition $x_{-T}$.  
\goodbreak

Let us start by describing the dynamics in a \nbh\ of $\xdet_t$. The main
idea is that for $\sigma$ sufficiently small, the typical spreading of
paths around $\xdet_t$ should be related to the variance $v(t)$ of the
solution of \eqref{SDE}, linearized around $\xdet_t$. This variance is
given by
\begin{equation}
\label{sres10a}
v(t) =
\frac{\sigma^2}{\eps} \int_{-T}^t \e^{2\balpha(t,s)/\eps} \6s,
\qquad\text{where $\balpha(t,s)=\int_s^t \ba(u)\6u$.}
\end{equation}
The variance is equal to zero at time $-T$, but behaves asymptotically like
$\sigma^2/\abs{2\ba(t)}$. In fact, if we define the function
\begin{equation}
\label{sres10b}
\z(t) \defby \frac 1{2\abs{\ba(-T)}} \e^{2\balpha(t,-T)/\eps} +
\frac1{\eps} \int_{-T}^t \e^{2\balpha(t,s)/\eps} \6s,
\end{equation}
then $v(t)$ differs from $\sigma^2\z(t)$ by a term that becomes negligible
as soon as $\abs{\balpha(t,-T)}$ is larger than a constant times
$\eps\abs{\log\eps}$. $\z(t)$ has the advantage to be bounded away from
zero for all $t$, which avoids certain technical problems in the proofs. 
We shall show that 
\begin{equation}
\label{sres11}
\z(t) \asymp \frac1{t^2\vee a_0\vee\eps^{2/3}} 
\qquad \text{for $\abs{t}\leqs T$.}
\end{equation}
We introduce the set
\begin{equation}
\label{sres12}
\cB(h) = \bigsetsuch{(x,t)}{\abs{t}\leqs T, \abs{x-\xdet_t}<h\sqrt{\z(t)}},
\end{equation}
and denote by $\tau_{\cB(h)}$ the first exit time of $x_t$ from $\cB(h)$. 

\begin{theorem}[Motion near the stable equilibrium branches]
\label{thm_snearxdet}
There exists a constant $h_0$, depending only on $f$, such that
\begin{itemiz}
\item	if $-T\leqs t\leqs -(\sqrt{a_0}\vee\eps^{1/3})$ and $h<h_0t^2$,
then 
\begin{equation}
\label{sres13}
\bigprobin{-T,x_{-T}}{\tau_{\cB(h)}<t} 
\leqs C(t,\eps) \exp\biggset{-\frac12 \frac{h^2}{\sigma^2}
\biggbrak{1-\Order{\eps}-\biggOrder{\frac h{t^2}}}};
\end{equation}
\item	if $-(\sqrt{a_0}\vee\eps^{1/3})\leqs t\leqs T$ and
$h<h_0(a_0\vee\eps^{2/3})$,
then 
\begin{equation}
\label{sres14}
\bigprobin{-T,x_{-T}}{\tau_{\cB(h)}<t} 
\leqs C(t,\eps) \exp\biggset{-\frac12 \frac{h^2}{\sigma^2}
\biggbrak{1-\Order{\eps}-\biggOrder{\frac h{a_0\vee\eps^{2/3}}}}}.
\end{equation}
\end{itemiz}
In both cases,
\begin{equation}
\label{sres15}
C(t,\eps) = \frac1{\eps^2} \abs{\balpha(t,-T)} + 2.
\end{equation}
\end{theorem}

We give the proof in Subsection \ref{ssec_symnear}. 
This result has several consequences. Observe first that the exponential
factors in \eqref{sres13} and \eqref{sres14} are very small as soon as
$h$ is significantly larger than $\sigma$. The prefactor $C(t,\eps)$ (which,
unlike the exponent, we do not believe to be optimal) leads to
subexponential corrections, which are negligible as soon as $h/\sigma >
\Order{\abs{\log\eps}}$. It mainly accounts for the fact that the
probability for a path to leave $\cB(h)$ increases slowly with time. 
The theorem shows that the typical spreading of paths around $\xdet_t$
is of order 
\begin{equation}
\label{sres16}
\sigma\sqrt{\z(t)} \asymp
\frac{\sigma}{\abs{t}\vee\sqrt{a_0}\vee\eps^{1/3}}.
\end{equation}
If $\sigma\ll a_0\vee\eps^{2/3}$, we may choose $h\gg\sigma$ for all times,
and thus the probability of leaving a \nbh\ of $\xdet_t$, let alone approach
the other stable branch, is exponentially small (in
$\sigma^2/(a_0\vee\eps^{2/3})^2$). On the other hand, if $\sigma$ is not so
small, \eqref{sres13} can still be applied to show that a transition is
unlikely to occur before a time of order $-\sqrt\sigma$. Figure
\ref{fig3} illustrates this phenomenon by showing typical paths for two
different noise intensities. 

Let us now assume that $\sigma$ is sufficiently large to allow for a
transition, and  examine the transition regime in more detail. We will
proceed in two steps. First we will estimate the probability of {\em not}\/
reaching the saddle at $x=0$ during a time interval $[t_0,t_1]$. The
symmetry of $f$ implies that for any $t\geqs t_1$ and $x_0>0$, 
\begin{align}
\label{sres17}
\nonumber
\bigprobin{t_0,x_0}{x_{t}<0} 
&= \frac12\mskip1.5mu\bigprobin{t_0,x_0}{\exists s\in(t_0,t) \colon x_s=0}
\\
\nonumber
&= \frac12 - \frac12\mskip1.5mu\bigprobin{t_0,x_0}{x_s>0\;\forall
s\in[t_0,t]} \\
&\geqs \frac12 - \frac12\mskip1.5mu\bigprobin{t_0,x_0}{x_s>0\;\forall
s\in[t_0,t_1]}.
\end{align}
In the second step, we will show, independently, that paths are likely to
leave a \nbh\ of $x=0$ after time $\sqrt\sigma$. Thus if the
probability of not reaching $x=0$ is small, the probability of making
a transition from the positive well to the negative one will be close
to $1/2$ (it can never exceed $1/2$ because of the symmetry). This
does not exclude, of course, that paths frequently switch back and
forth between the two potential wells during the time interval
$[-\sqrt\sigma,\sqrt\sigma\mskip1.5mu]$. But it shows that
\eqref{sres17} can indeed be interpreted as a lower bound on the 
transition probability.

Let $\delta>0$ be a constant such that 
\begin{equation}
\label{sres18}
x\sdpar{f}{xx}(x,t)\leqs 0
\qquad
\text{for $\abs{x}\leqs\delta$ and $\abs{t}\leqs T$.}
\end{equation}
Our hypotheses on $f$ imply that such a $\delta$ of order one always exists.
In some special cases, for instance if $f(x,t)=a(t)x-x^3$, $\delta$ may be
chosen arbitrarily large. 
\goodbreak

\begin{theorem}[Transition regime]
\label{thm_strans}
Let $c_1>0$ be a constant and assume $c_1^2\sigma\geqs a_0\vee
\eps^{2/3}$. Choose times $-T\leqs t_0\leqs t_1 \leqs T$ with
$t_1\in[-c_1\sqrt\sigma, c_1\sqrt\sigma\mskip1.5mu]$, and let $h>2\sigma$ be
such that $\xdet_s+h\sqrt{\z(s)}<\delta$ for all
$s\in[t_0,t_1]$. Then, if $c_1$ is sufficiently small and
$x_0\in(0,\xdet_{t_0}+\frac12h\sqrt{\z(t_0)}]$,    
\begin{multline}
\label{sres19}
\qquad
\bigprobin{t_0,x_0}{x_s>0\;\forall s\in[t_0,t_1]} \\
\leqs \frac52 \Bigpar{\frac{\abs{\balpha(t_1,t_0)}}\eps + 1} 
\e^{-\bar\kappa h^2/\sigma^2} + 
\mskip 1.5mu 2 \exp\Bigset{-\bar\kappa \frac1{\log(h/\sigma)}
\frac{\alpha(t_1,-c_1\sqrt\sigma)}\eps},
\qquad
\end{multline}
where $\bar\kappa$ is a positive constant, and $\alpha(t,s)=\int_s^t
a(u)\6u$.  
\end{theorem}

The proof is given in Subsection \ref{ssec_symtrans}. 
The first term in \eqref{sres19} is an upper bound on the probability that
$x_t$ escapes \lq\lq upward\rq\rq. Indeed, our hypotheses on $f$ do not
exclude that other stable equilibria exist for sufficiently large $x$, which
might trap escaping trajectories. 
The second term bounds the probability of $x_s$ remaining between $0$ and
$\xdet_s+h\sqrt{\z(s)}$ for $-c_1\sqrt\sigma\leqs s\leqs t_1$. This estimate
lies at the core of our argument, and can be understood as follows. Assume
$x_s$ starts near $\xdet_s$. It will perform a certain number of excursions
to attempt reaching the saddle at $x=0$. Each excursion requires a typical
time of order $\Delta s$, such that $\alpha(s+\Delta s,s)\asymp\eps$ (that
is, $a(s)\Delta s\asymp\eps)$, in the sense that the probability of
reaching $0$ before time $s+\Delta s$ is small. After an unsuccessful
excursion, $x_s$ may exceed $\xdet_s$, but will return typically after
another time of order $\Delta s$. Thus the total number of trials during
the time interval $[-c_1\sqrt\sigma,t_1]$ is of order
$\alpha(t_1,-c_1\sqrt\sigma)/\eps$. Under the hypotheses of the theorem, the
probability of not reaching the saddle during one excursion is of order
one, and thus the total number of trials determines the exponent in
\eqref{sres19}.

\begin{figure}
 \centerline{\psfig{figure=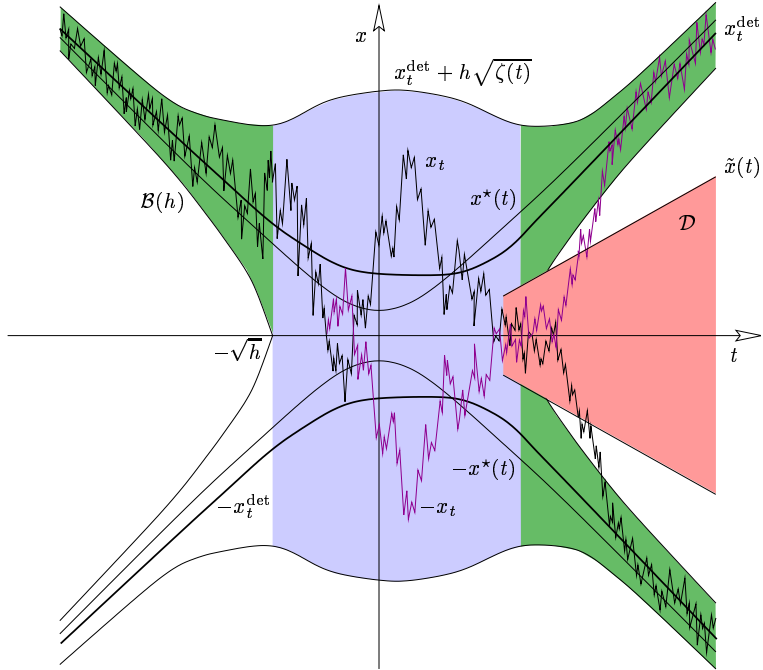,width=100mm,clip=t}}
 \caption[]
 {A typical path $x_t$ of the SDE \eqref{SDE} in the symmetric case, shown
 in a \nbh\ of the origin, where the potential barrier reaches its minimal
 height. We show here a situation where the noise intensity is large enough
 to allow a transition.  The potential wells at $\pm x^\star(t)$ behave like
 $\pm(\sqrt{a_0}\vee\abs{t})$. The deterministic solution $\xdet_t$
 starting near the right-hand potential well tracks $x^\star(t)$ at a
 distance at most of order $(\eps/a_0)\vee\smash{\eps^{1/3}}$, and never
 approaches the saddle at $x=0$ closer than
 $\smash{\sqrt{a_0}\vee\eps^{1/3}}$. The path $x_t$ is likely to stay in
 the set $\cB(h)$ up to time $-\smash{\sqrt h}$ when $h\gg\sigma$. Between
 times $-\smash{\sqrt\sigma}$ and $\smash{\sqrt\sigma}$, the path is likely
 to reach the origin. It may continue to jump back and forth between the
 potential wells up to time $\smash{\sqrt\sigma}$, but is likely to leave a
 \nbh\ $\cD$ of the saddle for times slightly larger than $\sqrt\sigma$.
 For each realization $\w$ such that $x_t(\w)$ reaches the saddle at a time
 $\tau$, there is a realization $\w'$ such that $x_t(\w')=-x_t(\w)$ is the
 mirror image of $x_t$ for $t\geqs\tau$, which explains why the probability
 to choose one well or the other after the transition region is close to
 $1/2$.}
\label{fig4}
\end{figure}
Before discussing the choice of the parameters giving an optimal bound
in Theorem~\ref{thm_strans}, let us first state the announced second
step, namely the claim that the paths are likely to escape from the
saddle after $t=c_1\sqrt\sigma$. For $\kappa\in(0,1)$, let us
introduce the set  
\begin{equation}
\label{sres19a}
\cD(\kappa) = \Bigsetsuch{(x,t)\in[-\delta,\delta]\times[c_1\sqrt\sigma,T]}
{\frac{f(x,t)}x > \kappa a(t)}.
\end{equation}
The upper boundary of $\cD(\kappa)$ is a function
$\tx(t)=\sqrt{1-\kappa}\,(1-\Order{t})x^\star(t)$. Let $\tau_{\cD(\kappa)}$
denote the first exit time of $x_t$ from $\cD(\kappa)$.

\begin{theorem}[Escape from the saddle]
\label{thm_sescape}
Let\/ $0<\kappa<1$ and assume $c_1^2\sigma\geqs
a_0\vee\eps^{2/3}$. Then there exist constants $c_2 \geqs c_1$ and 
$C_0>0$ such that 
\begin{equation}
\label{sres20}
\bigprobin{t_2,x_2}{\tau_{\cD(\kappa)}\geqs t}\leqs 
C_0 \Bigpar{\frac{t}{\sqrt\sigma}}^2 
\frac{\e^{-\kappa
\alpha(t,t_2)/2\eps}}{\sqrt{1-\e^{-\kappa\alpha(t,t_2)/\eps}}}, 
\end{equation}
for all $(x_2,t_2)\in\cD(\kappa)$ with $t_2 \geqs c_2\sqrt\sigma$. 
\end{theorem}

The proof is adapted from the proof of the similar Theorem 2.9 in
\cite{BG}. Compared to that result, we have sacrificed a factor $2$ in the
exponent, in order to get a weaker condition on $\sigma$. We discuss the
changes in the proof in Subsection \ref{ssec_symesc}.

For the moment, let us consider $t_2=c_2\sqrt\sigma$. We want to
choose a $t$ such that $\alpha(t,t_2)\geqs \eps\abs{\log\sigma}$. Since
$\alpha(t,t_2)$ is larger than a constant times $t_2^2(t-t_2)$, 
it suffices to choose a $t$ of order 
$\sqrt\sigma (1+\eps\abs{\log\sigma}/\sigma^{3/2})$
for~\eqref{sres20} to become small. Hence, after waiting for a time
of that order, we find
\begin{equation}
\label{sres20a}
\bigprobin{t_2,x_2}{\tau_{\cD(\kappa)}\geqs t} 
\leqs \text{\it const}\; \abs{\log\sigma} \sigma^{\kappa/2},
\end{equation}
which shows that most trajectories will have left $\cD(\kappa)$ by
time $t$, see \figref{fig4}. 

\goodbreak
It remains to show that the paths are likely to approach either
$x^\star(t)$ or $-x^\star(t)$ after leaving $\cD(\kappa)$. Let us first
consider the solution $\xhatdet_t$ of the deterministic differential
equation~\eqref{sres5} with initial time $t_2\geqs c_2\sqrt\sigma$ and
initial condition $\abs{x_2 - x^\star(t)}\leqs c\mskip1.5mu t_2$ for some
small constant $c>0$. Here we need to choose $c$ small in order to arrange
for $\ha(t)=\partial_xf(\xhatdet_t,t) \asymp -t^2$ which allows us to
proceed as in our investigation of the motion for $t\leqs -(\sqrt a_0 \vee
\eps^{1/3})$, cf.~Subsections~\ref{ssec_symdet} and~\ref{ssec_symnear}.
Under these assumptions, $\xhatdet_t$ approaches a neighbourhood of $\pm
x^\star(t)$ exponentially fast and then tracks the equilibrium branch at
distance $\eps/t^2$. As before, one can show that the path $x_t$ of the
solution of the SDE~\eqref{SDE} with the same initial condition is likely
to remain in a strip around $\xhatdet_t$ of width scaling with 
$\sigma/\sqrt{\ha(t)} \asymp \sigma/t$. So if a path $x_t$ leaves
$\cD(\kappa)$ at time $\tau_{\cD(\kappa)}$, then this path is likely to
approach $x^\star(t)$, if $x_{\tau_{\cD(\kappa)}}$ is positive, and
$-x^\star(t)$, otherwise. Note that
$\abs{x^\star(\tau_{\cD(\kappa)})-x_{\tau_{\cD(\kappa)}}}$ has  to be
smaller than $c\mskip1.5mu \tau_{\cD(\kappa)}$, which restricts the
possible values for $\kappa$. Therefore, we choose $\kappa$ small enough to
guarantee that $x^\star(t) - \tx(t) \leqs c \mskip1.5mu t$ for all $t\geqs
c_1\sqrt\sigma$. Finally note that paths which are not in $\cD(\kappa)$ at
time $c_2\sqrt\sigma$ but are not further away from $\pm x^\star$ than
$c\mskip1.5mu t$ at some time $t$ will also approach the corresponding
equilibrium branch.
\goodbreak

Let us now discuss the choice of the parameters in \eqref{sres19} giving an
optimal bound. For $\sigma\geqs (a_0\vee \eps^{2/3})/c_1^2$,
Theorem~\ref{thm_snearxdet} shows that up to time $t$ slightly less than
$-\sqrt{\sigma}$, the paths are concentrated around $\xdet_t$. Therefore,
here we should choose an initial time $t_0$ slightly before 
$-\sqrt{\sigma}$. In order for the second term in~\eqref{sres19} to be
small, we want to  choose $t_1-(-c_1\sqrt\sigma)$ as large as possible.
Note, however, that $h\sqrt{\z(s)}$ has to be smaller than $\delta-\xdet_s$
for all $s$. Since the order of $\z(s)$ is increasing for $s\leqs -(\sqrt
a_0 \vee\eps^{1/3})$, it turns out to be more advantageous to choose $t_1$
negative and of order $-\sqrt\sigma$, say $t_1=-\frac12 c_1\sqrt\sigma$. In
this case, $\abs{\alpha(t_1,-c_1\sqrt\sigma)}$ is larger than a constant
times $\smash{\sigma^{3/2}}$ (independently of $a_0$), and this does not
improve significantly for larger admissible $t_1$. At the same time, this
choice allows us to take $h\asymp \delta\sqrt\sigma$. We find
\begin{align}
\label{sres21}
\Bigprobin{t_0,x_0}{x_s>0\;\forall 
s\in[t_0,c_1\sqrt\sigma\,]}
&\leqs \Bigprobin{t_0,x_0}{x_s>0\;\forall 
s\in[t_0,-\tfrac12 c_1\sqrt\sigma\,]} \\
\nonumber
&\leqs \biggpar{\frac{\abs{t_0}}{\sqrt\sigma}}^3 \frac{\sigma^{3/2}}{\eps}
\e^{-\Order{\delta^2/\sigma}} 
+ \exp\biggset{-\frac{\text{\it const}}{\log(\delta^2/\sigma)} 
\frac{\sigma^{3/2}}{\eps}}.
\end{align}
Consider first the generic case $\delta\asymp 1$. The second term in
\eqref{sres21} becomes small as soon as  $\sigma/(\abs{\log
\sigma})^{2/3}\gg \eps^{2/3}$ holds in addition to the general
condition $c_1^2\sigma\geqs a_0\vee\eps^{2/3}$. The first term is
small as long as $\sigma=\Order{1/\log(\sigma/\eps^{2/3})}$. We thus
obtain the following regimes:   
\begin{itemiz}
\item	for $\sigma\leqs a_0\vee\eps^{2/3}$, the transition probability is
exponentially small in $\sigma^2/(a_0\vee\eps^{2/3})^2$;
\item	for $\sigma\geqs a_0/c_1^2$ with $(\eps\abs{\log\eps})^{2/3}
\ll \sigma \ll 1/\abs{\log\eps}$, the probability of a transition
between the wells is exponentially close to $1/2$, with an exponent
given essentially (up to logarithmic corrections) by 
\begin{equation}
\label{sres22}
\frac{\sigma^{3/2}}\eps \wedge \frac1\sigma;
\end{equation}
\item	for $\sigma\geqs 1/\abs{\log\eps}$, the paths become so poorly
localized that it is no longer meaningful to speak of a transition
probability.
\end{itemiz}
\eqref{sres22} shows that the transition probability becomes optimal for
$\sigma\asymp\eps^{2/5}$. For larger values of the noise intensity, the
possibility of paths escaping \lq\lq upward\rq\rq\ becomes sufficiently
important to decrease the transition probability. However, if the function
$f$ is such that $\delta$ can be chosen arbitrarily large, the second term
in \eqref{sres22} can be removed without changing the first one (up to
logarithmic corrections) by taking $\delta^2=\sigma/\eps$, for
instance. In that case, transitions between the wells become the more
likely the larger the ratio $\sigma/\eps^{2/3}$ is.

One should note that a typical path will reach maximal values of the
order $\sigma\sqrt{\z(0)} \asymp \sigma/(\sqrt{a_0}\vee\eps^{1/3})$. Thus,
due to the flatness of the potential near $t=0$, if $\sigma$ is larger than
$\sqrt{a_0}\vee\eps^{1/3}$, the spreading of the paths during the
transition interval is larger than the maximal distance between the
wells away from the transition. In general we cannot exclude that
paths escape to other attractors, if the potential has more than two wells. 

It may be surprising that the order of the transition probability is
independent of $a_0$ as soon as $\sigma>a_0$. Intuitively, one would rather
expect this probability to depend on the ratio $a_0^2/\sigma^2$, because of
Kramers' law. The fact that this is not the case illustrates the necessity
of a good understanding of dynamical effects (as opposed to a quasistatic
picture). Although the potential barrier is smallest between the times
$-\sqrt{a_0}$ and $\sqrt{a_0}$, the paths have more opportunities to reach
the saddle during larger time intervals. The optimal time interval turns
out to have a length of the order $\sqrt\sigma$, which corresponds to the
regime where diffusive behaviour prevails over the influence of the drift.


\subsection{Asymmetric case}
\label{ssec_rasym}

We consider in this subsection the SDE \eqref{SDE} in the case of $f$ being
periodic in $t$ and admitting two stable equilibrium branches, but without
the symmetry assumption. Instead, we want each of the potential wells to
approach the saddle once in every time period, but at different times for
the left-hand and the right-hand potential well. A typical example of such a
function is 
\begin{equation}
\label{ares1}
f(x,t) = x - x^3 + \lambda(t)
\qquad\text{with}\qquad 
\lambda(t) = -(\lc-a_0) \cos 2\pi t.
\end{equation}
Here $\lc = 2/(3\sqrt3)$ is defined by the fact that $f$ has two stable
equilibria if and only if $\abs{\lambda}<\lc$. Observe that $\sdpar fx$
vanishes at $x=\pm\xc=\pm1/\sqrt3$, and
\begin{equation}
\label{ares1a}
\begin{split}
f(\xc+y,t) &= \lc -(\lc-a_0) \cos 2\pi t - \sqrt3 y^2 - y^3 \\
&= a_0 + 2\pi^2(\lc-a_0)t^2 + \Order{t^4} - \sqrt3 y^2 - y^3.
\end{split}
\end{equation}  
Here the function $f(\xc,t)$ plays the role that $a(t)$ played in the
symmetric case, and near $t=0$ the right-hand potential well and the saddle
behave like $\xc\pm 3^{-1/4}\sqrt{f(\xc,t)}$, while the left-hand potential
well is isolated. Near $t=1/2$, a similar close encounter takes place
between the saddle and the left-hand potential well.

We will consider a more general class of functions $f:\R^2\to\R$, which we
assume to satisfy the following hypotheses:
\goodbreak
\begin{itemiz}
\item	{\it Smoothness:} 
$f \in \cC^3(\cM,\R)$, where $\cM = [-d,d\mskip1.5mu]\times\R$ and
$d>0$ is a constant; 

\item	{\it Periodicity:}
$f(x,t+1) = f(x,t)$ for all $(x,t)\in\cM$;

\item	{\it Equilibrium branches:}
There exist continuous functions $x^\star_- < x^\star_0 < x^\star_+$ from
$\R$ to $[-d,d\mskip1.5mu]$ with the property that $f(x,t)=0$ in $\cM$ if
and only if $x=x^\star_{\pm}(t)$ or $x=x^\star_0(t)$; the zeroes of $f$
should be isolated in the following sense: for every $\delta>0$, there
should exist a constant $\rho>0$ such that, if
$\abs{x-x^\star_\pm(t)}\geqs\delta$ and $\abs{x-x^\star_0(t)}\geqs\delta$,
then $\abs{f(x,t)}\geqs\rho$.\footnote{Since $f$ depends on a small
parameter $a_0$, we want to avoid that $f(x,t)$ approaches zero elsewhere
but near the three equilibrium branches, even when $a_0$ becomes small.}

\item	{\it Stability:}
The equilibrium branches $x^\star_{\pm}$ are stable and the equilibrium
branch $x^\star_0$ is unstable, that is, for all $t\in\R$,  
\begin{equation}
\label{ares2}
\begin{split}
a^\star_\pm(t) &\defby \sdpar fx(x^\star_{\pm}(t),t) < 0 \\
a^\star_0(t) &\defby \sdpar fx(x^\star_0(t),t) > 0.
\end{split}
\end{equation}

\item	{\it Behaviour near $t=0$:}
We want $x^\star_+$ and $x^\star_0$ to come close at integer times. Here
the natural assumption is that we have an \lq\lq avoided saddle--node
bifurcation\rq\rq, that is, there exists an $\xc\in(-\delta,\delta)$ such
that
\begin{equation}
\label{ares3}
\begin{split}
\sdpar f{xx}(\xc,0) &< 0 \\
\sdpar f{x}(\xc,t) &= \Order{t^2} \\
f(\xc,t) &= a_0 + a_1 t^2 + \Order{t^3},
\end{split}
\end{equation}
where $a_1>0$ and $\sdpar f{xx}(\xc,0)$ are fixed (of order one), while
$a_0 = a_0(\eps) = \orderone{\eps}$ is a positive small parameter. These
assumptions imply that $x^\star_+(t)$ reaches a local minimum at a time
$t^\star_+=\Order{a_0}$, and $x^\star_0(t)$ reaches a local maximum at a
possibly different time $t^\star_0=\Order{a_0}$. We can assume that for a
sufficiently small constant $T>0$, the three equilibrium branches and the
linearization of $f$ around them satisfy 
\begin{align}
\nonumber
x^\star_+(t) - \xc &\asymp 
\begin{cases}
\phantom-\sqrt{a_0} & \text{for $\abs{t}\leqs\sqrt{a_0}$} \\
\phantom-\abs{t} & \text{for $\sqrt{a_0}\leqs\abs{t}\leqs T$}
\end{cases}
\quad\;
&
a^\star_+(t) &\asymp 
\begin{cases}
-\sqrt{a_0} & \text{for $\abs{t}\leqs\sqrt{a_0}$} \\
-\abs{t} & \text{for $\sqrt{a_0}\leqs\abs{t}\leqs T$}
\end{cases}
\\
\nonumber
\vrule height 25pt depth 20pt width 0pt
x^\star_0(t) - \xc &\asymp 
\begin{cases}
-\sqrt{a_0} & \text{for $\abs{t}\leqs\sqrt{a_0}$} \\
-\abs{t} & \text{for $\sqrt{a_0}\leqs\abs{t}\leqs T$}
\end{cases}
&
a^\star_0(t) &\asymp 
\begin{cases}
\phantom-\sqrt{a_0} & \text{for $\abs{t}\leqs\sqrt{a_0}$} \\
\phantom-\abs{t} & \text{for $\sqrt{a_0}\leqs\abs{t}\leqs T$}
\end{cases}
\\
\label{ares4c}
x^\star_-(t) - \xc &\asymp -1 
\qquad\quad\text{for $\abs{t}\leqs T$}
&
a^\star_-(t) &\asymp -1
\qquad\quad\text{for $\abs{t}\leqs T$.}
\end{align}

\item	{\it Behaviour near $t=\tc$:}
We want $x^\star_-$ and $x^\star_0$ to come close at some time
$\tc\in (T,1-T)$. This is achieved by assuming that similar relations as
\eqref{ares3}, but with opposite signs, hold at a point $(\xc',\tc)$.

\item	{\it Behaviour between the close encounters:}
To exclude the possibility of other almost-bifurcations, we require that 
$x^\star_+(t) - x^\star_0(t)$ and $x^\star_0(t) - x^\star_-(t)$, as well as
the derivatives \eqref{ares2}, are bounded away from zero for $T<t<\tc-T$
and $\tc+T<t<1-T$. 
\end{itemiz}
Note that a sufficient assumption for the requirements on the behaviour
near $(\xc',\tc)$ to hold is that $f(x,t+\frac12)=-f(-x,t)$ for all
$(x,t)$. 

We start by considering the deterministic equation
\begin{equation}
\label{ares5}
\eps\dtot{\xdet_t}t = f(\xdet_t,t).
\end{equation}
As in the symmetric case, it is sufficient to consider the dynamics in the
time interval $[-T,T]$, with an initial condition satisfying $\xdet_{-T} -
x^\star_+(-T)\asymp\eps$. The situation in the time interval $[\tc-T,\tc+T]$
can be described in exactly the same way.

\begin{theorem}[Deterministic case]
\label{thm_adet}
The solution $\xdet_t$ and the curve $x_+^\star(t)$ cross once and only
once during the time interval $[-T,T]$. This crossing occurs at a time
$\tilde t$ satisfying $\tilde t -t^\star_+ \asymp
(\eps/\sqrt{a_0}\mskip1.5mu)\wedge\sqrt\eps$. There exists a constant
$c_0>0$ such that  
\begin{equation}
\label{ares6}
\xdet_t - x_+^\star(t) \asymp
\begin{cases}
\vrule height 14pt depth 14pt width 0pt
\phantom{-{}}\dfrac\eps{\abs{t}} 
&\text{for $-T\leqs t\leqs -c_0(\sqrt{a_0}\vee\sqrt\eps\mskip1.5mu)$} \\
\vrule height 14pt depth 12pt width 0pt
-\dfrac\eps{\abs{t}} 
&\text{for $c_0(\sqrt{a_0}\vee\sqrt\eps\mskip1.5mu)\leqs t\leqs T$,} 
\end{cases}
\end{equation}
and thus $\xdet_t-\xc\asymp\abs{t}$ in these time intervals. 
For $\abs{t}\leqs c_0(\sqrt{a_0}\vee\sqrt\eps\mskip1.5mu)$, 
\begin{equation}
\label{ares7}
\xdet_t -\xc \asymp
\begin{cases}
\sqrt{a_0} 
&\text{if $a_0\geqs\eps$} \\
\sqrt{\eps} 
&\text{if $a_0\leqs\eps$.} 
\end{cases}
\end{equation}
The linearization of $f$ at $\xdet_t$ satisfies
\begin{equation}
\label{ares8}
\ba(t) \defby \sdpar fx(\xdet_t,t) \asymp -(\abs{t}\vee
\sqrt{a_0}\vee\sqrt\eps\mskip1.5mu). 
\end{equation}
Moreover, \eqref{ares5} admits a particular solution $\xhatdet_t$ tracking
the {\/\em unstable} equilibrium branch $x_0^\star(t)$. It satisfies
analogous relations, namely, $\xhatdet_t$ and $x_0^\star(t)$ cross once at
a time $\hat t$ satisfying $\hat t-t^\star_0\asymp -(\tilde t-t^\star_+)$,
and \eqref{ares6}, \eqref{ares7} and \eqref{ares8} hold for $\xhatdet$ and
$x^\star_0(t)$, but with opposite signs.
\end{theorem}

The proof is similar to the proof of Theorem \ref{thm_sdet}, and we
comment on a few minor differences in Subsection \ref{ssec_asymdet}. 
Note that \eqref{ares7} implies that $\xdet_t$ never approaches the saddle
at $x^\star_0(t)$ closer than a distance of order $\sqrt\eps$.

We return now to the SDE \eqref{SDE} with $\sigma>0$. We will denote by
$\xdet_t$ and $x_t$, respectively, the solutions of \eqref{ares5} and
\eqref{SDE} with the same initial condition $x_{-T}$ satisfying $x_{-T} -
x^\star_+(-T)\asymp\eps$. We introduce again the function
\begin{equation}
\label{ares10}
\z(t) \defby \frac 1{2\abs{\ba(-T)}} \e^{2\balpha(t,-T)/\eps} +
\frac1{\eps} \int_{-T}^t \e^{2\balpha(t,s)/\eps} \6s,
\qquad\text{where $\balpha(t,s)=\int_s^t \ba(u)\6u$,}
\end{equation}
which behaves, in this case, like 
\begin{equation}
\label{ares11}
\z(t) \asymp \frac1{\abs{t}\vee\sqrt{a_0}\vee\sqrt\eps} 
\qquad \text{for $\abs{t}\leqs T$.}
\end{equation}
We define once more the set
\begin{equation}
\label{ares12}
\cB(h) = \bigsetsuch{(x,t)}{\abs{t}\leqs T, \abs{x-\xdet_t}<h\sqrt{\z(t)}}, 
\end{equation}
and denote by $\tau_{\cB(h)}$ the first exit time of $x_t$ from $\cB(h)$. 

\begin{theorem}[Motion near the stable equilibrium branches]
\label{thm_anearxdet}
There exists a constant $h_0$, depending only on $f$, such that
\begin{itemiz}
\item	if $-T\leqs t\leqs -(\sqrt{a_0}\vee\sqrt\eps\mskip1.5mu)$ and
$h<h_0\abs{t}^{3/2}$, then 
\begin{equation}
\label{ares13}
\bigprobin{-T,x_{-T}}{\tau_{\cB(h)}<t} 
\leqs C(t,\eps) \exp\biggset{-\frac12 \frac{h^2}{\sigma^2}
\biggbrak{1-\Order{\eps}-\biggOrder{\frac h{t^{3/2}}}}};
\end{equation}
\item	if $-(\sqrt{a_0}\vee\sqrt\eps\mskip1.5mu)\leqs t\leqs T$ and
$h<h_0(a_0^{3/4}\vee\eps^{3/4})$,
then 
\begin{equation}
\label{ares14}
\bigprobin{-T,x_{-T}}{\tau_{\cB(h)}<t} 
\leqs C(t,\eps) \exp\biggset{-\frac12 \frac{h^2}{\sigma^2}
\biggbrak{1-\Order{\eps}-\biggOrder{\frac h{a_0^{3/4}\vee\eps^{3/4}}}}}.
\end{equation}
\end{itemiz}
In both cases,
\begin{equation}
\label{ares15}
C(t,\eps) = \frac1{\eps^2} \abs{\balpha(t,-T)} + 2.
\end{equation}
\end{theorem}

This result is proved in exactly the same way as Theorem
\ref{thm_snearxdet}. It has similar consequences, only with different
values of the exponents. The typical spreading of paths around $\xdet_t$ is
of order
\begin{equation}
\label{ares16}
\sigma\sqrt{\z(t)} \asymp
\frac{\sigma}{\sqrt{\abs{t}}\vee a_0^{1/4}\vee\eps^{1/4}}.
\end{equation}
If $\sigma\ll a_0^{3/4}\vee\eps^{3/4}$, the probability of leaving a \nbh\
of $\xdet_t$, or making a transition to the other stable equilibrium
branch, is exponentially small (in
$\smash{\sigma^2/(a_0^{3/2}\vee\eps^{3/2})}$). On the other hand, if
$\sigma$ is not so small, \eqref{ares13} can still be applied to show that
a transition is unlikely to occur before a time of order $-\sigma^{2/3}$. 

Let us now assume that $\sigma$ is sufficiently large for a transition to
take place, i.e.\ that $\sigma\geqs a_0^{3/4}\vee\eps^{3/4}$. We want to
give an upper bound on the probability {\em not}\/ to make a transition.
Let us introduce levels $\delta_0 < \delta_1 < \xc < \delta_2$ such that 
\begin{equation}
\label{ares17}
\begin{split}
f(x,t) \asymp -1
&\qquad\qquad
\text{for $\delta_0\leqs x\leqs \delta_1$ and $\abs{t}\leqs T$} \\
\sdpar f{xx}(x,t) \leqs 0
&\qquad\qquad
\text{for $\delta_1\leqs x\leqs \delta_2$ and $\abs{t}\leqs T$.}
\end{split}
\end{equation} 
\goodbreak
Here $\delta_0$ and $\delta_1$ are always of order $1$ (in fact, we must
have $\delta_0>x^\star_-(t)$ for all $t\in[-T,T]$), and we think of
$\delta_0$ as being in the basin of attraction of $x^\star_-$. Our
hypotheses imply that a $\delta_2$ of order one satisfying \eqref{ares17}
always exists, but $\delta_2$ may be chosen arbitrarily large in particular
cases such as $f(x,t) = x - x^3 + \lambda(t)$. 

\begin{figure}
 \centerline{\psfig{figure=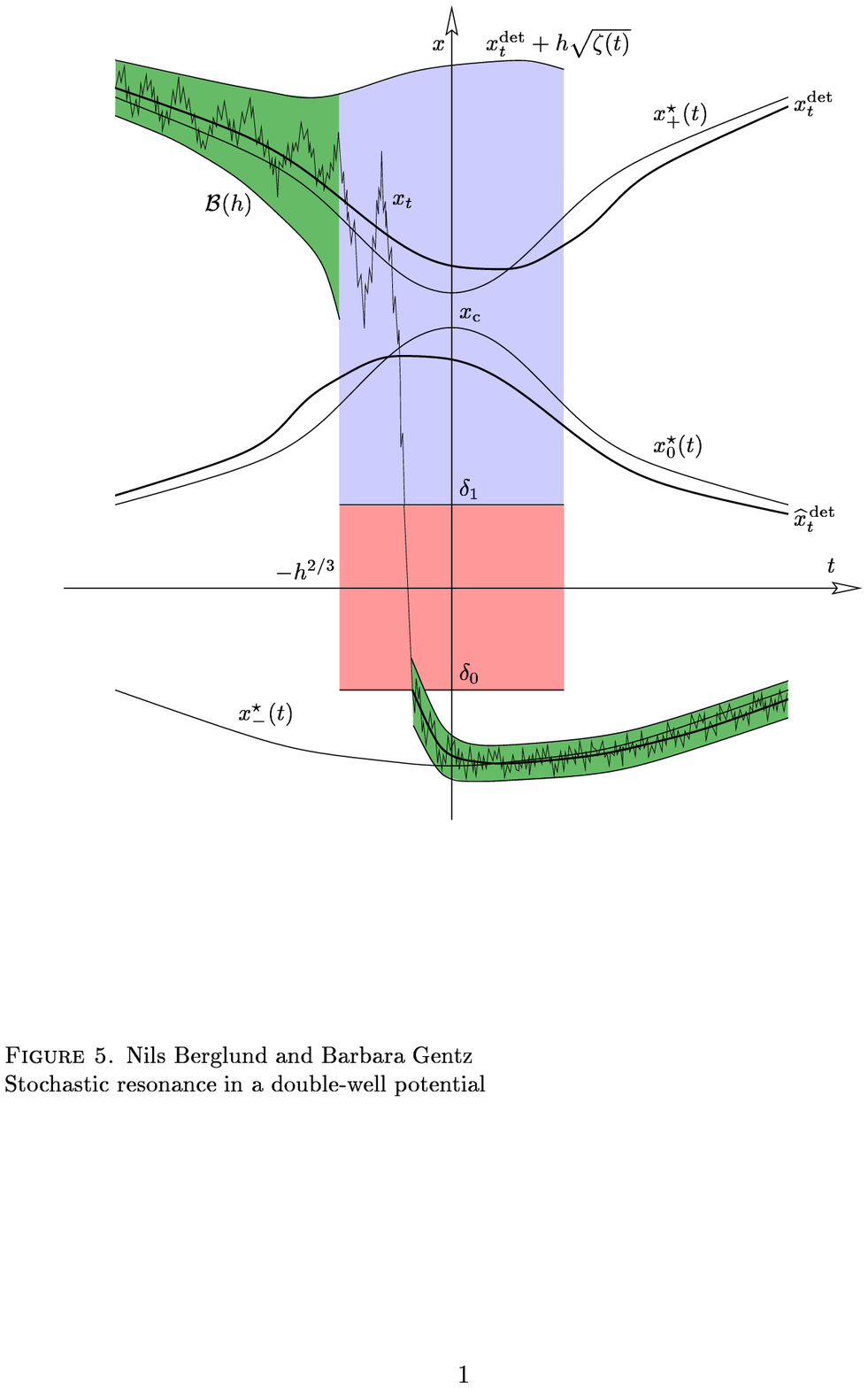,width=100mm,clip=t}}
 \caption[]
 {A typical path $x_t$ of the SDE \eqref{SDE} in the asymmetric case, shown
 near $t=0$, where the right-hand well at $x^\star_+$ approaches the saddle
 at $x^\star_0$. We show  again a situation where the noise intensity is
 large enough to allow a transition. The deterministic solution $\xdet_t$
 starting near the right-hand potential well tracks $x^\star_+(t)$ at a
 distance at most of order $(\smash{\eps/\sqrt{a_0})\wedge\sqrt{\eps}}$,
 and never approaches the saddle at $x^\star_0$ closer than
 $\smash{\sqrt{a_0}\vee\sqrt{\eps}}$. The path $x_t$ is likely to stay in
 the set $\cB(h)$ up to time $-\smash{h^{2/3}}$ when $h\gg\sigma$. Between
 times $-\smash{h^{2/3}}$ and $\smash{h^{2/3}}$, the path is likely to
 reach the saddle. Although it may fall back into the right-hand potential
 well, it is likely to finally overcome the potential barrier and reach a
 level $\delta_1$ of order $1$ below the saddle, after which it quickly
 reaches a lower level $\delta_0$. The distance between $\delta_1$ and
 $\delta_0$ can be much larger than in this picture. Finally, the path will
 track the deterministic solution starting in $x=\delta_0$, which
 approaches a \nbh\ of order $\eps$ of the left-hand potential well at
 $x^\star_-$.}
\label{fig5}
\end{figure}
The non-transition probability can be estimated by distinguishing three
cases:
\begin{itemiz}
\item	Either $x_t$, starting in $x_0>\xc$ at some $t_0<0$, never reaches
$\delta_1$. The probability of this event can be shown to be small in a
similar way as in Theorem \ref{thm_strans}, the main difference being that
due to the asymmetry, we can do better than estimating the probability not
to reach the saddle. 
\item	If $x_t$ reaches $\delta_1$, one can estimate in a very
simple way the probability not to reach $\delta_0$ as well, using the fact
that the drift term is bounded away from zero. 
\item	If $x_t$ reaches $\delta_0$, Theorem 2.3 in \cite{BG} shows that
$x_t$ is likely to reach a small \nbh\ of $x^\star_-$ as well.
\end{itemiz}

\begin{theorem}[Transition regime]
\label{thm_atrans}
Let $c_1$ and $c_2$ be positive constants and assume that
$\smash{c_1^{3/2}}\sigma\geqs \smash{a_0^{3/4}\vee \eps^{3/4}}$. Choose
times $-T\leqs t_0\leqs t_1 \leqs t\leqs T$ with $t_1\in[-c_1\sigma^{2/3},
c_1\sigma^{2/3}]$ and $t\geqs t_1+c_2\eps$. Let $h>2\sigma$ be such that
$\xdet_s+h\sqrt{\z(s)}<\delta_2$ for all $s\in[t_0,t_1]$. Then, for 
sufficiently small $c_1$, sufficiently large $c_2$ and all
$x_0\in(\delta_1,\xdet_{t_0}+\frac12h\sqrt{\z(t_0)}]$,   
\begin{align}
\nonumber
\bigprobin{t_0,x_0}{x_s>\delta_0\;\forall s\in[t_0,t]} 
\leqs{}& \frac52 \Bigpar{\frac{\abs{\balpha(t_1,t_0)}}\eps + 1} 
\e^{-\kappa h^2/\sigma^2}  \\
\nonumber &{}+ \frac32 
\exp\Bigset{-\kappa \frac1{\log(h/\sigma)\vee\abs{\log\sigma}}
\frac{\widehat\alpha(t_1,-c_1\sigma^{2/3})}\eps} \\
&{}+ \e^{-\kappa/\sigma^2},
\label{ares19}
\end{align}
where $\kappa$ is a positive constant, and $\widehat\alpha(t,s)=\int_s^t
\sdpar fx(\xhatdet_u,u)\6u$.  
\end{theorem}

The proof is given in Subsection \ref{ssec_asymtrans}. 
The three terms on the right-hand side of~\eqref{ares19} bound,
respectively, the probability that 
$x_t$ escapes through the upper boundary $\xdet_s+h\sqrt{\z(s)}$, the
probability that $x_t$ reaches neither the upper boundary nor $\delta_1$,
and the probability that $x_t$ does not reach $\delta_0$ when starting on
$\delta_1$ (\figref{fig5}). The crucial term is the second one. 

Let us now discuss the optimal choice of parameters.  If we choose $t_1 =
-\frac12 c_1\sigma^{2/3}$, we can take $h\asymp\tilde\delta_2\sigma^{1/3}$,
where $\tilde\delta_2=\delta_2-\xc$, and we get the estimate
\begin{multline}
\label{ares21}
\qquad
\Bigprobin{t_0,x_0}{x_s>\delta_0\;\forall 
s\in[t_0,t]} \\
\leqs \frac{t_0^2}\eps \e^{-\Order{\tilde\delta_2^2/\sigma^{4/3}}} 
+ \exp\biggset{-\frac{\text{\it const}}{\log(\tilde\delta_2^2/\sigma^{4/3})
\vee\abs{\log\sigma}} 
\frac{\sigma^{4/3}}{\eps}} + \e^{-\kappa/\sigma^2}.
\qquad
\end{multline}
As in the symmetric case, when $\delta_2\asymp 1$, we obtain the following
regimes:
\begin{itemiz}
\item	for $\sigma\leqs a_0^{3/4}\vee\eps^{3/4}$, the transition
probability is exponentially small in
$\sigma^2/(a_0^{3/4}\vee\eps^{3/4})^2$;
\item	for $a_0^{3/4}\vee\eps^{3/4} \ll \sigma \ll
(1/\abs{\log\eps})^{3/4}$, the transition probability is exponentially
close to $1$, with an exponent given essentially (up to logarithmic
corrections) by
\begin{equation}
\label{ares22}
\frac{\sigma^{4/3}}\eps \wedge \frac1{\sigma^{4/3}};
\end{equation}
\item	for $\sigma\geqs 1/\abs{\log\eps}$, the paths become so poorly
localized that it is no longer meaningful to speak of a transition
probability.
\end{itemiz}
The transition probability becomes optimal for $\sigma\asymp\eps^{3/8}$. 
Note, once again, that the exponent is independent of $a_0$. 

If the function $f$ is such that $\delta_2$ can be chosen arbitrarily
large, the second term in \eqref{ares22} can be removed without changing
the first one (up to logarithmic corrections) by taking
$\tilde\delta_2^2=\sigma^{8/3}/\eps$ for instance.  If $\sigma >
\smash{a_0^{1/4}\vee\eps^{1/4}}$, the paths may become extremely
delocalised in the transition zone, and could escape to other attractors.


\newpage
\section{Symmetric case}
\label{sec_sym}

We consider in this section the nonlinear SDE
\begin{equation}
\label{s1}
\6 x_\t = \frac{1}{\eps} f(x_\t,\t) \6\t +
\frac{\sigma}{\sqrt{\eps}} \6 W_\t,
\end{equation}
where $f$ satisfies the hypotheses given at the beginning of Subsection
\ref{ssec_rsym}. By rescaling $x$, we can arrange for $\sdpar{f}{xxx}(0,0)
= -6$, so that Taylor's formula allows us to write 
\begin{equation}
\label{s3}
\begin{split}
f(x,t) &= x\bigbrak{a(t)+g_0(x,t)} \\
\sdpar{f}{x}(x,t) &= a(t) + g_1(x,t),
\end{split}
\end{equation}
where $g_0, g_1\in\cC^3$ satisfy
\begin{equation}
\label{s4}
\begin{split}
g_0(x,t) &= \bigbrak{-1+r_0(x,t)}x^2 \\
g_1(x,t) &= \bigbrak{-3+r_1(x,t)}x^2, 
\end{split}
\end{equation}
with continuously differentiable functions $r_0, r_1$
satisfying $r_0(0,0)=r_1(0,0)=0$.

The implicit function theorem shows the existence, for small $t$, of
an equilibrium curve 
\begin{equation}
\label{s5}
x^\star(t) = \bigbrak{1+\bigOrder{\sqrt{a(t)}\,}} \sqrt{a(t)}.
\end{equation}
For small $t$, the curve $x^\star(t)$ behaves like
$\sqrt{a_0}\bigbrak{1+\Order{(t/\sqrt{a_0})^2}+\Order{\sqrt{a_0}}}$,
and it admits a quadratic minimum at some time $t^\star=\Order{a_0}$. 
Thus we can choose a constant $T\in(0,1/2)$ such that 
\begin{equation}
\label{s6}
x^\star(t) \asymp 
\begin{cases}
\sqrt{a_0} & \text{for $\abs{t}\leqs\sqrt{a_0}$} \\
\abs{t} & \text{for $\sqrt{a_0}\leqs\abs{t}\leqs T$} \\
1 & \text{for $T\leqs t \leqs 1-T$.} 
\end{cases}
\end{equation}


\subsection{Deterministic case}
\label{ssec_symdet}

In this subsection we consider the deterministic equation
\begin{equation}
\label{sdet1}
\eps\dtot{x_t}t = f(x_t,t).
\end{equation}
As already mentioned, Tihonov's theorem allows us to restrict the analysis
to the time interval $[-T,T]$, and to assume that $x_{-T} -
x^\star(-T)\asymp\eps$. 

\begin{remark}
During the time interval $[-T,T]$, the process $x_t$ crosses the
equilibrium branch $x^\star(t)$ once and only once, the time $\tilde
t$ of the crossing satisfying $\tilde t \geqs t^\star$. This fact
is due to the property that $x_t$ is strictly decreasing when lying
above $x^\star(t)$, and strictly increasing when lying below. Let
$\tilde t_1 < \tilde t_2 < \dots$ be the times of the successive
crossings of $x_t$ and $x^\star(t)$ in $[-T,T]$. Then $x_t$ is
decreasing between $-T$ and $\tilde t_1$ (since $x_{-T} - x^\star(-T)
>0$), increasing for $\tilde t_1<t<\tilde t_2$, and so on. Thus
$x^\star(t)$ must be increasing for $t$ slightly larger than $\tilde
t_1$, and decreasing for $t$ slightly larger than $\tilde t_2$. Since,
by assumption, $x^\star(t)$ is decreasing on $[-T,t^\star)$ and
increasing on $(t^\star,T]$, this implies that $\tilde t_1\geqs
t^\star$ and $\tilde t_2>T$. Therefore, there is at most one
crossing. We shall see below that $x_t$ and $x^\star(t)$ actually
cross and we will also determine the order of that time $\tilde t$.
\end{remark}

We consider now the difference $y_t = x_t - x^\star(t)$. It satisfies the
equation
\begin{equation}
\label{sdet2}
\eps\dtot yt = a^\star(t)y + b^\star(y,t) - \eps\dtot{x^\star}t,
\end{equation}
where Taylor's formula, \eqref{s3} and \eqref{s4} yield the relations
\begin{align}
\label{sdet3a}
a^\star(t) &= -2 a(t) \bigbrak{1+\bigOrder{\sqrt{a(t)}\,}}
\asymp 
\begin{cases}
-a_0 & \text{for $\abs{t}\leqs\sqrt{a_0}$} \\
-t^2 & \text{for $\sqrt{a_0}\leqs\abs{t}\leqs T$} 
\end{cases} \\
\label{sdet3b}
\vrule height 15pt depth 15pt width 0pt
b^\star(y,t) &= - \bigpar{3 x^\star(t)+y} y^2 
\bigbrak{1+\Order{x^\star(t)+y}} \\
\label{sdet3c}
\dtot{x^\star}t (t) &\asymp 
\begin{cases}
-1 & \text{for $-T\leqs t\leqs -\sqrt{a_0}$} \\
\dfrac{(t-t^\star)}{\sqrt{a_0}} &\text{for $\abs{t}\leqs\sqrt{a_0}$} \\
1 & \text{for $\sqrt{a_0}\leqs t\leqs T$,}
\end{cases}
\end{align}
with $t^\star = t^\star(a_0) = \Order{a_0}$. 

We start by giving a technical result that we will need several times.

\begin{lemma}
\label{lem_sdet1}
Let $\ta(t)$ be a continuous function satisfying $\ta(t)\asymp-(\beta\vee
t^2)$ for $\abs{t}\leqs T$, where $\beta=\beta(\eps)\geqs 0$. Let
$\chi_0\asymp 1$, and define $\talpha(t,s)=\int_s^t \ta(u)\6u$. Then 
\begin{equation}
\label{sdet4}
\chi_0 \e^{\talpha(t,-T)/\eps} + \frac1\eps\int_{-T}^t
\e^{\talpha(t,s)/\eps}\6s \asymp 
\begin{cases}
\vrule height 14pt depth 17pt width 0pt
\dfrac1{\beta\vee\eps^{2/3}} &\text{for $\abs{t}\leqs
\sqrt\beta\vee\eps^{1/3}$} \\
\vrule height 17pt depth 14pt width 0pt
\dfrac1{t^2} &\text{for $\sqrt\beta\vee\eps^{1/3}\leqs\abs{t}\leqs T$.}
\end{cases}
\end{equation}
\end{lemma}
\begin{proof}
To prove the lemma, we take advantage of the fact that the expression on
the left-hand side of~\eqref{sdet4} is the solution of an ordinary
differential equation. By the semi-group property, we may consider
separately the following regimes: For $\beta\geqs\eps^{2/3}$, we
distinguish the cases $t\in[-T,-T/2]$, $[-T/2,-\sqrt\beta\mskip1.5mu]$,
$[-\sqrt\beta,-\eps/a_0]$, $[-\eps/a_0,\eps/a_0]$,
$[\eps/a_0,\sqrt\beta\mskip1.5mu]$, $[\sqrt\beta,T]$   and for
$\beta<\eps^{2/3}$, we deal separately with  $t\in[-T,-T/2]$,
$[-T/2,-\eps^{1/3}]$, $[-\eps^{1/3},-\sqrt\beta\mskip1.5mu]$,
$[-\sqrt\beta,\sqrt\beta\mskip1.5mu]$, $[\sqrt\beta,\eps^{1/3}]$,
$[\eps^{1/3},T]$.
On each of these time intervals the claimed behaviour follows easily
by elementary calculus, see also the similar result~\cite[Lemma~4.2]{BG}.
\end{proof}

\begin{prop}
\label{prop_sdet1}
There exists a constant $c_0>0$ such that the solution of \eqref{sdet2} with
initial condition $y_{-T}\asymp\eps$ satisfies
\begin{equation}
\label{sdet5}
y_t \asymp \frac\eps{t^2} 
\qquad\text{for $-T\leqs t\leqs -c_0(\sqrt{a_0}\vee\eps^{1/3})$.}
\end{equation}
\end{prop}
\begin{proof}
Let $c_0\geqs 1$ and $c_1>T^2y_{-T}/\eps$ be constants to be chosen
later, and denote by $\tau$ the first exit time of $y_t$ from the strip
$0<y_t<c_1\eps/t^2$. Set $t_0=-c_0(\sqrt{a_0}\vee\eps^{1/3})$. Then, 
for $-T\leqs t\leqs\tau\wedge t_0$, we get from \eqref{sdet3b} and
\eqref{s6} that 
\begin{equation}
\label{sdet1:1}
\frac{\abs{b^\star(y,t)}}{yt^2} 
\leqs M \frac{3x^\star(t)+y}{\abs t} \frac{y}{\abs t}
\leqs M' \biggpar{1+c_1\frac{\eps}{\abs{t}^3}}c_1\frac{\eps}{\abs{t}^3}
\leqs M' \biggpar{1+\frac{c_1}{c_0^3}}\frac{c_1}{c_0^3},
\end{equation}
for some constants $M,M'>0$. The relations \eqref{sdet3a} and
\eqref{sdet3c} yield the existence of constants $c_\pm>0$ such that
$a^\star(t)\leqs-c_-t^2$ and $-\dtot{}tx^\star(t)\leqs c_+$ for
$t\in[-T,t_0]$. From \eqref{sdet2} and \eqref{sdet1:1} we obtain
\begin{equation}
\label{sdet1:3}
\eps\dtot yt \leqs -c_- t^2 y \biggbrak{1 -
\frac{M'}{c_-}\biggpar{1+\frac{c_1}{c_0^3}}\frac{c_1}{c_0^3}} + \eps c_+.
\end{equation}
For any given $c_1$, we can choose $c_0$ large enough for the term in
brackets to be larger than $1/2$. Then, by Lemma \ref{lem_sdet1},
there exists a constant $c_2=c_2(c_+,c_-)>0$ such that 
\begin{equation}
\label{sdet1:2}
y_t \leqs y_{-T} \e^{-c_-(t^3+T^3)/6\eps} +\mskip1.5mu c_+
\int_{-T}^t \e^{-c_-(t^3-s^3)/6\eps}\6s \leqs c_2 \frac{\eps}{t^2}
\end{equation}
for all $t\in[-T,\tau\wedge t_0]$. Therefore, if $c_1>c_2$, then
$\tau\geqs t_0$ follows. 

The lower bound can be obtained in exactly the same way.
\end{proof}

For the remainder of this subsection, let $t_0=-c_0(\sqrt{a_0}\vee
\eps^{1/3})$ with $c_0$ chosen according to the preceding proposition. Note
that this proposition implies that $x_t \asymp x^\star(t) \asymp \abs{t}$
for $-T \leqs t \leqs t_0$ and, in particular, that $y_{t_0} \asymp
(\eps/a_0)\wedge \eps^{1/3}$.

We now consider the dynamics for $\abs{t}\leqs \abs{t_0}$, starting
with the case of $a_0$ not too small, i.\,e., the case of $y_{t_0} \asymp
\eps/a_0$.

\begin{prop}
\label{prop_sdet2}
There exists a constant $\gamma_0>0$, depending only on $f$ and
$y_{t_0}$, such that, when $a_0 \geqs \gamma_0\eps^{2/3}$, then
\begin{equation}
\label{sdet6}
y_t = C_1(t)(t^\star-t) + C_2(t) 
\qquad\text{with}\qquad 
C_1(t)\asymp\frac{\eps}{a_0^{3/2}},
\quad 
C_2(t)\asymp\frac{\eps^2}{a_0^{5/2}}
\end{equation}
for all $\abs{t}\leqs \abs{t_0}$.
\end{prop}

\begin{proof}
Again, we will only show how to obtain an upper bound, since the
corresponding lower bound can be established in exactly the same way.

First we fix a constant
$c_1 > a_0y_{t_0}/\eps + 2(t^\star-t_0)/\sqrt{a_0} + 4\eps/(c_-a_0^{3/2})$.
We denote by $\tau$ the first exit time of $y_t$ from the strip
$\abs{y_t} < c_1\eps/a_0$. For $t_0 \leqs t \leqs\tau\wedge\abs{t_0}$,
we have 
\begin{equation}
\label{sdet2:4}
\frac{\abs{b^\star(y,t)}}{\abs{y}} \leqs M
\bigpar{3x^\star(t)+\abs{y}}\abs{y}
\leqs M'a_0 \biggpar{1+c_1\frac\eps{a_0^{3/2}}} c_1\frac\eps{a_0^{3/2}}
\end{equation}
with constants $M,M'>0$. Choosing $\gamma_0$ large enough, we get
\begin{equation}
\label{sdet2:5}
\eps\dtot yt 
\leqs -\frac{c_-}{2}a_0 y - c_-\frac{\eps}{\sqrt{a_0}}(t-t^\star),
\end{equation}
which implies
\begin{align}
\label{sdet2:1}
\nonumber
y_t &\leqs y_{t_0} \e^{-c_-a_0(t-t_0)/2\eps} 
- \frac{c_-}{\sqrt{a_0}} \int_{t_0}^t \e^{-c_- a_0(t-s)/2\eps}
(s-t^\star) \6s\\  
&= \frac{2\eps}{a_0^{3/2}}(t^\star-t) + \eta(\eps) \e^{-c_-
a_0(t-t_0)/2\eps} 
+ \mskip1.5mu\frac{4\eps^2}{c_-a_0^{5/2}}
\end{align}
by partial integration. Here
\begin{equation}
\label{sdet2:1a}
\eta(\eps) = y_{t_0} -2\frac{\eps}{a_0}
\biggpar{\frac{t^\star-t_0}{\sqrt{a_0}} + \frac{2\eps}{c_- a_0^{3/2}}}
\end{equation}
satisfies $\eta(\eps) = \Order{\eps/a_0}$. 

We want to estimate the contribution of the middle term on the
right-hand side of~\eqref{sdet2:1}. Assume first that $\eta(\eps)>0$
and consider $t \leqs t^\star$. By convexity,
\begin{equation}
\label{sdet2:1b}
\e^{-c_- a_0(t-t_0)/2\eps}
\leqs \frac{t^\star-t}{t^\star-t_0} + \e^{-c_- a_0(t^\star-t_0)/2\eps}.
\end{equation}
Now, 
\begin{equation}
\label{sdet2:1c}
\frac{2\eps}{a_0^{3/2}} + \eta(\eps) \frac{1}{t^\star-t_0}
\asymp \frac{\eps}{a_0^{3/2}}.
\end{equation}
Since $X\e^{-X}\to 0$ as $X\to\infty$, we also have
\begin{equation}
\label{sdet2:1d}
\eta(\eps)\e^{-c_- a_0(t^\star-t_0)/2\eps} + \frac{4\eps^2}{c_-a_0^{5/2}}
\asymp \frac{\eps^2}{a_0^{5/2}},
\end{equation}
provided $\gamma_0$ is large enough. This shows the existence of
constants $\overbar C_1\geqs2$ and $\overbar C_2>0$ such that
\begin{equation}
\label{sdet2:1e}
y_t \leqs \overbar C_1 \frac{\eps}{a_0^{3/2}} (t^\star-t) 
          + \overbar C_2 \frac{\eps^2}{a_0^{5/2}}
\qquad\text{for $t\leqs t^\star\wedge\tau$.}
\end{equation}
For $t \geqs t^\star$,
\begin{equation}          
\label{sdet2:1f}
\e^{-c_- a_0(t-t_0)/2\eps}
\leqs \e^{-c_- a_0(t^\star-t_0)/2\eps}
\end{equation}
is immediate, and~\eqref{sdet2:1d} shows that~\eqref{sdet2:1e} also
holds for $t^\star \leqs t\leqs\tau$. Note that in the case
$\eta(\eps)\leqs0$, ~\eqref{sdet2:1e} holds trivially.
Since $y_t <c_1 \eps/a_0$ is a direct consequence of~\eqref{sdet2:1}
and our choice of $c_1$, $\tau\geqs\abs{t_0}$ follows, and, therefore,
the upper bound~\eqref{sdet2:1e} holds for all $|t|\leqs|t_0|$.  
\end{proof}

Note that the result~\eqref{sdet6} implies that $y_t$ changes sign at a
time $t^\star+\Order{\eps/a_0}$, which shows that $x_t$ actually crosses
$x^\star(t)$ at a time $\tilde t$ satisfying $\tilde t-t^\star \asymp
\eps/a_0$. For large enough $\gamma_0$, the proposition also shows that
$x_t \asymp \sqrt{a_0}$ for $\abs{t} \leqs \abs{t_0}$ and that
$y_{\abs{t_0}} \asymp -\eps/a_0$.  

We consider now the case $a_0<\gamma_0\eps^{2/3}$ with $\gamma_0$ from
Proposition~\ref{prop_sdet2}. Without loss of generality, we may
assume that $\gamma_0\geqs1$. 

\begin{prop}
\label{prop_sdet3}
Assume that $a_0<\gamma_0\eps^{2/3}$. Then, for any fixed $t_1\asymp
\eps^{1/3}$, 
\begin{equation}
\label{Sdet7}
x_t \asymp \eps^{1/3} 
\qquad\text{for $t_0\leqs t\leqs t_1$,}
\end{equation}
and $x_t$ crosses $x^\star(t)$ at a time $\tilde t$ satisfying
$\tilde t \asymp \eps^{1/3}$. 
\end{prop}
\begin{proof}
In order to show~\eqref{Sdet7}, we rescale space and time in the
following way:
\begin{equation}
\label{sdet3:1}
x = \eps^{1/3} a_1^{1/6} z, 
\qquad
t = \eps^{1/3} a_1^{-1/3} s.
\end{equation}
Let $s_0 = \eps^{-1/3}a_1^{1/3} t_0$. 
Then $z_{s_0}\asymp 1$, and $z$ satisfies the differential equation
\begin{equation}
\label{sdet3:2}
\dtot zs = \ta(s,\eps) z + \bigbrak{-1+r_0(\eps^{1/3} a_1^{1/6} z,
\eps^{1/3} a_1^{-1/3} s)}z^3,
\end{equation}
\goodbreak
\noindent where
\begin{equation}
\label{sdet3:3}
\ta(s,\eps) = \frac{a(\eps^{1/3} a_1^{-1/3} s)}{\eps^{2/3}a_1^{1/3}} 
= \ta_0 + s^2 + \Order{\eps^{1/3} s^3}, 
\qquad
\ta_0 = \frac{a_0}{\eps^{2/3}a_1^{1/3}} \leqs \frac{\gamma_0}{a_1^{1/3}}.
\end{equation}
\eqref{sdet3:2} is a perturbation of order $\eps^{1/3}$ of the Bernoulli
equation 
\begin{equation}
\label{sdet3:4}
\dtot zs = \ta(s) z - z^3, 
\qquad \ta(s) = \ta_0 + s^2.
\end{equation}
Using Gronwall's inequality, one easily shows that on an $s$-time scale of
order $1$, the solution of \eqref{sdet3:2} differs by $\Order{\eps^{1/3}}$
from the solution of \eqref{sdet3:4}, which is
\begin{equation}
\label{sdet3:5}
z_s = \frac{z_{s_0} \e^{\talpha(s,s_0)}} 
           {\biggpar{1+2z_{s_0}^2 \displaystyle\int_{s_0}^s 
            \e^{2\talpha(u,s_0)}\6u}^{1/2}}, 
\qquad
\text{where $\talpha(s,s_0) = \int_{s_0}^s \ta(u)\6u$.}
\end{equation}
This function is bounded away from zero, and remains of order one for $s$ of
order one, which shows that $x_t\asymp\eps^{1/3}$ on $[t_0,t_1]$.

Since $x_t \asymp \eps^{1/3}$ and $x^\star(t) \asymp \sqrt{a_0} \vee
\abs{t}$, $x_t$ and $x^\star(t)$ necessarily cross at some time
$\tilde t \asymp \eps^{1/3}$. 
\end{proof}

Note that the above proposition also implies bounds on $y_t$, namely,
$y_t = \Order{\eps^{1/3}}$ for $t_0 \leqs t \leqs t_1$, and
there exist constants $\tilde c_+ > \tilde c_- >0$ such that
\begin{equation}
\label{sdet3:5a}
y_t \asymp
\begin{cases}
\phantom{-}
\eps^{1/3} & \text{for $t_0 \leqs t \leqs \tilde c_-\eps^{1/3}$,} \\
\phantom{-}
0 & \text{for $t = \tilde t $,} \\
- \eps^{1/3} & \text{for $\tilde c_+\eps^{1/3} \leqs t \leqs t_1$.} \\
\end{cases}
\end{equation}

Gathering the results for $a_0 \geqs \gamma_0 \eps^{2/3}$ and $a_0
<\gamma_0 \eps^{2/3}$, we see that there exists a time $t_1\asymp
(\sqrt{a_0}\vee\eps^{1/3})$ such that $y_{t_1}\asymp -\eps/t_1^2$. By
enlarging $c_0$ if necessary, we may assume that $t_1 = c_0
(\sqrt{a_0}\vee\eps^{1/3})$. 

\begin{prop}
\label{prop_sdet4}
On the interval $[t_1,T]$, 
\begin{equation}
\label{Sdet8}
y_t \asymp -\frac{\eps}{t^2}.
\end{equation}
\end{prop}
\begin{proof}
The proof is similar to the one of Proposition \ref{prop_sdet1}. 
\end{proof}

Note that the previous result implies $x_t \asymp x^\star(t) \asymp
t$, provided $c_0$ is large enough.

So far, we have proved that for $t\in[-T,T]$, $x_t$ tracks
$x^\star(t)$ at a distance of order 
\begin{equation}
\label{Sdet9}
\frac\eps{t^2} \wedge \frac{\eps}{a_0} \wedge \eps^{1/3},
\end{equation}
and that the two curves cross at a time $\tilde t$ satisfying $\tilde
t-t^\star\asymp(\eps/a_0)\wedge\eps^{1/3}$. Let us now examine the
behaviour of the linearization 
\begin{equation}
\label{Sdet10}
\ba(t) = \sdpar fx(x_t,t),
\end{equation}
which will determine the behaviour of orbits starting close to the
particular solution $x_t$.

\begin{prop}
\label{prop_sdet5}
For all $t\in[-T,T]$ and all $a_0=\orderone{\eps}$, 
\begin{equation}
\label{Sdet11}
\ba(t) \asymp - (t^2\vee a_0\vee \eps^{2/3}).
\end{equation}
\end{prop}
\begin{proof}
By Taylor's formula we get
\begin{align}
\label{sdet5:1}
\nonumber
\ba(t) &= \sdpar fx(x^\star(t)+y_t,t) \\
&= a^\star(t) - \Bigbrak{6+\bigOrder{x^\star(t)+y_t}}
\Bigpar{x^\star(t)+\frac12y_t}y_t.
\end{align}
Consider first the case $a_0\geqs\gamma_0\eps^{2/3}$ for $\gamma_0$ large
enough. Equation \eqref{sdet3a} implies that $a^\star(t)\leqs -c_-(a_0\vee
t^2)$ for a constant $c_->0$. On the other hand, \eqref{Sdet9} shows that
the second term on the right-hand side of~\eqref{sdet5:1} is bounded
in absolute value by $c_+\eps/a_0$ for a constant $c_+>0$. Thus if
$\gamma_0>(c_+/c_-)^{2/3}$ we obtain that $\ba(t)\asymp
a^\star(t)\asymp -(a_0\vee t^2)$.  

We consider next the case $a_0<\gamma_0\eps^{2/3}$. For $\abs{t}\geqs
c_0\eps^{1/3}$, the above argument can be repeated. The non-trivial
case occurs for $\abs{t}<c_0\eps^{1/3}$. By rescaling variables as in
Proposition \ref{prop_sdet3}, we obtain that  
\begin{equation}
\label{sdet5:2}
\ba(t) = \eps^{2/3} a_1^{1/3} \bigbrak{\ta(s,\eps) - 3z_s^2 +
\Order{\eps^{1/3}}}.
\end{equation}
We have to show that $\ba(t) \asymp -\eps^{2/3}$ which is equivalent
to $\ta(s,\eps) - 3z_s^2 \asymp -1$ for $s$ of order one. The lower
bound is trivial as $\ta(s,\eps) \geqs 0$ and $z_s \asymp 1$. In order
to show the upper bound, first note that for $t\leqs 0$, we have
$x_t>x^\star(t)$ which implies $\ba(t)<a^\star(t)$. Therefore, it is
sufficient to consider $s\geqs0$. Taking into account the
expression~\eqref{sdet3:5}, we find that showing the upper bound
amounts to showing that  
\begin{equation}
\label{sdet5:3}
\Bigpar{\text{\it const}{}+\frac{\ta(s)}{z_{s_0}^2}} \e^{-2\talpha(s,s_0)} +
2\ta(s)\int_{s_0}^s\e^{-2\talpha(s,u)}\6u < 3.
\end{equation}
Since $\abs{s_0}$ is proportional to $c_0$, choosing {\it a priori} a
large enough $c_0$ also makes $\abs{s_0}$ large. Thus it is in fact
sufficient to verify that  
\begin{equation}
\label{sdet5:4}
2\ta(s)\int_{-\infty}^s\e^{-2\talpha(s,u)}\6u < 3
\end{equation}
for all $s\geqs 0$. Optimizing the left-hand side with respect to
$\ta_0\geqs0$ and $s$ shows that we may assume $\ta_0=0$ and
that~\eqref{sdet5:4} holds.
\end{proof}


\subsection{The random motion near the stable equilibrium branches}
\label{ssec_symnear}
 
We now consider the SDE
\begin{equation}
\label{snear1}
\6x_t = \frac1\eps f(x_t,t) \6t + \frac\sigma{\sqrt{\eps}} \6W_t,
\qquad x_{-T} = x_0,
\end{equation} 
on the time interval $[-T,T]$, where we assume $x_0-x^\star(-T) \asymp
\eps$.
In order to compare the solution $x_t$ with the solution $\xdet_t$ of
the corresponding deterministic equation~\eqref{sdet1}, we introduce
the difference $y_t = x_t - \xdet_t$, which satisfies the SDE 
\begin{equation}
\label{snear2}
\6y_t = \frac1\eps \bigbrak{\ba(t)y_t + \bb(y_t,t)}\6t +
\frac\sigma{\sqrt{\eps}} \6W_t,
\qquad y_{-T} = 0,
\end{equation}
where $\ba(t)$ is the linearization~\eqref{Sdet10} of $f$ along
$\xdet_t$, and Taylor's formula yields the relations
\begin{equation}
\label{snear3}
\begin{split}
\bb(y,t) &= -\bigbrak{1+\Order{\xdet_t+\abs{y}}+\Order{t}}
(3\xdet_t+y) y^2, \\ 
\abs{\bb(y,t)} &\leqs{} M (\xdet_t+\abs{y}) y^2
\end{split}
\end{equation}
whenever $\abs{t}\leqs T$ and $\xdet_t+\abs{y}\leqs d$, where $M$ is a
positive constant. Let us first consider the linearization of
\eqref{snear2}, namely
\begin{equation}
\label{snear4}
\6y^0_t = \frac1\eps \ba(t) y^0_t \6t +
\frac\sigma{\sqrt{\eps}} \6W_t,
\qquad y^0_{-T} = 0.
\end{equation}
The random variable $y^0_t$ is Gaussian with expectation zero and
variance
\begin{equation}
\label{snear5}
v(t) = \frac{\sigma^2}{\eps} \int_{-T}^t \e^{2\balpha(t,s)/\eps} \6s, 
\qquad \text{where\ }\balpha(t,s) = \int_s^t \ba(u) \6u.
\end{equation}
Lemma \ref{lem_sdet1} and Proposition \ref{prop_sdet5} imply that 
\begin{equation}
\label{snear6}
\z(t) \defby \frac 1{2\abs{\ba(-T)}} \e^{2\balpha(t,-T)/\eps} +
\frac1{\eps} \int_{-T}^t \e^{2\balpha(t,s)/\eps} \6s
\asymp \frac1{t^2\vee a_0\vee\eps^{2/3}}.
\end{equation}
Thus $v(t)$ is of order $\sigma^2/(t^2\vee a_0\vee\eps^{2/3})$, except for
$t$ very close to $-T$. We now show that $y^0_t$ is likely to remain in a
strip of width proportional to $\sqrt{\z(t)}$. 

\begin{prop}
\label{prop_snear1}
For $-T\leqs t\leqs T$, 
\begin{equation}
\label{snear7}
\Bigprobin{-T,0}{\sup_{-T\leqs s\leqs t}
\frac{\abs{y^0_s}}{\sqrt{\z(s)}}\geqs h}
\leqs C(t,\eps) \exp \Bigset{-\frac12 \frac{h^2}{\sigma^2}
\bigpar{1-\Order{\eps}}},
\end{equation}
where
\begin{equation}
\label{snear8}
C(t,\eps) = \frac{\abs{\balpha(t,-T)}}{\eps^2} + 2.
\end{equation}
\end{prop}
\begin{proof}
Let $-T=u_0<u_1<\dots<u_K=t$, with some $K>0$, be a partition of $[-T,t]$.
In \cite[Lemma 3.2]{BG}, we show that the probability \eqref{snear7} is
bounded above by 
\begin{equation}
\label{snear1:1}
2\sum_{k=1}^K P_k, 
\qquad\text{where\ }
P_k = \exp\Bigset{-\frac12 \frac{h^2}{\sigma^2} \inf_{u_{k-1}\leqs s\leqs
u_k} \frac{\z(s)}{\z(u_k)} \e^{2\balpha(u_k,s)/\eps}}.
\end{equation}
Now we choose the partition by requiring that 
\begin{equation}
\label{snear1:2}
\balpha(u_{k},u_{k-1}) = -2\eps^2 
\qquad \text{for\ }
1\leqs k < K = \biggintpartplus{\frac{\abs{\balpha(t,-T)}}{2\eps^2}}. 
\end{equation}
Since $\ba(s)< 0$, we have $\z'(s)=[2\ba(s)\z(s)+1]/\eps \leqs1/\eps$, and
thus
\begin{equation}
\label{snear1:3}
\inf_{u_{k-1}\leqs s\leqs u_k} \frac{\z(s)}{\z(u_k)} 
\geqs \frac1{\z(u_k)} \inf_{u_{k-1}\leqs s\leqs u_k} 
\Bigbrak{\z(u_k) - \frac{u_k-s}{\eps}} 
= 1 - \frac{u_k-u_{k-1}}{\eps\z(u_k)}.
\end{equation}
If $k$ is such that $\abs{u_k}\geqs\sqrt{a_0}\vee\eps^{1/3}$,
then by \eqref{snear1:2} and \eqref{Sdet11}, there is a constant $c_-$ such
that
\begin{equation}
\label{snear1:4}
2\eps^2 \geqs c_-\int_{u_{k-1}}^{u_k}s^2\6s 
\geqs \frac{c_-}{6} \mskip1.5mu u_k^2 (u_k-u_{k-1}),
\end{equation}
and hence by \eqref{snear6} (choosing the same $c_-$ for brevity of
notation)
\begin{equation}
\label{snear1:5}
\frac{u_k-u_{k-1}}{\z(u_k)} 
\leqs \frac{12\eps^2}{c_- u_k^2} \frac{u_k^2}{c_-} 
= \Order{\eps^2}.
\end{equation}
For all other $k$, we have 
\begin{equation}
\label{snear1:6}
2\eps^2 \geqs c_-(a_0\vee\eps^{2/3})(u_k-u_{k-1}) 
\qquad\Rightarrow\qquad
\frac{u_k-u_{k-1}}{\z(u_k)} \leqs \frac{2\eps^2}{c_-^2}.
\end{equation}
In both cases, we find 
\begin{equation}
\label{snear1:7}
P_k \leqs \exp\Bigset{-\frac12 \frac{h^2}{\sigma^2}
\bigpar{1-\Order{\eps}}},
\end{equation}
which leads to the result, using the definition of $K$.
\end{proof}

Let us now compare solutions of the linear equation \eqref{snear4} and the
nonlinear equation \eqref{snear2}. We introduce the events
\begin{align}
\label{snear9a}
\Omega_t(h) &= \bigsetsuch{\w}{\abs{y_s}\leqs h\sqrt{\z(s)}\;\;\forall
s\in[-T,t]} \\
\label{snear9b}
\Omega^0_t(h) &= \bigsetsuch{\w}{\abs{y^0_s}\leqs h\sqrt{\z(s)}\;\;\forall
s\in[-T,t]}. 
\end{align}

\begin{notation}
\label{n_snear}
For two events $\Omega_1$ and $\Omega_2$, we write $\Omega_1\subas\Omega_2$
if\/ $\fP$-almost all $\w\in\Omega_1$ belong to $\Omega_2$.
\end{notation}

\begin{prop}
\label{prop_snear2}
There exists a constant $\varrho$, depending only on $f$, such that
\begin{itemiz}
\item	if $-T\leqs t\leqs -(\sqrt{a_0}\vee\eps^{1/3})$ and $h<t^2/\varrho$,
then 
\begin{equation}
\label{snear10}
\Omega^0_t(h) \subas \Omega_t\Bigpar{\Bigbrak{1+\varrho\frac h{t^2}}h};
\end{equation}
\item	if $-(\sqrt{a_0}\vee\eps^{1/3})\leqs t\leqs T$ and
$h<(a_0\vee\eps^{2/3})/\varrho$,
then 
\begin{equation}
\label{snear11}
\Omega^0_t(h) \subas \Omega_t\Bigpar{\Bigbrak{1+\varrho\frac
h{a_0\vee\eps^{2/3}}}h}.
\end{equation}
\end{itemiz}
\end{prop}
\begin{proof}
The proof is based on the fact that the variable $z_s = y_s - y^0_s$
satisfies the relation
\begin{equation}
\label{snear2:1}
z_s = \frac1\eps \int_{-T}^s \e^{\balpha(s,u)/\eps} \bb(y_u,u) \6u.
\end{equation}
Consider first the case $-T\leqs t\leqs -(\sqrt{a_0}\vee\eps^{1/3})$. Let
$\varrho>0$ be a constant to be chosen later, and set $\delta=\varrho
h/t^2<1$. We define the first exit time
\begin{equation}
\label{snear2:2}
\tau = \inf\bigsetsuch{s\in[-T,t]}{\abs{z_s}\geqs\delta h\sqrt{\z(s)}}
\in [-T,t]\cup \{\infty\}. 
\end{equation}
Pick any $\w\in A\defby \Omega^0_t(h)\cap\setsuch{\omega}{\tau(\w)<\infty}$
and
$s\in[-T,\tau(\w)]$. Then we have 
\begin{equation}
\label{snear2:3}
\abs{y^0_u(\w)}\leqs h\sqrt{\z(u)}, \qquad
\abs{y_u(\w)}\leqs (1+\delta)h\sqrt{\z(u)} < 2h\sqrt{\z(u)}
\end{equation}
for all $u\in[-T,s]$. From \eqref{s6} and \eqref{snear6}, we obtain the
existence of a constant $c_+>0$ such that 
\begin{equation}
\label{snear2:4}
\xdet_u \leqs c_+\abs{u}, \qquad
\abs{y_u(\w)} < 2h \frac{\sqrt{c_+}}{\abs u}
\end{equation}
for these $u$. Hence, by \eqref{snear3} we get the estimate
\begin{equation}
\label{snear2:5}
\abs{\bb(y_u,u)} 
< M \Bigpar{c_+\abs{u}+2h\frac{\sqrt{c_+}}{\abs{u}}} 4h^2\frac{c_+}{u^2}
\leqs 4M \frac{h^2c_+^2}{\abs{s}} \Bigpar{1+\frac{2h}{\sqrt{c_+}s^2}} 
\end{equation}
and thus, by \eqref{snear2:1} and Lemma \ref{lem_sdet1}, 
\begin{equation}
\label{snear2:6}
\abs{z_s} 
< 4M \frac{h^2c_+^2}{\abs{s}} \Bigpar{1+\frac{2h}{\sqrt{c_+}s^2}} 
\frac1\eps \int_{-T}^s \e^{\balpha(s,u)/\eps}\6u 
\leqs 4M \frac{h^2c_+^3}{\abs{s}^3} \Bigpar{1+\frac{2h}{\sqrt{c_+}s^2}},
\end{equation}
where we use again the same $c_+$ for brevity of notation. Using
\eqref{snear6} once again, we arrive at the bound 
\begin{equation}
\label{snear2:7}
\frac{\abs{z_s}}{h\sqrt{\z(s)}} < 4M\frac{c_+^3}{\sqrt{c_-}} \frac
h{s^2} \Bigpar{1+\frac{2}{\sqrt{c_+}}\frac{h}{s^2}}.
\end{equation}
Now we choose 
\begin{equation}
\label{snear2:8}
\varrho = \frac{2}{\sqrt{c_+}} \vee 8M \frac{c_+^3}{\sqrt{c_-}},
\end{equation}
which implies 
\begin{equation}
\label{snear2:9}
\frac{\abs{z_s}}{h\sqrt{\z(s)}} < \frac\varrho2 \frac h{s^2} 
\Bigpar{1+\varrho\frac h{s^2}} \leqs \frac\delta2(1+\delta)<\delta
\end{equation}
for all $s\in[-T,\tau(\omega)]$, by the definition of $\delta$. Hence
$\abs{z_{\tau(\w)}}<\delta h\sqrt{\z(\tau(\w))}$ for almost all $\w\in
A$. Since we have $\abs{z_{\tau(\w)}}=\delta h\sqrt{\z(\tau(\w))}$ whenever
$\tau(\omega)<\infty$, we conclude that $\fP(A)=0$, and thus
$\tau(\w)=\infty$ for almost all $\w\in\Omega^0_t(h)$, which implies that
$\abs{y_s(\w)}\leqs(1+\delta)h\sqrt{\z(s)}$ for $-T\leqs s\leqs t$ and these
$\w$. This completes the proof of \eqref{snear10}. 

The proof of \eqref{snear11} is almost the same. In the case $-(\sqrt{a_0}
\vee \eps^{1/3}) \leqs t \leqs T$, we take $\delta=\varrho
h/(a_0\vee\eps^{2/3})$. The estimate \eqref{snear2:4} has 
to be replaced by  
\begin{equation}
\label{snear2:10}
\xdet_u \leqs c_+(\abs{u}\vee\sqrt{a_0}\vee\eps^{1/3}), \qquad
\abs{y_u(\w)} \leqs 
2h \frac{\sqrt{c_+}}{\abs u\vee\sqrt{a_0}\vee\eps^{1/3}}, 
\end{equation}
and thus we get, instead of \eqref{snear2:6}, the bound
\begin{equation}
\label{snear2:11}
\abs{z_s} 
\leqs 4M \frac{h^2c_+^3}{(\sqrt{a_0}\vee\eps^{1/3}) 
(s^2\vee a_0\vee\eps^{2/3})} 
\biggpar{1+\frac{2h}{\sqrt{c_+}(a_0\vee\eps^{2/3})}}. 
\end{equation}
The remainder of the proof is similar.
\end{proof}

Now, the preceding two propositions immediately imply
Theorem~\ref{thm_snearxdet}, as Proposition~\ref{prop_snear1} shows
the desired behaviour for the approximation by a Gaussian process and
Proposition~\ref{prop_snear2} allows to extend this result to the
original process.


\subsection{The transition regime}
\label{ssec_symtrans}

We consider now the regime of $\sigma$ sufficiently large to allow for
transitions from one stable equilibrium branch to the other. Here $\xdet_t$
is the solution of the deterministic equation \eqref{sdet1} with the same
initial condition $\xdet_{-T}$ as in the previous sections, which tracks
$x^\star(t)$ at distance at most $\Order{\eps^{1/3}}$. $x_t$ denotes a
general solution of the SDE \eqref{s1}. Our aim is to establish an upper
bound for the probability of {\em not\/} reaching the axis $x=0$, which, by
using symmetry, will allow us to estimate the transition probability.
Let $\delta>0$ be the constant defined in \eqref{sres18}, i.e.\ by
\begin{equation}
\label{strans1}
x\sdpar{f}{xx}(x,t)\leqs 0
\qquad
\text{for $\abs{x}\leqs\delta$ and $\abs{t}\leqs T$.}
\end{equation}

The basic ingredient of our estimate is the following comparison lemma
which allows us to linearize the stochastic differential equations
under consideration and, therefore, to investigate Gaussian
approximations to our processes. The lemma gives conditions under which
relations between initial conditions carry over to the sample
paths. 

\begin{lemma}
\label{lem_strans}
Fix some initial time $t_0\in[-T,T]$. 
We consider the following processes on $[t_0,T]$:
\begin{itemiz}
\item	the solution $\xdet_t$ of the deterministic differential
equation~\eqref{sdet1} with initial condition $\xdet_{t_0} \in [0,\delta]$;
\item	the solution $x_t$ of the SDE~\eqref{snear1} with initial condition
$x_{t_0} \in [\xdet_{t_0},\delta]$;
\item	the difference $y_t=x_t-\xdet_t$, which satisfies $y_{t_0} =
x_{t_0}-\xdet_{t_0} \geqs 0$;
\item	 the solution $y^0_t$of the linear SDE
\begin{equation}
\label{strans2}
\6y^0_t = \frac1\eps \ta(t)  y^0_t \6t 
+ \frac{\sigma}{\sqrt{\eps}} \6W_t, 
\qquad
\text{where $\ta(t) = \sdpar{f}{x}(\xdet_t,t)$}
\end{equation}
with initial condition $y^0_{t_0} \in[y_{t_0}, \delta-\xdet_{t_0}]$.
\end{itemiz}
If\/ $0 \leqs y^0_s \leqs \delta-\xdet_s$ for all $s\in[t_0,t]$, then
$y_s\leqs y^0_s$ for those $s$. Similarly, if\/ $0 \leqs y_s \leqs
\delta-\xdet_s$ for all $s\in[t_0,t]$, then $y^0_s\geqs y_s$ for those
$s$. The result remains true when $t$ is replaced by a stopping time.
\end{lemma}
\begin{proof}
The hypothesis \eqref{strans1} implies that for all
$y\in[0,\delta-\xdet_s]$, 
\begin{equation}
\label{lstrans:1}
f(\xdet_s+y, s) \leqs f(\xdet_s,s) + \ta(s) y.
\end{equation}
Let $\tau = \inf\setsuch{s\in[t_0,t]}{y_s\not\in[0,\delta-\xdet_s]}
\in [t_0,t] \cup\{\infty\}$. For $t_0\leqs s\leqs\tau$, the variable
$z_s=y_s-y^0_s$ satisfies 
\begin{align}
\label{lstrans:2}
\nonumber
z_s &= z_{t_0} + \frac1\eps 
\int_{t_0}^s \bigbrak{f(\xdet_u+y_u,u) - f(\xdet_u,u) - \ta(u)y^0_u}\6u\\
&{}\leqs z_{t_0} + \frac1\eps \int_{t_0}^s \ta(u) z_u \6u, 
\qquad
z_{t_0}\leqs 0.
\end{align}
Applying Gronwall's inequality, we obtain 
\begin{equation}
\label{lstrans:3}
z_s \leqs z_{t_0} \e^{\talpha(s,t_0)/\eps} \leqs 0
\qquad
\forall s\in[t_0,\tau\wedge t],
\end{equation}
where $\talpha(s,t_0) = \int_{t_0}^s \ta(u)\6u$.
This proves the result for $t_0\leqs s\leqs \tau\wedge t$. Now if $y_s$ is
negative, the result is trivially satisfied, and if $y_s$ becomes positive
again, the above argument can be repeated. Note that $y_\tau \leqs
\delta - \xdet_\tau$ is immediate. This proves the first assertion,
and the second assertion can be proved directly, without use of $\tau$.   
\end{proof}

We will now proceed as follows. Let $\z(t)$ be the function defined in
\eqref{snear6}, and let $h$ be such that $\xdet_s+h\sqrt{\z(s)}<\delta$ for
all $s\in[t_0,t]$. Given $x_0\in(0,\delta)$, we can write 
\begin{align}
\label{strans3}
\nonumber
\bigprobin{t_0,x_0}{x_s>0 \;\; \forall s\in[t_0,t]} 
\leqs{}& \Bigprobin{t_0,x_0}{\sup_{t_0\leqs s\leqs t}
\frac{x_s-\xdet_s}{\sqrt{\z(s)}} > h} \\ 
&{}+
\bigprobin{t_0,x_0}{0<x_s\leqs\xdet_s+h\sqrt{\z(s)} \;\; \forall
s\in[t_0,t]}.
\end{align}
We will estimate these two terms separately. 
The first event is similar to the event we have examined in the previous
subsection, but we need here an estimate valid for all times, even when
$\sigma$ is not very small, whereas the previous result is only useful
for $\sigma \leqs t^2 \vee a_0 \vee \eps^{2/3}$. We will show the
following.

\begin{prop}
\label{prop_strans1}
Assume $0\leqs x_0\leqs \xdet_{t_0}+\frac12 h\sqrt{\z(t_0)}$. Then 
\begin{equation}
\label{strans4}
\Bigprobin{t_0,x_0}{\sup_{t_0\leqs s\leqs t}
\frac{x_s-\xdet_s}{\sqrt{\z(s)}} > h} 
\leqs \frac52 \Bigpar{\frac{\abs{\balpha(t,t_0)}}\eps + 1} 
\e^{-\bar\kappa h^2/\sigma^2},
\end{equation}
where $\bar\kappa$ is a positive constant and $\balpha(t,t_0) = \int_{t_0}^t
\ba(s)\6s$. 
\end{prop}
\goodbreak
\begin{proof}\hfill
\begin{enum}
\item  We define a partition $t_0=u_0<u_1<\dots<u_K=t$ of the interval
$[t_0,t]$ by requiring
\begin{equation}
\label{strans1:1}
\abs{\balpha(u_k,u_{k-1})} = \eps 
\qquad
\text{for\ }
1\leqs k < K = \biggintpartplus{\frac{\abs{\balpha(t,t_0)}}{\eps}}. 
\end{equation}
Note that similar arguments as in the proof of Proposition
\ref{prop_snear1} yield 
\begin{equation}
\label{strans1:2}
\frac{u_{k+1}-u_k}{\z(u_k)} = \Order{\eps}
\qquad\text{for all $k$.}
\end{equation}
Now let $\rho_k=\frac12 h\sqrt{\z(u_k)}$ and $y_s=x_s-\xdet_s$ as
usual. Define 
\begin{align}
\label{strans1:3}
\nonumber
Q_k = \sup_{y_k\leqs\rho_k} \biggbrak{
&\Bigprobin{u_k,y_k}{\sup_{u_k\leqs s\leqs u_{k+1}}
\frac{y_s}{\sqrt{\z(s)}}>h} \\
&{}+ \Bigprobin{u_k,y_k}{\sup_{u_k\leqs s\leqs u_{k+1}}
\frac{y_s}{\sqrt{\z(s)}}\leqs h, \;
y_{u_{k+1}} > \rho_{k+1}}
}, 
\end{align}
for $0 \leqs k < K-1$, and
\begin{equation}
\label{strans1:3b}
Q_{K-1} = \sup_{y_{K-1}\leqs\rho_{K-1}} 
\Bigprobin{u_{K-1},y_{K-1}}{\sup_{u_{K-1}\leqs s\leqs u_{K}}
\frac{y_s}{\sqrt{\z(s)}}>h}.
\end{equation}
Then
\begin{align}
\nonumber
\vphantom{\sum_{k=0}^{K-1}}
&\Bigprobin{t_0,x_0}{\sup_{t_0\leqs s\leqs t}
\frac{x_s-\xdet_s}{\sqrt{\z(s)}} > h} \\
\vphantom{\sum_{k=0}^{K-1}}\nonumber
&\qquad\qquad\leqs{} \Bigprobin{t_0,y_{t_0}}{\sup_{t_0\leqs s\leqs u_1}
\frac{y_s}{\sqrt{\z(s)}} > h} 
+ \Bigprobin{t_0,y_{t_0}}{\sup_{t_0\leqs s\leqs u_1}
\frac{y_s}{\sqrt{\z(s)}} \leqs h, \; y_{u_1}>\rho_1} \\
\vphantom{\sum_{k=0}^{K-1}}\nonumber
&\qquad\qquad\phantom{{}\leqs{}} 
{}+ \Bigexpecin{t_0,y_{t_0}}{\indexfct{y_{u_1}\leqs\rho_1} 
\Bigprobin{u_1,y_{u_1}}{\sup_{u_1\leqs s\leqs t}
\frac{y_s}{\sqrt{\z(s)}} > h}} \\
\label{strans1:4}
\vphantom{\sum_{k=0}^{K-1}}
&\qquad\qquad\leqs{} \dots \leqs \sum_{k=0}^{K-1} Q_k. 
\end{align}

\item	In order to estimate $Q_k$, we introduce the stochastic process
$(y^{(k)}_s)_{s\in[u_k,u_{k+1}]}$ defined by 
\begin{equation}
\label{strans1:5}
y^{(k)}_s = \rho_k \e^{\balpha(s,u_k)/\eps} +
\frac{\sigma}{\sqrt\eps}\int_{u_k}^s \e^{\balpha(s,u)/\eps}\6W^{(k)}_u, 
\end{equation}
where $(W^{(k)}_u)_{u\in[u_k,u_{k+1}]}$ is the Brownian motion
$W^{(k)}_u = W_u-W_{u_k}$. Note that $y^{(k)}$ is the solution of the
SDE~\eqref{strans2} with initial condition \smash{$y^{(k)}_{u_k} =
\rho_k$} at time $u_k$. We define  
the stopping times 
\begin{equation}
\label{strans1:6}
\begin{split}
\tau^0 &= \inf\bigsetsuch{s\in[u_k,u_{k+1}]}{y^{(k)}_s=0} \\
\tau^+ &= \inf\bigsetsuch{s\in[u_k,u_{k+1}]}{y^{(k)}_s=h\sqrt{\z(s)}\,} 
\end{split}
\end{equation}
describing the time when $y^{(k)}_s$ either reaches the $t$-axis or the
upper
boundary $h \sqrt{\z(s)}$. Now, Lemma~\ref{lem_strans} implies that if
$y_{u_k}\leqs\rho_k$, then $y_s\leqs y^{(k)}_s$ for $u_k\leqs
s\leqs\tau^0\wedge\tau^+$. This shows that 
\begin{equation}
\label{strans1:7}
Q_k \leqs \bigprobin{u_k,\rho_k}{\tau^0<u_{k+1}} +
\bigprobin{u_k,\rho_k}{\tau^+<u_{k+1}} + 
\bigprobin{u_k,\rho_k}{y^{(k)}_{u_{k+1}} > \rho_{k+1}}
\end{equation}
for $0 \leqs k < K-1$, and
\begin{equation}
\label{strans1:7b}
Q_{K-1} \leqs \bigprobin{u_{K-1},\rho_{K-1}}{\tau^0<u_{K}} +
\bigprobin{u_{K-1},\rho_{K-1}}{\tau^+<u_{K}}.
\end{equation}
Each of these terms depends only on $y^{(k)}$, and can be easily estimated. 
Let 
\begin{equation}
\label{strans1:8}
v^{(k)}_{u_{k+1}} = \frac{\sigma^2}{\eps}
\int_{u_k}^{u_{k+1}}\e^{2\balpha(u_{k+1},u)/\eps}\6u
\end{equation}
denote the variance of $y^{(k)}_{u_{k+1}}$. Then by symmetry (as in
\eqref{sres17}), we have 
\begin{align}
\label{strans1:9}
\nonumber
\bigprobin{u_k,\rho_k}{\tau^0<u_{k+1}} 
&= 2\mskip1.5mu \bigprobin{u_k,\rho_k}{y^{(k)}_{u_{k+1}}<0} \\
\nonumber
&= \frac2{\sqrt{2\pi}}
\int_{-\infty}^{-{\rho_k\e^{\balpha(u_{k+1},u_k)/\eps}} 
\bigpar{v^{(k)}_{u_{k+1}}}^{-1/2}} \e^{-z^2/2} \6z \\
\nonumber
&\leqs \exp\biggset{-\frac12 
\frac{\rho_k^2\e^{2\balpha(u_{k+1},u_k)/\eps}}{v^{(k)}_{u_{k+1}}}}\\
&\leqs \exp\biggset{-\frac18 \e^{-2} \frac{h^2}{\sigma^2} 
\frac{\z(u_k)}{v^{(k)}_{u_{k+1}}/\sigma^2}}. 
\end{align}
The second term on the right-hand side of~\eqref{strans1:7}
or~\eqref{strans1:7b}, respectively, can be estimated using the
symmetry (in distribution) of~\eqref{strans1:5} under the map
$\sigma\mapsto-\sigma$:
\begin{align}
\label{strans1:10} \nonumber
\bigprobin{u_k,\rho_k}{\tau^+<u_{k+1}} 
& = \bigprobin{u_k,\rho_k}{\exists s\in[u_k,u_{k+1}] \colon
y^{(k)}_s \geqs h\sqrt{\z(s)}} \\
\nonumber
& \leqs \bigprobin{u_k,\rho_k}{\exists s\in[u_k,u_{k+1}] \colon
y^{(k)}_s \leqs h \sqrt{\z(u_k)} \e^{\balpha(s,u_k)/\eps} - h
\sqrt{\z(s)}} \\
& \leqs \bigprobin{u_k,\rho_k}{\tau^0<u_{k+1}}. 
\end{align}
In order to estimate the third term on the right-hand side
of~\eqref{strans1:7}, we will use the fact that for $k<K-1$
\begin{align}
\label{strans1:11}
\nonumber
\z(u_{k+1}) &= \z(u_k)\e^{2\balpha(u_{k+1},u_k)/\eps} 
+ \frac1\eps \int_{u_k}^{u_{k+1}} \e^{2\balpha(u_{k+1},s)/\eps} \6s \\
&\geqs \z(u_k)\e^{-2} + \mskip1.5mu
\frac{1-\e^{-2}}2 \inf_{u_k\leqs s\leqs
u_{k+1}} \frac1{\abs{\ba(s)}}.
\end{align}
Proposition \ref{prop_sdet5} and \eqref{snear6} thus yield the existence of
a constant $c_->0$ such that
\begin{equation}
\label{strans1:12}
\frac{\z(u_{k+1})}{\z(u_k)} 
\geqs \e^{-2} + \mskip1.5mu
\frac{1-\e^{-2}}{2c_-^2}.
\end{equation}
This allows us to estimate (for $k<K-1$)
\begin{align}
\label{strans1:13}
\nonumber
\bigprobin{u_k,\rho_k}{y^{(k)}_{u_{k+1}} > \rho_{k+1}} 
&\leqs \frac12 \exp\biggset{-\frac12
\frac{(\rho_{k+1}-\rho_k\e^{\balpha(u_{k+1},u_k)/\eps})^2}{v^{(k)}_{u_{k+1}}}}
\\ 
&= \frac12 \exp\Biggset{-\frac{h^2}8 \frac{\z(u_k)}{v^{(k)}_{u_{k+1}}}
\Biggpar{\sqrt{\frac{\z(u_{k+1})}{\z(u_k)}} - \e^{-1}}^2}.
\end{align}

\item	The estimates \eqref{strans1:9}, \eqref{strans1:10} and
\eqref{strans1:13}, inserted in~\eqref{strans1:7}
and~\eqref{strans1:7b}, imply that 
\begin{equation}
\label{strans1:14}
Q_k \leqs \frac 52 \exp\Bigset{-\kappa_k\frac{h^2}{\sigma^2}}, 
\end{equation}
with 
\begin{equation}
\label{strans1:15}
\kappa_k = \frac18 \frac{\z(u_k)}{v^{(k)}_{u_{k+1}}/\sigma^2} \e^{-2} 
\Biggbrak{1\wedge \Biggpar{\sqrt{1+\frac{\e^{2}-1}{2c_-^2}}-1}^2\,}. 
\end{equation}
By \eqref{strans1:8}, for each $k$, there exists a
$\theta_k\in[\e^{-2},1]$ such that 
\begin{equation}
\label{strans1:16}
\frac{v^{(k)}_{u_{k+1}}}{\sigma^2} 
= \frac{(u_{k+1}-u_k)}\eps\theta_k. 
\end{equation}
Together with \eqref{strans1:2}, this implies that $\kappa_k\asymp1$ for
all $k$, and thus the result follows from \eqref{strans1:4} with
$\bar\kappa=\inf_k\kappa_k$. 
\qed
\end{enum}
\renewcommand{\qed}{}
\end{proof}

We now give an estimate of the second term in \eqref{strans3}. The
Markov property implies that we will obtain an upper bound by starting at 
time $-c_1\sqrt\sigma$.  

\begin{prop}
\label{prop_strans2}
There exist constants $c_1>0$ and $\bar\kappa>0$ such that, if
$c_1^2\sigma \geqs a_0\vee\eps^{2/3}$ and $h>2\sigma$, then
\begin{multline}
\label{strans5}
\qquad\qquad
\bigprobin{-c_1\sqrt\sigma,x_0}{0<x_s\leqs\xdet_s+h\sqrt{\z(s)} \;\; \forall
s\in[-c_1\sqrt\sigma,t_1]} \\
\leqs 2 \exp\biggset{-\bar\kappa \mskip3mu \frac1{\log(h/\sigma)}
\frac{\alpha(t_1,-c_1\sqrt\sigma)}\eps}
\qquad\qquad
\end{multline}
holds with $\alpha(t_1,-c_1\sqrt\sigma) = \int_{-c_1\sqrt\sigma}^t
a(s)\6s$, for  $-c_1\sqrt\sigma \leqs t_1 \leqs c_1\sqrt{\sigma}$ and all
initial conditions $x_0$ satisfying $0\leqs x_0\leqs
\xdet_{-c_1\sqrt\sigma}+h\sqrt{\z(-c_1\sqrt\sigma)}$.
\end{prop}
%
\begin{proof} \hfill
\begin{enum}
\item	Let $\varrho = \varrho(h/\sigma) \geqs 1$ and define a partition
$-c_1\sqrt\sigma=u_0 < \dots < u_K=t_1$ of $[-c_1\sqrt\sigma,t_1]$ by 
\begin{equation}
\label{strans2:1}
\alpha(u_k,u_{k-1}) = \varrho\eps
\qquad \text{for\ }
1\leqs k<K =
\biggintpartplus{\frac{\alpha(t_1,-c_1\sqrt\sigma)}{\varrho\eps}}.
\end{equation}
We would like to control the probability of not reaching the $t$-axis during
the time interval $[u_k,u_{k+1}]$. Let 
\begin{equation}
\label{strans2:2}
Q_k = \sup_{0<x_k\leqs \xdet_{u_k}+h\sqrt{\z(u_k)}} 
\bigprobin{u_k,x_k}{0<x_s\leqs\xdet_s+h\sqrt{\z(s)} \;\; \forall
s\in[u_k,u_{k+1}]}.
\end{equation}
Then the probability on the left-hand side of~\eqref{strans5} is
\begin{align}
\label{strans2:3}
\nonumber
&\bigprobin{-c_1\sqrt\sigma,x_0}{0<x_s\leqs\xdet_s+h\sqrt{\z(s)} \;\;
\forall
s\in[-c_1\sqrt\sigma,t_1]}\\
\nonumber&\qquad
=\Bigexpecin{-c_1\sqrt\sigma,x_0}{\indexfct{0<x_s\leqs\xdet_s+h\sqrt{\z(s)}
\;\; \forall s\in[-c_1\sqrt\sigma,u_{K-1}]} \\
\nonumber
&\qquad\hphantom{{}=\E^{\mskip1.5mu -c_1\sqrt\sigma,x_0}\Bigl\{\null\Bigr.}
\bigprobin{u_{K-1},x_{u_{K-1}}}{0<x_s\leqs\xdet_s+h\sqrt{\z(s)} 
\;\; \forall s\in[u_{K-1},u_K]}} \\
\nonumber
&\qquad
{}\leqs Q_{K-1} \,
\bigprobin{-c_1\sqrt\sigma,x_0}{0<x_s\leqs\xdet_s+h\sqrt{\z(s)} \;\; \forall
s\in[-c_1\sqrt\sigma,u_{K-1}]} \\
&\qquad
{}\leqs \dots \leqs \prod_{k=0}^{K-1} Q_k.
\end{align}
If we manage to estimate each $Q_k$ by a constant less than $1$ (say,
$1/2$), then the probability will be exponentially small in $K$.  In the
sequel, we shall estimate $Q_k$ uniformly in $k=0,\dots K-2$, and bound
$Q_{K-1}$ by $1$, since the last interval of the partition may be too small
to get a good bound. So let $k<K-1$ from now on.

\item	We consider first the case $0<x_k\leqs\xdet_{u_k}$. We define the
process $(x^{(k)}_s)_{u_k\leqs s\leqs u_{k+1}}$ as the solution of the
linearized SDE
\begin{equation}
\label{strans2:4}
\6x^{(k)}_s = a(s) x^{(k)}_s \6s + \frac{\sigma}{\sqrt\eps} \6W^{(k)}_s, 
\qquad x^{(k)}_{u_k} = x_k, 
\end{equation}
where $(W^{(k)}_s)_{s\in[u_k,u_{k+1}]}$ is the Brownian motion
$W^{(k)}_s = W_s - W_{u_k}$. Let $v^{(k)}_{u_{k+1}}$ denote the variance of
$x^{(k)}_{u_{k+1}}$. Then 
\begin{align}
\label{strans2:5}
\nonumber
\e^{-2\alpha(u_{k+1},u_k)/\eps} v^{(k)}_{u_{k+1}} 
&= \frac{\sigma^2}\eps \int_{u_k}^{u_{k+1}} \e^{-2\alpha(u,u_k)/\eps} \6u \\
\nonumber
&\geqs \frac{\sigma^2}2 \inf_{u_k\leqs u\leqs u_{k+1}} \frac1{a(u)}
\bigbrak{1-\e^{-2\alpha(u_{k+1},u_k)/\eps}} \\
&\geqs \frac{1-\e^{-2\varrho}}2 \frac{\sigma^2}{a(u_k)\vee a(u_{k+1})}. 
\end{align}
We can now apply Lemma \ref{lem_strans} in the particular case
$\xdet_s\equiv0$ to show that if $0<x_s\leqs\delta$ for $s\in[u_k,u_{k+1}]$,
then $x^{(k)}_s \geqs x_s$ in the same interval. We thus obtain
\begin{equation}
\label{strans2:6a} 
\bigprobin{u_k,x_k}{0<x_s\leqs\xdet_s+h\sqrt{\z(s)} \;\; \forall
s\in[u_k,u_{k+1}]} 
\leqs \bigprobin{u_k,x_k}{x^{(k)}_s >0 \;\; \forall
s\in[u_k,u_{k+1}]}.
\lower1\jot\hbox{\strut}
\end{equation}
The probability on the right-hand side satisfies
\begin{align}
\nonumber
\bigprobin{u_k,x_k}{x^{(k)}_s >0 \;\; \forall s\in[u_k,u_{k+1}]}
&= 1- 2\mskip1.5mu \bigprobin{u_k,x_k}{x^{(k)}_{u_{k+1}} < 0}\\
\label{strans2:6b}
&= 2\mskip1.5mu \bigprobin{u_k,x_k}{x^{(k)}_{u_{k+1}} \geqs 0} -1,
\lower1\jot\hbox{\strut}
\end{align}
yielding
\begin{align}
\label{strans2:6}
\nonumber
&\bigprobin{u_k,x_k}{0<x_s\leqs\xdet_s+h\sqrt{\z(s)} \;\; \forall
s\in[u_k,u_{k+1}]} \\
\nonumber&\qquad
\leqs \frac2{\sqrt{2\pi}} \int_0^{{x_k\e^{\alpha(u_{k+1},u_k)/\eps}} 
\bigpar{v^{(k)}_{u_{k+1}}}^{-1/2}} \e^{-z^2/2} \6z \\
\nonumber&\qquad
\leqs \frac2{\sqrt{2\pi}}  
\frac{x_k}{\sqrt{\e^{-2\alpha(u_{k+1},u_k)/\eps} v^{(k)}_{u_{k+1}}}} \\
&\qquad
\leqs \frac2{\sqrt{\pi}} \sqrt{\frac1{1-\e^{-2\varrho}}} 
\frac{\sqrt{a(u_k)\vee a(u_{k+1})}}\sigma \,\xdet_{u_k}.
\end{align}
By making $c_1$ small enough, we can guarantee that this bound is smaller
than some imposed constant of order $1$, say $1/2$. This shows that the
length of $[u_k,u_{k+1}]$ has been chosen large enough that the probability
of reaching the $t$-axis during this time interval is appreciable.

\item	We examine now the case $\xdet_{u_k} < x_k < \xdet_{u_k} +
h\sqrt{\z(s)}$. We introduce a time $\tilde u_k\in(u_k,u_{k+1})$, defined by
\begin{equation}
\label{strans2:7}
\alpha(\tilde u_k,u_k) = \frac12\varrho\eps.
\end{equation}
Our strategy will be to show that $x_t$ is likely to cross $\xdet_t$ before
time $\tilde u_k$, which will allow us to use the previous result. 
Proposition \ref{prop_sdet5} implies the existence of a constant $L>0$ such
that 
\begin{equation}
\label{strans2:8}
\frac12\varrho\eps = \int_{u_k}^{\tilde u_k} a(u)\6u 
\leqs L \int_{u_k}^{\tilde u_k} (-\ba(u))\6u
= L \abs{\balpha(\tilde u_k,u_k)}.
\end{equation}
Let $(y^{(k)}_s)_{u_k\leqs s\leqs u_{k+1}}$ be the solution of the
linear SDE
\begin{equation}
\label{strans2:8a}
\6y^{(k)}_s = \ba(s) y^{(k)}_s \6s + \frac{\sigma}{\sqrt\eps} \6W^{(k)}_s, 
\qquad y^{(k)}_{u_k} = y_k = x_k - \xdet_{u_k}, 
\end{equation}
where $(W^{(k)}_s)_{s\in[u_k,u_{k+1}]}$ is again the Brownian motion
$W^{(k)}_s = W_s - W_{u_k}$. The variance of $y^{(k)}_{\tilde u_k}$ is 
\begin{equation}
\label{strans2:9}
\tilde v^{(k)}_{\tilde u_k} = 
\frac{\sigma^2}\eps \int_{u_k}^{\tilde u_k} 
\e^{2\balpha(\tilde u_k,s)/\eps}\6s 
\geqs \frac{1-\e^{-\varrho/L}}2
\frac{\sigma^2}{\abs{\ba(u_k)}\vee\abs{\ba(\tilde u_k)}}. 
\end{equation}
Lemma \ref{lem_strans} shows that if $\xdet_s\leqs x_s\leqs\delta$ on the
interval $[u_k,\tilde u_k]$, then $x_s-\xdet_s \leqs y^{(k)}_s$ on that
interval. If we introduce the stopping time
\begin{equation}
\label{strans2:10}
\tau_k = \inf\bigsetsuch{s\in[u_k,\tilde u_k]}{x_s=\xdet_s}
\in[u_k,u_{k+1}]\cup\{\infty\},  
\end{equation}
then we have
\begin{align}
\nonumber
&\bigprobin{u_k,x_k}{0<x_s\leqs\xdet_s+h\sqrt{\z(s)} \;\; \forall
s\in[u_k,u_{k+1}]} \\
\nonumber&\quad
\leqs \bigprobin{u_k,x_k}{\xdet_s<x_s\leqs\xdet_s+h\sqrt{\z(s)} 
\;\; \forall s\in[u_k,\tilde u_k]} \\
&\quad
\phantom{\leqs{}} + 
\Bigexpecin{u_k,x_k}{\indexfct{\tau_k<\tilde u_k} \mskip1.5mu
\bigprobin{\tau_k,\xdet_{\tau_k}}{0<x_s\leqs\xdet_s+h\sqrt{\z(s)} 
\;\; \forall s\in[\tau_k,u_{k+1}]}}. 
\label{strans2:11}
\end{align}
The second term on the right-hand side can be bounded, as in
\eqref{strans2:6}, by 
\begin{equation}
\label{strans2:11a}
\frac2{\sqrt{\pi}} \sqrt{\frac1{1-\e^{-\varrho}}}
\frac{\sqrt{a(u_k)\vee a(u_{k+1})}}\sigma 
\sup_{u_k\leqs s\leqs \tilde u_k} \xdet_s.
\end{equation}
Using \eqref{strans2:9}, the first term on the right-hand side of
\eqref{strans2:11} can be estimated in the following way:
\begin{align}
\nonumber
&\bigprobin{u_k,x_k}{\xdet_s<x_s\leqs\xdet_s+h\sqrt{\z(s)} 
\;\; \forall s\in[u_k,\tilde u_k]} \\
\nonumber&\qquad
\vrule height 12pt depth 7pt width 0pt
\leqs \bigprobin{u_k,y_k}{y^{(k)}_s > 0  
\;\; \forall s\in[u_k,\tilde u_k]} \\
\nonumber&\qquad
\leqs \frac2{\sqrt{2\pi}} \frac{y_k \e^{\balpha(\tilde u_k,u_k)/\eps}} 
{\bigpar{\tilde v^{(k)}_{\tilde u_k}}^{1/2}} \\
&\qquad
\leqs 
\frac2{\sqrt{\pi}} \sqrt{\frac1{1-\e^{-\varrho/L}}} 
\sqrt{\abs{\ba(u_k)}\vee\abs{\ba(\tilde u_k)}} \,\frac h\sigma 
\sqrt{\z(u_k)} \e^{-\varrho/2L}.
\label{strans2:12}
\end{align}
Using Proposition \ref{prop_sdet5} and \eqref{snear6}, it is easy to show
that the expression $(\abs{\ba(u_k)}\vee\abs{\ba(\tilde u_k)})\z(u_k)$ is
uniformly bounded by a constant independent of $k$ and $\eps$.  The sum of
\eqref{strans2:11a} and of the last term in \eqref{strans2:12} provides an
upper bound for $Q_k$.  

\item	Using the fact that for $\abs{u}\leqs c_1\sqrt\sigma$ and
$c_1^2\sigma \geqs a_0 \vee \eps^{2/3}$, one has
$a(u)=\Order{c_1^2\sigma}$ and $\xdet_u=\Order{c_1\sqrt{\sigma}}$, we arrive
at the bound
\begin{equation}
\label{strans2:13}
Q_k \leqs C \Bigpar{c_1^2 + \frac h\sigma \e^{-\varrho/2L}},
\end{equation}
where the constant $C$ can be chosen independent of $\varrho$ because
$\varrho\geqs 1$ by assumption. Thus if we choose $c_1^2\leqs 1/4C$ and
$\varrho = 2L\log(4Ch/\sigma) \vee 1$, we obtain that $Q_k\leqs 1/2$ for
$k=0,\dots,K-2$. This yields
\begin{equation}
\label{strans2:14}
\prod_{k=0}^{K-1} Q_k \leqs 2\frac1{2^K} 
\leqs 2\exp\Bigset{-(\log2)\mskip1.5mu 
\frac{\alpha(t_1,-c_1\sqrt\sigma)}{\varrho\eps}}, 
\end{equation}
and the result follows from our choice of $\varrho$. 
\qed
\end{enum}
\renewcommand{\qed}{}
\end{proof}

Now the proof of Theorem~\ref{thm_strans} follows from~\eqref{strans3} and
the two preceding propositions, where we use the Markov property to \lq\lq
restart\rq\rq\ at time $-c_1\sqrt\sigma$ before applying
Proposition~\ref{prop_strans2}. 


\subsection{Escape from the saddle}
\label{ssec_symesc}

In this subsection, we investigate the behaviour of the random motion
$x_t$ given by the SDE~\eqref{s1} for $t\geqs t_2 \geqs c_1\sqrt\sigma$,
i.\,e., after the transition regime. We want to show that $x_t$ is
likely to leave a suitably defined neighbourhood of the saddle within
time $\Order{\eps\abs{\log\sigma}/t_2^2}$. The proof of
Theorem~\ref{thm_sescape} is very similar to the proof
of~\cite[Theorem~2.9]{BG}, and for the sake of brevity, we will
refrain from giving all the details. Instead, we will discuss how to
proceed and then focus on those parts which need to be modified.

From now on, we will assume that $t\geqs t_2\geqs c_1\sqrt\sigma$ and that
$\sigma$ is large enough in order to allow for transitions, i.\,e.,
$\sigma\geqs a_0\vee \eps^{2/3}$. We want to estimate the first exit
time $\tau_{\cD(\kappa)}$ of $x_t$ from the set
\begin{equation}
\label{sesc:1}
\cD(\kappa) = \biggsetsuch{(x,t)\in[-\delta,\delta]\times[c_1\sqrt\sigma,T]}
{\frac{f(x,t)}{x} > \kappa a(t)},
\end{equation}
where $\kappa\in(0,1)$ is a constant. Note that the upper boundary $\tx(t)$
of $\cD(\kappa)$ satisfies  $\tx(t) = \sqrt{1-\kappa}\mskip1.5mu
(1-\Order{t}) x^\star(t)$. Our first step towards estimating
$\tau_{\cD(\kappa)}$ is to estimate the first exit time $\tau_\cS$ from a
smaller strip $\cS$, defined by
\begin{equation}
\label{sesc:2}
\cS = \biggsetsuch{(x,t)\in[-\delta,\delta]\times[c_1\sqrt\sigma,T]}
{\abs{x} < \frac h{\sqrt{a(s)}}},
\end{equation}
where we will choose $h$ later. Note that  $h<\text{\it const}\;\sigma$ for
some (small) constant would assure $\cS\subset\cD(\kappa)$. We will not
impose such a restrictive condition on $h$ but replace $\cS$ by
$\cS\cap\cD(\kappa)$ in case $\cS$ is not a subset of $\cD(\kappa)$. The
following proposition gives our estimate on the first exit time from $\cS$.

\begin{prop}
\label{p_sescape}
Let $t_2\geqs c_1\sqrt\sigma$ and $(x_2,t_2)\in\cS$. Then
there exists a constant $L>0$ such that for any $\mu>0$, we have 
\begin{equation}
\label{sesc:3}
\bigprobin{t_2,x_2}{\tau_\cS \geqs t} 
\leqs \Bigpar{\frac h\sigma}^\mu
\exp\biggset{-\frac{\mu}{1+\mu}\frac{\alpha(t,t_2)}{\eps} 
\Bigbrak{1-\BigOrder{\frac1{\mu\log(h/\sigma)}}}}
\end{equation}
under the condition
\begin{equation}
\label{sesc:4}
\Bigpar{\frac h\sigma}^{3+\mu} 
\Bigpar{1 + (1+\mu) \frac{\eps}{t_2^3}\log\frac h\sigma} 
\leqs L\mskip1.5mu \frac{t_2^4}{\sigma^2}.
\end{equation}
\end{prop}
\goodbreak
\begin{proof}
The proof follows along the lines of the one of
\cite[Proposition~4.7]{BG}, the main difference being the quadratic
behaviour $a_-t^2\leqs a(t)\leqs a_+t^2$ of $a$ in our case as opposed
to the linear one in~\cite{BG}. 

We start by defining a partition $t_2=u_0<\dots<u_K=t$ of $[t_2,t]$,
given by
\begin{equation}
\label{sesc:5}
\alpha(u_k,u_{k-1}) = (1+\mu)\frac\eps2 \log\frac{h^2}{\sigma^2},
\qquad \text{for } 
1\leqs k < K =
\biggintpartplus{\frac{2\alpha(t,t_2)}{(1+\mu)\eps\log(h^2/\sigma^2)}}.
\end{equation}
On each interval $[u_k,u_{k+1}]$, we consider a Gaussian approximation
$(x^{(k)}_t)_{t\in[u_k,u_{k+1}]}$ of $x_t$, defined by 
\begin{equation}
\label{sesc:6}
\6 x^{(k)}_t = \frac1\eps a(t)x^{(k)}_t \6 t 
+ \frac{\sigma}{\sqrt\eps} \6 W^{(k)}_t
\qquad x^{(k)}_{u_k}=x_{u_k},
\end{equation}
where $W^{(k)}_t = W_t - W_{u_k}$.
If $\abs{x_s}\sqrt{a(s)}\leqs h$ for all $s\in[u_k,u_{k+1}]$, then by
\eqref{s3} and \eqref{s4}, there is a constant $M>0$ such that 
\begin{equation}
\label{sesc:7}
\begin{split}
\abs{x_s-x^{(k)}_s} 
&\leqs \frac1\eps \int_{u_k}^s \abs{g_0(x_u,u)x_u} \e^{\alpha(s,u)/\eps} \6
u \\
&\leqs M \frac{h^3}{a(u_k)^{3/2}} \frac1{a(u_k)}
\e^{\alpha(u_{k+1},u_k)/\eps} \leqs \frac{h}{\sqrt{a(s)}} 
\end{split}
\end{equation}
for  all $s\in[u_k,u_{k+1}]$, provided the condition
\begin{equation}
\label{sesc:8}
h^2 \leqs \frac{a_-^2}M \sqrt{\frac{a(u_k)}{a(u_{k+1})}}
\e^{-\alpha(u_{k+1},u_k)/\eps} t_2^4
\end{equation}
holds for all $k$. Now,
\begin{equation}
\label{sesc:9}
\frac{a(u_{k+1})}{a(u_k)}
\leqs 1 + \frac{ca_+}{a_-^3}\frac{\alpha(u_{k+1},u_k)}{t_2^3} 
\Bigpar{1+\frac{a_+}{a_-^2}\frac{\alpha(u_{k+1},u_k)}{t_2^3}},
\end{equation}
where $c$ is a constant satisfying $0\leqs a'(t) \leqs ct$ for all
$t\in[0,T]$. This shows that there exists a constant $L>0$ such that
the condition~\eqref{sesc:8} 
is satisfied whenever 
\begin{equation}
\label{sesc:10}
\Bigpar{\frac h\sigma}^{3+\mu} 
\Bigpar{1 + \frac{\alpha(u_{k+1},u_k)}{t_2^3}}
\leqs L\mskip1.5mu \frac{t_2^4}{\sigma^2}.
\end{equation}
This condition is equivalent to~\eqref{sesc:4}.

Assume $\abs{x_{u_k}} \sqrt{a(u_k)} \leqs h$ for the moment. Then,
\begin{align}
\nonumber
\Bigprobin{u_k,x_{u_k}}{\sup_{u_k\leqs s\leqs u_{k+1}}
\abs{x_s}\sqrt{a(s)} \leqs h}
&\leqs \Bigprobin{u_k,x_{u_k}}{\abs{x^{(k)}_{u_{k+1}}}\sqrt{a(u_{k+1})}
\leqs
2h} \\
\label{sesc:12}
&\leqs \frac{4h}{\sqrt{2\pi v^{(k)}_{u_{k+1}} a(u_{k+1})}},
\end{align}
where $v^{(k)}_{u_{k+1}}$ denotes the variance of
$x^{(k)}_{u_{k+1}}$. By partial integration, we find
\begin{equation}
\label{sesc:13}
v^{(k)}_{u_{k+1}} = \frac{\sigma^2}{\eps} \int_{u_k}^{u_{k+1}}
\e^{2\alpha(u_{k+1},s)/\eps} \6 s
\geqs \frac{\sigma^2} {a(u_{k+1})}
\Bigbrak{\e^{2\alpha(u_{k+1},u_k)/\eps}-1}.
\end{equation}
Now, the Markov property yields
\begin{equation}
\label{sesc:11}
\bigprobin{t_2,x_2}{\tau_\cS \geqs t}
= \Bigprobin{t_2,x_2}{\sup_{t_2\leqs s\leqs t}
\abs{x_s}\sqrt{a(s)} \leqs h} 
\leqs \prod_{k=0}^{K-1} \biggpar{\frac4{\sqrt{2\pi}}
\frac{h}{\sqrt{v^{(k)}_{u_{k+1}} a(u_{k+1})}}\wedge 1},
\end{equation}
and the bound~\eqref{sesc:3} follows by a straightforward calculation. 
\end{proof}

The preceding proposition shows that a path starting in $\cS$ is
likely to leave $\cS$ after a short time. We want to show that such a
path (or any path starting in $\cD(\kappa)\setminus\cS$) is also
likely to leave $\cD(\kappa)$. For this purpose, we will again compare $x_t$
to a Gaussian approximation, given by
\begin{equation}
\label{sesc:14}
\6x^0_t = \frac1\eps a_0(t) x^0_t \6t + \frac\sigma{\sqrt\eps}\6W_t,
\end{equation}
where $a_0(t)=\kappa a(t)$, so that $f(x,t)/x \geqs a_0(t)$ in
$\cD(\kappa)$. Assume that $x_{t_2}>0$. Then $x_s\geqs x^0_s$ holds as
long as $x_s$ neither leaves $\cD(\kappa)$ nor crosses the $t$-axis,
cf.~\cite[Lemma~4.8]{BG}. Therefore we can proceed as follows. Once a
path is in $\cD(\kappa)\setminus\cS$, there are two
possibilities. Either, $x^0_s$ does not return to zero, or it does. If
$x^0_t$ does not return to zero, then it is likely to leave $\cD(\kappa)$
via
the upper boundary and so is $x_t$. So we are left with the case of
$x^0_t$ returning to zero. This event has a small but not negligible
probability. Note that if $x^0_t$ returns to zero, then $x_t$ is still
non-negative. If $x_t$ has nevertheless left $\cD(\kappa)$, we are done. If
not, $x_t$ is either in $\cS$ or in $\cD(\kappa)\setminus\cS$. Since
we may assume that, after a short time, $x_t$ is in
$\cD(\kappa)\setminus\cS$ again, we can repeat the above argument.  

Making the above-said precise, we obtain an integral equation for an
upper bound on the probability that $x_u$ does not leave $\cD(\kappa)$
up to time $t$, which will be solved by iterations. We will cite the
integral equation from~\cite{BG}, as the general arguments leading to
it do not require adaptation. Let us first introduce the necessary
notations. 
We choose $h=K\sigma$ for some (possibly large) constant $K>0$. For
$\kappa\in(0,1)$, we choose $\mu>0$ in such a way that
\begin{equation}
\label{sesc:14a}
\frac12 \frac{\mu}{1+\mu}\Bigbrak{1-\BigOrder{\frac{1}{\mu\log K}}}
\leqs \kappa \leqs
\frac{\mu}{1+\mu}\Bigbrak{1-\BigOrder{\frac{1}{\mu\log K}}}
\end{equation}
for all large enough $K$. Note that choosing $\kappa$ too close to $1$
requires large $\mu$ and is thus not desirable. Since we want to apply
Proposition~\ref{p_sescape} on the first exit time from $\cS$ with
$t_2\geqs c_2\sqrt\sigma$ for a suitably chosen $c_2$,
Condition~\eqref{sesc:4} must be satisfied. Therefore, we choose
$c_2=c_2(K)$ large enough for
\begin{equation}
\label{sesc:14b}
K^{3+\mu}\Bigpar{1+\frac{1+\mu}{c_2^3}\mskip1.5mu \log K}
\leqs Lc_2^4
\end{equation}
to hold. Now, set
\begin{equation}
\label{sesc:15}
g(t,s) = \frac{\e^{-\kappa\alpha(t,s)/\eps}}
{\sqrt{1-\e^{-2\kappa\alpha(t,s)/\eps}}},
\end{equation}
and 
\begin{equation}
\label{sesc:16}
C = \max\Bigset{\frac{\tx(t)\sqrt{\kappa a(t)}}{\sqrt\pi \sigma},1} 
\qquad\text{and}\qquad
c = \sqrt{\pi\kappa} \Bigpar{\frac h\sigma}^{1+\mu} 
\e^{-\kappa h^2/\sigma^2}.
\end{equation}
For $t_2\leqs s\leqs t\leqs T$, let $Q^{(0)}_t(s)\equiv 1$, and define
$Q^{(n)}_t(s)$ for $n\geqs1$ by 
\begin{equation}
\label{sesc:17}
Q^{(n)}_t(s) = C g(t,s) + c\e^{-\kappa\alpha(t,s)/\eps} 
+ c\int_s^t Q^{(n-1)}_t(u) \frac{a(u)}\eps \e^{-\kappa\alpha(u,s)/\eps}
\bigbrak{1+g(u,s)} \6 u.
\end{equation}
Then, \cite[(4.95) and (4.107)]{BG} show that for any $n$ and $s\geqs
t_2 \geqs c_2\sqrt\sigma$,
\begin{equation}
\label{sesc:18}
\sup_{x\colon(x,s)\in\cS} \bigprobin{s,x}{\tau_{\cD(\kappa)}\geqs t} 
\leqs 2\Bigpar{\frac{h}{\sigma}}^\mu \e^{-\kappa\alpha(t,s)/\eps} 
+ \kappa\Bigpar{\frac{h}{\sigma}}^\mu \int_s^t Q^{(n)}_t(u)
\frac{a(u)}{\eps} \e^{-\kappa\alpha(u,s)/\eps} \6 u
\end{equation}
and
\begin{equation}
\label{sesc:19}
\sup_{x\colon(x,s)\in\cD(\kappa)\setminus\cS}
\bigprobin{s,x}{\tau_{\cD(\kappa)}\geqs t}  
\leqs Q^{(n)}_t(s).
\end{equation}
Next we estimate $Q^{(n)}_t$ by showing that 
\begin{equation}
\label{sesc:20}
Q^{(n)}_t(s) \leqs  C g(t,s) + a_n \e^{-\kappa\alpha(t,s)/2\eps} + b_n
\qquad\forall\mskip1.5mu n
\end{equation}
holds with $a_1=c$, $b_1=3c/\kappa$ in the case $n=1$, and with
\begin{align}
\label{sesc:21a}
a_n &= c\Bigpar{1+\frac{4C}{\kappa}} \sum_{j=0}^{n-2}
\Bigpar{\frac{6c}{\kappa}}^j + c\Bigpar{\frac{6c}{\kappa}}^{n-1} 
\leqs \Bigpar{1+\frac{4C}{\kappa}} \frac {c}{1-6c/\kappa} 
\leqs 2c \Bigpar{1+\frac{4C}{\kappa}}
\\
\label{sesc:21b}
b_n &= \biggpar{\frac {3c}{\kappa}}^n
\end{align}
\goodbreak
\noindent
for $n>1$, provided $6c/\kappa \leqs 1/2$. Note that the latter
imposes a condition on $K=h/\sigma$. To obtain the bound~\eqref{sesc:20}, we
proceed as in~\cite{BG}, the only difference lying in the term $a_n
\e^{-\kappa\alpha(t,s)/2\eps}$, where we sacrifice a factor of $2$ in
the exponent in order to gain a smaller coefficient $a_n$. Our choice
of $a_n$ yields a less restrictive condition on $h/\sigma$, namely we
only need
\begin{equation}
\label{sesc22}
\frac{12\sqrt\pi}{\sqrt\kappa} \Bigpar{\frac{h}{\sigma}}^{1+\mu}
\e^{-\kappa h^2/\sigma^2} \leqs 1,
\end{equation}
which is satisfied whenever $K=h/\sigma$ is large enough.

Now, \eqref{sesc:21a} and \eqref{sesc:21b} imply that for $K$
and $c_2(K)$ large enough,
\begin{equation}
\label{sesc23}
\sup_{x\colon(x,t_2)\in\cD(\kappa)}
\bigprobin{t_2,x}{\tau_{\cD(\kappa)}\geqs t}  
\leqs C_0  \Bigpar{\frac{t}{\sqrt\sigma}}^2 
\frac{\e^{-\kappa\alpha(t,t_2)/2\eps}}
{\sqrt{1-\e^{-\kappa\alpha(t,t_2)/\eps}}}
\qquad\text{for all $t\geqs t_2\geqs c_2\sqrt\sigma$}
\end{equation}
with some constant $C_0$. This completes our outline of the proof of
Theorem~\ref{thm_sescape}. 


\section{Asymmetric case}
\label{sec_asym}

We consider in this section the nonlinear SDE
\begin{equation}
\label{a1}
\6 x_\t = \frac{1}{\eps} f(x_\t,\t) \6\t +
\frac{\sigma}{\sqrt{\eps}} \6 W_\t,
\end{equation}
where $f$ satisfies the hypotheses given at the beginning of Subsection
\ref{ssec_rasym}. By rescaling $x$, we can arrange for $\sdpar{f}{xx}(\xc,0)
= -2$, so that Taylor's formula allows us to write 
\begin{equation}
\label{a2}
\sdpar fx(\xc+\tx,t) = \sdpar fx(\xc,t) + \tx\bigbrak{-2+r_1(\tx,t)} 
\end{equation}
where $r_1\in\cC^1$ and $r_1(0,0)=0$. Since $\sdpar fx(\xc,t)=\Order{t^2}$
by assumption, $\sdpar fx(x,t)$ vanishes on a curve
$\bx(t)=\xc+\Order{t^2}$. We further obtain that 
\begin{equation}
\label{a3}
\begin{split}
f(\bx(t)+z,t) &= f(\bx(t),t) + z^2\bigbrak{-1+r_0(z,t)} \\
f(\bx(t),t) &= f(\xc,t) + \Order{t^4} = a_0 + a_1 t^2 + \Order{t^3},
\end{split}
\end{equation}
where $r_0\in\cC^1$ and $r_0(0,0)=0$. Thus $f(x,t)$ vanishes on two curves
$x^\star_+(t)$ and $x^\star_0(t)$, which behave near $t=0$ like
$\xc\pm\sqrt{a_0+a_1 t^2}[1+\Order{\sqrt{a_0+a_1 t^2}\mskip1.5mu}]$, as
indicated in~\eqref{ares4c}. The behaviour of the linearization follows
from~\eqref{a2}. 


\subsection{Deterministic case}
\label{ssec_asymdet}

The proof of Theorem~\ref{thm_adet} follows closely the proof of Theorem
\ref{thm_sdet}, with some minor differences we comment on here. The
dynamics of $y_t = x_t-x^\star_+(t)$ is still governed by an equation of
the form 
\begin{equation}
\label{adet2}
\eps\dtot yt = a_+^\star(t)y + b_+^\star(y,t) - \eps\dtot{x_+^\star}t,
\end{equation}
but now Taylor's formula yields the relations
\begin{align}
\label{adet3a}
a_+^\star(t) &\asymp - (\sqrt{a_0}\vee \abs{t}) \\
\label{adet3b}
b_+^\star(y,t) &= - y^2 \bigbrak{1+\Order{\sqrt{a_0}} + \Order{t} +
\Order{y}}, 
\end{align}
while relation~\eqref{sdet3c} holds for the derivative of $x_+^\star$, with
$t^\star$ replaced by $t^\star_+$.  Lemma~\ref{lem_sdet1} becomes

\begin{lemma}
\label{lem_adet1}
Let $\ta(t)$ be a continuous function satisfying $\ta(t)\asymp-(\beta\vee
\abs{t})$ for $\abs{t}\leqs T$, where $\beta=\beta(\eps)\geqs 0$. Let
$\chi_0\asymp 1$, and define $\talpha(t,s)=\int_s^t \ta(u)\6u$. Then 
\begin{equation}
\label{adet4}
\chi_0 \e^{\talpha(t,-T)/\eps} + \frac1\eps\int_{-T}^t
\e^{\talpha(t,s)/\eps}\6s \asymp 
\begin{cases}
\vrule height 14pt depth 17pt width 0pt
\dfrac1{\beta\vee\sqrt\eps} &\text{for $\abs{t}\leqs
\beta\vee\sqrt\eps$} \\
\vrule height 17pt depth 14pt width 0pt
\dfrac1{\abs{t}} &\text{for $\beta\vee\sqrt\eps\leqs\abs{t}\leqs T$.}
\end{cases}
\end{equation}
\end{lemma}

Proposition~\ref{prop_sdet1} carries over with some obvious adjustments, and
shows the existence of a constant $c_0$ such that
\begin{equation}
\label{adet4a}
y_t \asymp \frac\eps{\abs t} 
\qquad
\text{for $-T\leqs t\leqs t_0 = -c_0(\sqrt{a_0}\vee\sqrt\eps\mskip1.5mu)$.}
\end{equation}
In particular, $y_{t_0}\asymp (\eps/\sqrt{a_0}\mskip1.5mu)\wedge\sqrt\eps$.
An adaptation of Proposition~\ref{prop_sdet2} yields the existence of a
constant $\gamma_0>0$ such that, for $a_0\geqs\gamma_0\eps$, 
\begin{equation}
\label{adet4b}
y_t = C_1(t)(t^\star_+-t) + C_2(t) 
\qquad\text{with}\qquad 
C_1(t)\asymp\frac{\eps}{a_0},
\quad 
C_2(t)\asymp\frac{\eps^2}{a_0^{3/2}}
\end{equation}
for all $\abs{t}\leqs \abs{t_0}$. This shows in particular that $y_t$
vanishes at a time $\tilde t$ satisfying $\tilde t - t^\star_+ \asymp
\eps/\sqrt{a_0}$. 
Proposition~\ref{prop_sdet3} is replaced by

\begin{prop}
\label{prop_adet3}
Assume that $a_0<\gamma_0\eps$. Then, for any fixed $t_1\asymp
\sqrt\eps$, 
\begin{equation}
\label{adet5}
x_t - \xc \asymp \sqrt\eps 
\qquad\text{for $t_0\leqs t\leqs t_1$,}
\end{equation}
and $x_t$ crosses $x_+^\star(t)$ at a time $\tilde t$ satisfying
$\tilde t \asymp \sqrt\eps$. 
\end{prop}
\goodbreak
\begin{proof}
Let $\tx_t = x_t - \xc$. 
We first observe that by Taylor's formula,
\begin{equation}
\label{padet3:1}
f(\xc+\tx,t) = f(\xc,t) + \tx\sdpar fx(\xc,t) + \tx^2
\bigbrak{-1+\Order{\tx}+\Order{t}}.
\end{equation}
This shows that 
\begin{equation}
\label{padet3:2}
\eps\dtot{\tx}t = a_0 + a_1 t^2 - \tx^2 
+ \Order{t^3} + \Order{t^2\tx} + \Order{t\tx^2} + \Order{\tx^3}.
\end{equation}
Thus, with the rescaling
\begin{equation}
\label{padet3:3}
\tx = a_1^{1/4}\sqrt\eps\, z, \qquad
t = a_1^{-1/4}\sqrt\eps\, s,
\end{equation}
we obtain that $z_t$ obeys a perturbation of order $\sqrt\eps$ of the
Riccati equation
\begin{equation}
\label{padet3:4}
\dtot zs = \tilde a_0 + s^2 - z^2, 
\qquad
\text{with $\tilde a_0 = \frac1{\sqrt{a_1}}\frac{a_0}\eps <
\frac{\gamma_0}{\sqrt{a_1}}$.}
\end{equation}
One easily shows that the solution satisfies $z_s\asymp 1$ for $s$ of order
1, and this property carries over to the perturbed equation with the help
of Gronwall's inequality. Finally, since $\tx_t\asymp\sqrt\eps$ and
$x^\star_+(t)-\xc\asymp \sqrt{a_0}\vee\abs{t}$, these curves necessarily
cross at a time $\tilde t\asymp\sqrt\eps$. 
\end{proof}
\goodbreak

The assertion on the existence of a particular solution $\xhatdet$
tracking the unstable equilibrium branch $x^\star_0(t)$ follows from the
observation that $z_s = x_{-s}$ satisfies the equation
\begin{equation}
\label{adet6}
\eps\dtot {z_s}s = -f(z_s,-s).
\end{equation}
This system admits $z^\star_0(s)=x^\star_0(-s)$ as a stable equilibrium
branch. Thus the same arguments as above can be used to show the existence
of a solution $z_s$ tracking $z^\star_0(s)$, with similar properties. 
Proposition~\ref{prop_sdet5} admits the following counterpart:

\begin{prop}
\label{prop_adet5}
For all $t\in[-T,T]$ and all $a_0=\orderone{\eps}$, 
\begin{align}
\label{adet7a}
\ba(t) &\defby \sdpar fx(\xdet_t,t) \asymp - (\abs{t}\vee
\sqrt{a_0}\vee\sqrt\eps\mskip1.5mu) \\
\label{adet7b}
\ha(t) &\defby \sdpar fx(\xhatdet_t,t) \asymp \abs{t}\vee
\sqrt{a_0}\vee\sqrt\eps.
\end{align}
\end{prop}
\goodbreak
\begin{proof}
The proof is a direct consequence of \eqref{a2} and the properties of
$\xdet_t$, and thus much simpler than the proof of
Proposition~\ref{prop_sdet5}.
\end{proof}

Finally, with Lemma~\ref{lem_adet1}, we immediately obtain
\begin{equation}
\label{adet9}
\z(t) \defby \frac 1{2\abs{\ba(-T)}} \e^{2\balpha(t,-T)/\eps} +
\frac1{\eps} \int_{-T}^t \e^{2\balpha(t,s)/\eps} \6s
 \asymp \frac1{\abs{t}\vee\sqrt{a_0}\vee\sqrt\eps}. 
\end{equation}


\subsection{The transition regime}
\label{ssec_asymtrans}

We consider now the regime of $\sigma$ sufficiently large to allow for
transitions from the potential well at $x^\star_+$ to the potential well at
$x^\star_-$, by passing over the saddle at $x^\star_0$. Here $\xdet_t$  and
$\xhatdet_t$ denote solutions of the deterministic equation 
\begin{equation}
\label{atrans1}
\eps\dtot xt = f(x,t)
\end{equation}
tracking, respectively, the stable equilibrium branch $x^\star_+(t)$ and
the unstable equilibrium branch $x^\star_0(t)$, while $x_t$ denotes a
general solution of the SDE~\eqref{a1}. Our aim is to establish an upper
bound for the probability {\em not}\/ to reach a level $\delta_0$ between
$x^\star_0(t)$ and $x^\star_-(t)$, situated at a distance of order $1$ from
both equilibria. \cite[Theorem 2.3]{BG} shows that if $x_t$ reaches
$\delta_0$, and $\delta_0$ is close enough to $x^\star_-(t)$ (but it may
still be at a distance of order $1$), then it is likely to reach a \nbh\ of
$x^\star_-(t)$ as well. 

Let $\delta_0 < \delta_1 < \xc < \delta_2$ be the constants satisfying 
\eqref{ares17}, that is, 
\begin{equation}
\label{atrans2}
\begin{split}
f(x,t) \asymp -1
&\qquad\qquad
\text{for $\delta_0\leqs x\leqs \delta_1$ and $\abs{t}\leqs T$} \\
\sdpar f{xx}(x,t) \leqs 0
&\qquad\qquad
\text{for $\delta_1\leqs x\leqs \delta_2$ and $\abs{t}\leqs T$.}
\end{split}
\end{equation} 
The basic ingredient of our estimate is the following analogue of Lemma
\ref{lem_strans}:

\begin{lemma}
\label{lem_atrans}
Fix some initial time $t_0\in[-T,T]$. 
We consider the following processes on $[t_0,T]$:
\begin{itemiz}
\item	the solution $\xdet_t$ of the deterministic differential
equation~\eqref{atrans1} with an initial condition $\xdet_{t_0} \in
[\delta_1,\delta_2]$, such that $\xdet_t\geqs\delta_1$ for
all $t\in[t_0,T]$;
\item	the solution $x_t$ of the SDE~\eqref{a1} with an initial condition
$x_{t_0} \in [\xdet_{t_0},\delta_2]$;
\item	the difference $y_t=x_t-\xdet_t$, which satisfies $y_{t_0} =
x_{t_0}-\xdet_{t_0} \geqs 0$;
\item	 the solution $y^0_t$of the linear SDE
\begin{equation}
\label{atrans3}
\6y^0_t = \frac1\eps \tilde a(t)  y^0_t \6t 
+ \frac{\sigma}{\sqrt{\eps}} \6W_t, 
\qquad
\text{where $\tilde a(t) = \sdpar{f}{x}(\xdet_t,t)$}
\end{equation}
with initial condition $y^0_{t_0} \in[y_{t_0}, \delta_2-\xdet_{t_0}]$.
\end{itemiz}
If\/ $\delta_1 \leqs y^0_s +\xdet_s\leqs \delta_2$ for all $s\in[t_0,t]$,
then
$y_s\leqs y^0_s$ for those $s$. Similarly, if\/ $\delta_1 \leqs x_s \leqs
\delta_2$ for all $s\in[t_0,t]$, then $y^0_s\geqs y_s$ for those
$s$. The result remains true when $t$ is replaced by a stopping time.
\end{lemma}

We will proceed as follows. Let $\z(t)$ be the function defined in
\eqref{adet9}, and let $h$ be such that $\xdet_s+h\sqrt{\z(s)} < \delta_2$
for all $s\in[t_0,t]$. Given $x_0\in(\delta_1,\delta_2)$ and times
$t_0<t_1<t$ in $[-T,T]$, we consider the solution $x_t$ of the
SDE~\eqref{a1} with initial condition $x_{t_0}=x_0$. We introduce the
stopping time
\begin{equation}
\label{atrans4}
\tau = \inf\bigsetsuch{s\in[t_0,t_1]}{x_s\leqs\delta_1} 
\in[t_0,t_1] \cup\set{\infty}.
\end{equation}
We can thus write
\begin{align}
\nonumber
\bigprobin{t_0,x_0}{x_s>\delta_0 \;\; \forall s\in[t_0,t]}
\leqs{} &  \bigprobin{t_0,x_0}{x_s>\delta_1 \;\; \forall s\in[t_0,t_1]} \\
\label{atrans5}
&{}+ \Bigexpecin{t_0,x_0}{\indexfct{\tau\leqs t_1}
\bigprobin{\tau,\delta_1}{x_s>\delta_0 \;\; \forall s\in[\tau,t]}}. 
\end{align}
The first term can be further estimated by
\begin{align}
\nonumber
\bigprobin{t_0,x_0}{x_s>\delta_1 \;\; \forall s\in[t_0,t_1]}
\leqs{} &  \Bigprobin{t_0,x_0}{\sup_{t_0\leqs s\leqs t_1}
\frac{x_s-\xdet_s}{\sqrt{\z(s)}} > h} 
 \\
\label{atrans6}
&{}+
\bigprobin{t_0,x_0}{\delta_1<x_s\leqs\xdet_s+h\sqrt{\z(s)} 
\;\; \forall s\in[t_0,t_1]}.
\end{align}
The two summands in~\eqref{atrans6} can be estimated in a similar way as in
the symmetric case. The first one is dealt with in the following result.

\begin{prop}
\label{prop_atrans1}
Assume $\delta_1\leqs x_0\leqs \xdet_{t_0}+\frac12 h\sqrt{\z(t_0)}$. Then 
\begin{equation}
\label{atrans7}
\Bigprobin{t_0,x_0}{\sup_{t_0\leqs s\leqs t_1}
\frac{x_s-\xdet_s}{\sqrt{\z(s)}} > h} 
\leqs \frac52 \Bigpar{\frac{\abs{\balpha(t_1,t_0)}}\eps + 1} 
\e^{-\kappa h^2/\sigma^2},
\end{equation}
where $\kappa$ is a positive constant and $\balpha(t_1,t_0) =
\int_{t_0}^{t_1} \ba(s)\6s$. 
\end{prop}

The proof is almost the same as the proof of
Proposition~\ref{prop_strans1}.  Instead of \eqref{strans1:6}, we may
define $\tau^0$ and $\tau^+$ as the first times when $\smash{y^{(k)}_s}$
either reaches $\delta_1 - \xdet_s$ or the upper boundary $h\sqrt{\z(s)}$.
Then Lemma~\ref{lem_atrans} implies $y_s\leqs \smash{y^{(k)}_s}$ for
$s\in[u_k,u_{k+1}\wedge\tau^0\wedge\tau^+]$. However, when estimating the
probability that $\tau^0 < u_{k+1}$ as in \eqref{strans1:9}, it is
sufficient to use the fact that $\tau^0$ is larger than the first time
$\smash{y^{(k)}_s}$ reaches $0$.  Finally, \eqref{strans1:12} still holds
with the present definitions of $\z$ and $\ba$, because of~\eqref{adet9}
and Proposition~\ref{prop_adet5}. 

\goodbreak
Let us now examine the second term in~\eqref{atrans6}.

\begin{prop}
\label{prop_atrans2}
There exist constants $c_1>0$ and $\bar\kappa>0$ such that, if
$c_1^{3/2}\sigma \geqs a_0^{3/4}\vee\eps^{3/4}$ and $h>2\sigma$, then
\begin{multline}
\label{atrans8}
\qquad\qquad
\bigprobin{-c_1\sigma^{2/3},x_0}{\delta_1<x_s\leqs\xdet_s+h\sqrt{\z(s)}
 \;\; \forall s\in[-c_1\sigma^{2/3},t_1]} \\
\leqs \frac32 \exp\biggset{-\bar\kappa \mskip3mu \frac1{\log(h/\sigma)
\vee\abs{\log\sigma}}
\frac{\halpha(t_1,-c_1\sigma^{2/3})}\eps}
\qquad\qquad\qquad
\end{multline}
holds with $\halpha(t,s) = \int_s^t \ha(u)\6u$, for 
$-c_1\sigma^{2/3} \leqs t_1 \leqs c_1\sigma^{2/3}$ and all initial
conditions $x_0$ satisfying $\delta_1\leqs x_0\leqs
\xdet_{-c_1\sigma^{2/3}}+h\sqrt{\z(-c_1\sigma^{2/3})}$.
\end{prop}
\goodbreak
\begin{proof}\hfill
\begin{enum}
\item	Let $\xhatdet_t$ be the deterministic solution tracking the saddle
at $x^\star_0(t)$ and set $x_t=\xhatdet_t+z_t$. Then
\begin{equation}
\label{atrans2:1}
\6z_t = \frac1\eps \bigbrak{\ha(t)z_t + \hb(z_t,t)}\6t + 
\frac\sigma{\sqrt\eps}\6W_t,
\end{equation}
where \eqref{adet7b}, \eqref{adet9} and \eqref{atrans2} imply
\begin{equation}
\label{atrans2:2a}
\ha(t) \asymp \abs{t} \vee \sqrt{a_0} \vee \sqrt\eps 
\asymp \frac1{\z(t)} 
\end{equation}
and
\begin{equation}
\label{atrans2:2b}
\hb(z,t) \leqs 0 
\qquad\qquad\text{for $\xhatdet_t+z\in[\delta_1,\delta_2]$.}
\end{equation}
Let $\varrho = \varrho(h/\sigma) \geqs 1$ and define a partition
$-c_1\sigma^{2/3}=u_0 < \dots < u_K=t_1$ of $[-c_1\sigma^{2/3},t]$ by 
\begin{equation}
\label{atrans2:3}
\halpha(u_k,u_{k-1}) = \varrho\eps
\qquad \text{for\ }
1\leqs k<K = 
\biggintpartplus{\frac{\halpha(t_1,-c_1\sigma^{2/3})}{\varrho\eps}}.
\end{equation}
Let 
\begin{equation}
\label{atrans2:4}
Q_k = \sup_{\delta_1<\xhatdet_{u_k}+z_k\leqs \xdet_{u_k}+h\sqrt{\z(u_k)}} 
\bigprobin{u_k,z_k}{\delta_1<\xhatdet_s+z_s\leqs\xdet_s+h\sqrt{\z(s)}
 \;\; \forall s\in[u_k,u_{k+1}]}.
\end{equation}
Then we have, as in~\eqref{strans2:3}, 
\begin{equation}
\label{atrans2:5}
\bigprobin{-c_1\sigma^{2/3},x_0}{\delta_1<x_s\leqs\xdet_s+h\sqrt{\z(s)}
 \;\; \forall s\in[-c_1\sigma^{2/3},t_1]}
 \leqs \prod_{k=0}^{K-1} Q_k.
\end{equation}
The result will thus be proved if we manage to choose $\varrho$ in such a
way that $Q_k$ is bounded away from $1$ for $k=0,\dots,K-2$. 

\item	We will estimate the $Q_k$ in a similar way as in Proposition
\ref{prop_strans2}, but we shall distinguish three cases instead of two.
These cases correspond to $x_s$ reaching the levels $\xdet_s$, $\xhatdet_s$
and $\delta_1$. We introduce a subdivision $u_k < \smash{\tilde u_{k,1} < 
\tilde u_{k,2}} < u_{k+1}$ defined by 
\begin{equation}
\label{atrans2:6}
\halpha(\tilde u_{k,1},u_k) = \frac13 \varrho \eps, 
\qquad\qquad
\halpha(\tilde u_{k,2},u_k) = \frac23 \varrho \eps, 
\end{equation}
and stopping times
\begin{equation}
\label{atrans2:7}
\begin{split}
\tau_{k,1} &= \inf\bigsetsuch{s\in[u_k,\tilde u_{k,1}]}
{z_s\leqs\xdet_s-\xhatdet_s} \\
\tau_{k,2} &= \inf\bigsetsuch{s\in[u_k,\tilde u_{k,2}]}
{z_s\leqs 0}.
\end{split}
\end{equation}
Then we can write, similarly as in~\eqref{strans2:11}, 
\begin{align}
\label{atrans2:8}
&
\bigprobin{u_k,z_k}{\delta_1<\xhatdet_s+z_s\leqs\xdet_s+h\sqrt{\z(s)}
 \;\; \forall s\in[u_k,u_{k+1}]}\\
\nonumber&
\vrule height 13pt depth 7pt width 0pt
\leqs \bigprobin{u_k,z_k}{\xdet_s<\xhatdet_s+z_s\leqs\xdet_s+h\sqrt{\z(s)} 
\;\; \forall s\in[u_k,\tilde u_{k,1}]} \\
&\phantom{\leqs{}}
{}+ \Bigexpecin{u_k,z_k}{\indexfct{\tau_{k,1}<\tilde u_{k,1}} 
\bigprobin{\tau_{k,1},z_{\tau_{k,1}}}
{\delta_1<\xhatdet_s+z_s\leqs\xdet_s+h\sqrt{\z(s)}
 \;\; \forall s\in[\tau_{k,1},u_{k+1}]}}
\nonumber
\end{align}
The first term can be bounded by comparing with the solution of the
SDE~\eqref{a1} linearized around $\xdet_t$, with the help of
Lemma~\ref{lem_atrans}.  As in~\eqref{strans2:12}, we obtain the upper
bound
\begin{equation}
\label{atrans2:9}
\frac2{\sqrt{\pi}} \sqrt{\frac1{1-\e^{-2\varrho/3L}}} 
\sup_{u_k\leqs u\leqs\tilde u_{k,1}} 
\sqrt{\abs{\ba(u)}\z(u_k)} \,\frac h\sigma \e^{-\varrho/3L},
\end{equation}
where $L>0$ is a constant such that $\ha(u) \leqs L\abs{\ba(u)}$ for all
$u$. Now if $\tau_{k,1}<\tilde u_{k,1}$, we also have 
\begin{align}
\nonumber&
\mskip-12mu
\bigprobin{\tau_{k,1},z_{\tau_{k,1}}}
{\delta_1<\xhatdet_s+z_s\leqs\xdet_s+h\sqrt{\z(s)}
 \;\; \forall s\in[\tau_{k,1},u_{k+1}]}\\
\label{atrans2:10}
\vrule height 13pt depth 7pt width 0pt
&
\leqs \bigprobin{\tau_{k,1},z_{\tau_{k,1}}}
{0<z_s\leqs\xdet_s-\xhatdet_s+h\sqrt{\z(s)} 
\;\; \forall s\in[\tau_{k,1},\tilde u_{k,2}]} \\
\nonumber 
&\phantom{\leqs{}}
{}+ \Bigexpecin{\tau_{k,1},z_{\tau_{k,1}}}
{\indexfct{\tau_{k,2}<\tilde u_{k,2}} \\
&\qquad\qquad\qquad\quad
\bigprobin{\tau_{k,2},z_{\tau_{k,2}}}
{\delta_1<\xhatdet_s+z_s\leqs\xdet_s+h\sqrt{\z(s)}
 \;\; \forall s\in[\tau_{k,2},u_{k+1}]}}.
\nonumber
\end{align}
Comparing with the solution of the SDE~\eqref{a1} linearized around
$\xhatdet_t$, the first term can be bounded, as in~\eqref{strans2:6}, by 
\begin{equation}
\label{atrans2:11}
\frac2{\sqrt{\pi}} \sqrt{\frac1{1-\e^{-2\varrho/3}}} 
\sup_{u_k\leqs u\leqs\tilde u_{k,2}} 
\frac{\sqrt{\ha(u)}}\sigma 
\sup_{u_k\leqs u\leqs\tilde u_{k,2}}(\xdet_u-\xhatdet_u).
\end{equation}
This estimate shows that a path starting on $\xdet$ at time $\tau_{k,1}$ has
an appreciable probability to reach the saddle before time $\tilde u_{k,2}$.
Note, however, that we cannot obtain directly a similar estimate for the
probability to reach $\delta_1$ as well, which is why we restart the process
in $\xhatdet$.	

\item	In order to estimate the second summand in~\eqref{atrans2:10}, let
$z^{(k)}_s$ be the process starting in $0$ at time $\tau_{k,2}$ which
satisfies the linear equation
\begin{equation}
\label{atrans2:12}
\6z^{(k)}_s = \frac1\eps \ha(s) z^{(k)}_s \6s + \frac{\sigma}{\sqrt\eps}
\6W^{(k)}_s,
\end{equation}
with $W^{(k)}_s=W_s - W_{\tau_{k,2}}$. The variance $v^{(k)}_{u_{k+1}}$ of
$z^{(k)}_{u_{k+1}}$ satisfies, as in~\eqref{strans2:5}, 
\begin{equation}
\label{atrans2:13}
\e^{-2\halpha(u_{k+1},\tau_{k,2})/\eps} v^{(k)}_{u_{k+1}} 
\geqs \frac{1-\e^{-2\varrho/3}}2 \inf_{u_k\leqs u\leqs u_{k+1}}
\frac{\sigma^2}{\ha(u)}.
\end{equation}
Thus we obtain, using Lemma~\ref{lem_atrans}, 
\begin{align}
\nonumber&
\bigprobin{\tau_{k,2},z_{\tau_{k,2}}}
{\delta_1<\xhatdet_s+z_s\leqs\xdet_s+h\sqrt{\z(s)}
 \;\; \forall s\in[\tau_{k,2},u_{k+1}]}
\\
\nonumber&\qquad
\vrule height 13pt depth 7pt width 0pt
\leqs \bigprobin{\tau_{k,2},0}{z^{(k)}_{u_{k+1}}
\geqs-(\xhatdet_{u_{k+1}}-\delta_1)} 
\\
\nonumber&\qquad
{}= \frac1{\sqrt{2\pi}}
\int_{-(\xhatdet_{u_{k+1}}-\delta_1)(v^{(k)}_{u_{k+1}})^{-1/2}}^\infty
\e^{-z^2/2}\6z 
\\
\label{atrans2:14}&\qquad
{}\leqs \frac12 +  \frac1{\sqrt\pi}\sqrt{\frac1{1-\e^{-2\varrho/3}}}
\sup_{u_k\leqs u\leqs u_{k+1}} 
\frac{\sqrt{\ha(u)}}\sigma (\xhatdet_{u_{k+1}}-\delta_1) \e^{-\varrho/3}.
\end{align}
Here the introduction of the stopping time $\tau_{k,2}$ turns out to play a
crucial role. The above probability is indeed close to $1/2$ when 
$\varrho$ is larger than a constant times $\abs{\log\sigma}$, which shows
that once a path has reached the saddle, it also has about fifty percent
chance to reach the level $\delta_1$ in a time of order
$\eps\abs{\log\sigma}/\ha(u)$. 

\item	From \eqref{atrans2:9}, \eqref{atrans2:11} and \eqref{atrans2:14}
and the fact that $\varrho\geqs 1$, we obtain the existence of a constant
$C>0$ such that
\begin{align}
\label{atrans2:15}
Q_k \leqs \frac12 + \frac C\sigma \sup_{u_k\leqs u\leqs u_{k+1}}
\sqrt{\ha(u)} \Bigbrak{&\sqrt{\z(u_k)}\, h \e^{-\varrho/3L} \\
&{}+ \sup_{u_k\leqs u\leqs u_{k+1}} (\xdet_u-\xhatdet_u) 
+ (\xhatdet_{u_{k+1}}-\delta_1) \e^{-\varrho/3}}.
\nonumber
\end{align}
Since $\abs{u}\leqs c_1\sigma^{2/3}$, the properties of $\xdet$, $\xhatdet$,
$\z$ and $\ha$ imply the existence of another constant $C_1$ such that
\begin{equation}
\label{atrans2:16}
Q_k \leqs \frac12 + C_1 
\Bigbrak{\frac h\sigma \e^{-\varrho/3L} +
c_1^{3/2} + \sqrt{c_1}\, \frac{\e^{-\varrho/3}}{\sigma^{2/3}}}.
\end{equation}
Now choosing $c_1^{3/2}=1/(18 C_1)$ and 
$\varrho = 3L\log(18C_1h/\sigma) \vee 
3\log(18C_1\sqrt{c_1}\mskip1.5mu/\sigma^{2/3})
\vee 1$, we get $Q_k\leqs 2/3$ for $k=0,\dots,K-2$ and thus
\begin{equation}
\label{atrans2:17}
\prod_{k=0}^{K-1} Q_k \leqs \frac32\frac1{(3/2)^K} 
\leqs \frac32\exp\Bigset{-\Bigpar{\log\frac32}\mskip1.5mu 
\frac{\halpha(t_1,-c_1\sigma^{2/3})}{\varrho\eps}}, 
\end{equation}
and the result follows from our choice of $\varrho$. 
\qed
\end{enum}
\renewcommand{\qed}{}
\end{proof}

It remains to estimate the second term in~\eqref{atrans5}, describing the
probability not to reach $\delta_0$ when starting in $\delta_1$. This is
done by using the fact that, by assumption, the drift term is bounded away
from zero on the interval $[\delta_0,\delta_1]$. We will need to assume that
it can be bounded away from zero on a slightly larger interval, which is
possible by continuity of $f$.

\begin{prop}
\label{prop_atrans3}
Let $0<\rho\leqs \delta_1-\delta_0$ be a constant such that $f(x,t)\leqs
-f_0<0$ for $\delta_0\leqs x\leqs \delta_1+\rho$ and $\abs{t}\leqs T$. Then 
\begin{equation}
\label{atrans9}
\bigprobin{t_0,\delta_1}{x_s>\delta_0 \;\; \forall s\in[t_0,t_0+c\eps]} 
\leqs \e^{-\tilde\kappa/\sigma^2}
\end{equation}
holds for all $t_0\in[-T,T-c\eps]$, where $\tilde\kappa =
f_0\rho^2/4(\delta_1-\delta_0)$ and $c=2(\delta_1-\delta_0)/f_0$. 
\end{prop}
\goodbreak
\begin{proof}
Let $x^0_t$ be defined by 
\begin{equation}
\label{atrans3:1}
x^0_t = \delta_1 - \frac{f_0}\eps (t-t_0) + \frac\sigma{\sqrt\eps}
W_{t-t_0}, 
\qquad t\geqs t_0.
\end{equation}
By Gronwall's inequality, it is easy to see, as in Lemma~\ref{lem_atrans},
that if $\delta_0\leqs x^0_t\leqs \delta_1+\rho$ for all
$t\in[t_0,t_0+c\eps]$, then $x_t\leqs x^0_t$ for those $t$. We thus have 
\begin{multline}
\bigprobin{t_0,\delta_1}{x_s>\delta_0 \;\; \forall s\in[t_0,t_0+c\eps]} 
\leqs{}  \Bigprobin{t_0,\delta_1}{\sup_{t_0\leqs s\leqs t_0+c\eps} 
x^0_s + \frac{f_0}{\eps}(s-t_0) > \delta_1+\rho} \\
\label{atrans3:2}
{}+
\Bigprobin{t_0,\delta_1}{\delta_0<x^0_s<\delta_1+\rho-\frac{f_0}{\eps}(s-t_0)
 \;\; \forall s\in[t_0,t_0+c\eps]}. 
\end{multline}
Note, however, that for $s=t_0+c\eps$, 
\begin{equation}
\label{atrans3:3}
\delta_1+\rho-\frac{f_0}{\eps}(s-t_0) = \delta_1+\rho-2(\delta_1-\delta_0)
\leqs \delta_0,
\end{equation}
so that the second term in~\eqref{atrans3:2} is equal to zero. The first
term equals
\begin{equation}
\label{atrans3:4}
\Bigprobin{0,0}{\sup_{0\leqs s\leqs c\eps} \frac\sigma{\sqrt\eps} W_s >
\rho} 
\leqs \exp \Bigset{-\frac{\rho^2}{2c\sigma^2}}
\end{equation}
by Doob's submartingale inequality.
\end{proof}


\goodbreak

\bigskip\bigskip\noindent
{\small 
Nils Berglund \\ 
{\sc Department of Mathematics, ETH Z\"urich} \\ 
ETH Zentrum, 8092~Z\"urich, Switzerland \\
{\it E-mail address: }{\tt berglund@math.ethz.ch}

\bigskip\noindent
Barbara Gentz \\ 
{\sc Weierstra\ss\ Institute for Applied Analysis and Stochastics} \\
Mohrenstra{\ss}e~39, 10117~Berlin, Germany \\
{\it E-mail address: }{\tt gentz@wias-berlin.de}
}


\end{document}